\documentclass{amsart}

\usepackage[pdfauthor={Maxime Fortier Bourque},
    pdftitle={The holomorphic couch theorem},
    hidelinks
]{hyperref}

\usepackage{amssymb}
\usepackage{amsmath}
\usepackage{amsthm}
\usepackage{enumerate}
\usepackage{enumitem}
\usepackage{graphicx}
\usepackage{color}
\usepackage{footmisc}
\allowdisplaybreaks

\makeatletter
\newcommand{\myhypertarget}[2]{%
  \hypertarget{#1}{#2}%
    \protected@write\@mainaux{}{%
        \string\expandafter\string\gdef
          \string\csname\string\detokenize{#1}\string\endcsname{#2}%
    }%
  }
\newcommand{\myhyperlink}[1]{%
  \textup{\hyperlink{#1}{\csname #1\endcsname}}%
  }
\makeatother

\theoremstyle{plain}
\newtheorem{theorem}{Theorem}[section]
\newtheorem{proposition}[theorem]{Proposition}
\newtheorem{lemma}[theorem]{Lemma}
\newtheorem{corollary}[theorem]{Corollary}

\newtheorem{claim}[theorem]{Claim}
\newtheorem*{claim*}{Claim}

\theoremstyle{definition}
\newtheorem{definition}[theorem]{Definition}

\theoremstyle{remark}
\newtheorem*{remark}{Remark}

\newtheorem*{acknowledgements}{Acknowledgements}

\newcommand{\eps}{\varepsilon}
\newcommand{\vphi}{\varphi}
\newcommand{\wtilde}{\widetilde}
\newcommand{\what}{\widehat}

\DeclareMathOperator{\id}{id}
\DeclareMathOperator{\re}{Re}
\DeclareMathOperator{\im}{Im}
\DeclareMathOperator{\emb}{CEmb}
\DeclareMathOperator{\maps}{Map}
\DeclareMathOperator{\aut}{Aut}

\DeclareMathOperator{\Dil}{Dil}
\DeclareMathOperator{\res}{Res}

\DeclareMathOperator{\blob}{Blob}
\DeclareMathOperator{\ev}{ev}
\DeclareMathOperator{\dir}{Dir}
\DeclareMathOperator{\el}{EL}
\DeclareMathOperator{\lift}{lift}
\DeclareMathOperator{\hgt}{height}

\DeclareMathOperator{\Arg}{Arg}
\DeclareMathOperator{\deck}{Deck}

\DeclareMathOperator{\fell}{CL}
\DeclareMathOperator{\cat}{CAT}

\newcommand{\sect}{\rotatebox[origin=c]{90}{$\sphericalangle$}}

\newcommand{\HH}{\mathbb{H}}
\newcommand{\RR}{\mathbb{R}}
\newcommand{\DD}{\mathbb{D}}
\newcommand{\CC}{\mathbb{C}}
\newcommand{\ZZ}{\mathbb{Z}}
\newcommand{\NN}{\mathbb{N}}

\newcommand{\CHAT}{\widehat{\mathbb{C}}}
\newcommand{\m}{\mathfrak{m}}

\newcommand{\push}{\mathcal{P}}
\newcommand{\teich}{\mathcal{T}}
\newcommand{\qd}{\mathcal{Q}}

\newcommand{\dbar}{\overline \partial}

\begin{document}

\title{The holomorphic couch theorem}
\author{Maxime Fortier Bourque}
\thanks{Research partially supported by a Post\-gra\-dua\-te Scholarship from the Natural Sciences and Engineering Research Council of Canada and grant DMS 0905812 from the National Science Foundation.}
\address{Department of Mathematics, University of Toronto, 40 St. George Street, Toronto, ON, Canada M5S 2E4}
\email{mbourque@math.toronto.edu}

\begin{abstract}
We prove that if two conformal embeddings between Riemann surfaces with finite topology are homotopic, then they are isotopic through conformal embeddings. Furthermore, we show that the space of all conformal embeddings in a given homotopy class is homotopy equivalent to a point, a circle, a torus, or the unit tangent bundle of the codomain, depending on the induced homomorphism on fundamental groups. Quadratic differentials play a central role in the proof.
\end{abstract}

\maketitle

\tableofcontents

\section{Introduction} 

Loosely speaking, the \emph{$1$-parametric $h$-principle} is said to hold for a class of maps bet\-ween manifolds if the only obstructions to connecting two maps in the class through maps in the same class are topological \cite[p.60]{Eliashberg}. For example, the $1$-parametric $h$-principle holds for immersions of $S^2$ in $\RR^3$, so that the standard sphere can be turned inside out via immersions. This is known as Smale's paradox. Of course, there are situations where the $1$-parametric $h$-principle fails due to geometric obstructions. A famous example is Gromov's symplectic camel theorem, which says that one cannot move a closed $4$-dimensional ball through a hole in a wall in $\RR^4$ via symplectic embeddings if the ball is bigger than the hole. 

For us, a \emph{Riemann surface} is a connected $1$-dimensional complex manifold. A \emph{finite Riemann surface} is a Riemann surface with finitely generated fundamental group. Equivalently, a finite Riemann surface is a closed Riemann surface with finitely many points and closed disks removed. This should not be confused with the notion of a Riemann surface of \emph{finite type}, which is a closed Riemann surface with finitely many points removed. A \emph{conformal embedding} between Riemann surfaces is an injective holomorphic map. 

In this paper, we prove that the $1$-parametric $h$-principle holds for conformal embeddings between finite Riemann surfaces.

\begin{theorem}[The holomorphic couch theorem] \label{thm:CouchThm}
If two conformal embeddings between finite Riemann surfaces are homotopic, then they are isotopic through conformal embeddings.
\end{theorem}

We think of the codomain as a house and the domain as a couch that we want to move around in the house without changing its holomorphic structure. Hence the name ``holomorphic couch".

Given finite Riemann surfaces $X$ and $Y$, and a topolo\-gical embedding $h:X \to Y$, we define $\emb(X,Y,h)$ to be the set of all conformal embeddings homotopic to $h$. We equip this set with the compact-open topology, which is the same as the topo\-logy of uniform convergence on compact sets with respect to any metric inducing the topo\-logy given on $Y$. Theorem \ref{thm:CouchThm} is equivalent\footnote{Since $X$ is locally compact Hausdorff, a path $[0,1]\to \maps(X,Y)$ is the same as a homotopy $X \times [0,1] \to Y$ \cite[p.287]{Munkres}.} to the statement that $\emb(X,Y,h)$ is path-connected whenever it is non-empty. 

In fact, we prove a stronger result. Namely, we determine the homotopy type of the space $\emb(X,Y,h)$. The answer depends on the image of $h$ at the level of fundamental groups. We say that $h$ is \emph{trivial}, \emph{cyclic}, or \emph{generic} if the image of the induced homomorphism $\pi_1(h) : \pi_1(X,x) \to \pi_1(Y,h(x))$ is trivial, infinite cyclic, or non-abelian, respectively. If $Y$ is a torus, then the image of $\pi_1(h)$ can be isomorphic to $\ZZ^2$, but we do not need to distinguish this case. 

\begin{theorem}[Strong holomorphic couch] \label{thm:defret} Let $h:X \to Y$ be a topological embedding between finite Riemann surfaces such that $\emb(X,Y,h)$ is non-empty. Then $\emb(X,Y,h)$ is homotopy equivalent to either the unit tangent bundle of $Y$, a circle, or a point, depending on whether $h$ is trivial, cyclic, or generic. This is unless $Y$ is a torus and $h$ is non-trivial, in which case $\emb(X,Y,h)$ is homotopy equivalent to a torus.
\end{theorem}

If $h$ is generic, then $\emb(X,Y,h)$ is contractible. This is the most in\-te\-res\-ting case; the other cases are either analogous or corollaries (see Section \ref{sec:trivialcases}). The rest of the introduction outlines the proof of Theorem \ref{thm:defret} in the case that $h$ is generic, although we state some intermediate results in greater generality. 

\subsection{Ioffe's theorem}

A \emph{Teichm\"uller embedding} between finite Riemann surfaces $X$ and $Y$ is an injective continuous map $f:X\to Y$ for which there exists a constant $K \geq 1$ and non-zero integrable holomorphic quadratic differentials on $X$ and $Y$ that extend to be real and non-negative along the ideal boundary, such that $f$ is locally of the form $x+iy \mapsto Kx+iy$ in natural coordinates and $Y \setminus f(X)$ is a finite union of points and horizontal arcs. 

Our main tool is a theorem of Ioffe which says that extremal quasiconformal embeddings and Teichm\"uller embeddings are one and the same \cite{Ioffe}.

\begin{theorem}[Ioffe] \label{thm:IoffeIntro}
Let $f:X\to Y$ be a quasiconformal embedding between finite Riemann surfaces which is not conformal. Then $f$ has minimal dilatation (i.e. is closest to being conformal) in its homotopy class if and only if it is a Teichm\"uller embedding.
\end{theorem}

We reprove this theorem in Section \ref{sec:Ioffe}. In the special case where $X$ and $Y$ are closed, this reduces to Teich\-m\"uller's celebrated theorem, since an embedding between closed surfaces is a homeo\-morphism. What is different from Teichm\"uller's theorem, however, is that Teichm\"uller embeddings are not ne\-cessarily unique in their homotopy class---even after ruling out the obvious counterexamples. This is an important issue which we discuss next.

\subsection{Slit mappings}

A \emph{slit mapping} is a conformal Teichm\"uller embedding, i.e., with stretch factor $K=1$. In this case, the quadratic differential on $X$ is the pull-back of the quadratic differential on $Y$ by the slit mapping. 

We show that if $\emb(X,Y,h)$ contains a slit mapping, then eve\-ry ele\-ment of $\emb(X,Y,h)$ is a slit mapping and $\emb(X,Y,h)$ is naturally ho\-meo\-mor\-phic to a point, a compact interval, a circle, or a torus.

\begin{theorem}[Slit mappings are almost rigid] \label{thm:AlmostRigidity}
Let $h:X \to Y$ be any embedding between finite Riemann surfaces. Suppose that $\emb(X,Y,h)$ contains a slit mapping with respect to a quadratic differential $\psi$ on $Y$. Then every $f\in \emb(X,Y,h)$ is a slit mapping with respect to $\psi$. Moreover, for every $x \in X$, the evaluation map $\emb(X,Y,h) \to Y$ sending $f$ to $f(x)$ is a homeomorphism onto its image. This image is equal to $Y$ if $Y$ is a torus, is a compact horizontal arc or a point if $h$ is generic, and is a horizontal circle if $h$ is cyclic and $Y$ is not a torus.
\end{theorem}

A better way to say this is: except in the torus case, any conformal embedding homotopic to a slit mapping differs from the latter by a horizontal translation which can be performed gradually (Theorem \ref{thm:isotopythruslits}). An analogous statement holds for Teichm\"uller embeddings of dilatation $K>1$. 

This theorem corrects the wrong statement from \cite{Ioffe} that $\emb(X,Y,h)$ is a single point if $h$ is generic and $\emb(X,Y,h)$ contains a slit mapping. Indeed, there are easy examples showing that $\emb(X,Y,h)$ can be a non-degenerate interval (see Subsection \ref{subsec:nonunique}). 

We prove Theorem \ref{thm:AlmostRigidity} in Section \ref{sec:rigidity}. Observe that Theo\-rem \ref{thm:AlmostRigidity} implies Theo\-rem \ref{thm:defret} when\-e\-ver $\emb(X,Y,h)$ contains a slit mapping. If $h$ is generic and $\emb(X,Y,h)$ is non-empty but does not contain a slit mapping, the idea is to en\-lar\-ge $X$ un\-til it barely fits in $Y$, then appeal to Theorem \ref{thm:AlmostRigidity} for the enlarged surface.

\subsection{Modulus of extension}

Given a finite Riemann surface $X$ with non-empty ideal boundary, we define a $1$-parameter family of enlargements of $X$ as follows. We first choose an analytic parametrization $S^1 \to C$ of each ideal boundary component $C$ of $X$. Then, for every $r \in (0,\infty]$, we let $X_r$ be the bordered surface $X \cup \partial X$ with a copy of the cylinder $S^1 \times [0,r)$ glued to each ideal boundary component along $S^1 \times \{0\}$ via the fixed parametrization $S^1 \to C$. We also let $X_0=X$. 

We say that a cyclic embedding is \emph{parabolic} if its image on $\pi_1$ is generated by a loop around a puncture. Denote by \myhypertarget{h1}{(\textbf{H})} the hypothesis that:
\begin{itemize}
\item $h:X\to Y$ is a non-trivial and non-parabolic embedding between finite Riemann surfaces;
\item $X$ has non-empty ideal boundary\footnote{If $X$ has finite type, then Theorem \ref{thm:defret} is easy (see Subsection \ref{subsec:aut}).};
\item $\{X_r\}_{r \in [0,\infty]}$ is a $1$-parameter family of enlargements of $X$;
\item $\emb(X,Y,h)$ is non-empty.
\end{itemize}

Under hypothesis \myhyperlink{h1}, we define the \emph{modulus of extension $\m(f)$}  of any $f \in \emb(X,Y,h)$ as the supremum of the set of $r \in [0,\infty]$ such that $f$ extends to a conformal embedding of $X_r$ into $Y$. Montel's theorem in complex analysis implies that:
\begin{itemize}
\item for every $f\in \emb(X,Y,h)$ the supremum $\m(f)$ is achieved by a unique conformal embedding $\what f : X_{\m(f)} \to Y$ extending $f$;
\item $\emb(X,Y,h)$ is compact;
\item $\m$ is upper semi-continuous.
\end{itemize}

In particular, $\m$ achieves its maximum value over $\emb(X,Y,h)$. Using Ioffe's theorem, it is not too hard to show that if $\m$ attains its maximum at $f$, then $\what f$ is a slit mapping. We prove in Section \ref{sec:localmax} that the same holds if $f$ is only assumed to be a local maximum of $\m$.

\begin{theorem}[Characterization of local maxima] \label{thm:localmaxintro} Assume hypothesis \myhyperlink{h1}. Let $f\in \emb(X,Y,h)$ be a local ma\-ximum of $\m$ such that $\m(f)<\infty$ and let $\what f$ be the conformal extension of $f$ to $X_{\m(f)}$. Then $\what f$ is a slit mapping. Conversely, if $g: X_r \to Y$ is a slit mapping such that $g|_X$ is homotopic to $h$, then $g|_X$ is a global ma\-xi\-mum of $\m$. The set $M$ of all local maxima of $\m$ is homeomorphic to a point, a compact interval, a circle, or a torus, and $\m$ is constant on $M$.
\end{theorem}

The initial motivation for studying $\m$ was to think of it as a Morse function for the space $\emb(X,Y,h)$. In an ideal world, flowing along the gradient of $\m$ would yield a deformation retraction of $\emb(X,Y,h)$ into $M$. This does not make sense, however, since $\m$ is not even continuous unless it is constant equal to zero. In any case, the connectedness of $\emb(X,Y,h)$ is a direct consequence of Theorem \ref{thm:localmaxintro} (see Section \ref{sec:localmax}).

\begin{theorem} \label{thm:connected}
Under hypothesis \myhyperlink{h1}, the spa\-ce $\emb(X,Y,h)$ is connected.
\end{theorem}

To improve upon this, we show that there are no local obstructions to contractibility.

\subsection{Where can one point go?}

Throughout this subsection, we assume that:
\begin{itemize}
\item $h:X\to Y$ is a generic embedding between finite Riemann surfaces;
\item $X$ has non-empty ideal boundary;
\item $\emb(X,Y,h)$ is non-empty and does not contain any slit mapping;
\end{itemize}
which we call hypothesis \myhypertarget{h2}{(\textbf{H'})}. Note that $\myhyperlink{h2} \Rightarrow \myhyperlink{h1}$.

Given a point $x \in X$, we are interested in set of points in $Y$ where $x$ can be mapped by the elements of $\emb(X,Y,h)$. It is convenient to also keep track of how $x$ gets mapped to a given $y\in Y$ in the following sense. If $f \in \emb(X,Y,h)$, then by definition there exists a homotopy $$H:X \times [0,1] \to Y$$ from $h$ to $f$. Since $h$ is generic, the homotopy class rel endpoints of the path $t \mapsto H(x,t)$ from $h(x)$ to $f(x)$ does not depend on the particular choice of $H$ (see Lemma \ref{lem:lift}). Denote the homotopy class of that path by $\lift_x(f)$. If the point $x\in X$ is kept fixed, $\lift_x(f)$ represents an element of the universal cover of $Y$ based at $h(x)$. Since $Y$ has non-abelian fundamental group, its universal cover is conformally equivalent to the unit disk $\DD$. 

The map $\lift_x : \emb(X,Y,h) \to \DD$ is continuous, and we call its image $\blob(x,X,Y,h)$. The blob is simpler than the image of the evaluation map in much the same way as Teichm\"uller space is simpler than moduli space. Indeed, the blob is as simple as can be.

\begin{theorem}[The blob is a disk] \label{thm:blobisdisk}
Under hypothesis \myhyperlink{h2}, $\blob(x,X,Y,h)$ is homeo\-morphic to a closed disk for any $x \in X$.
\end{theorem}

The proof has four steps:
\begin{itemize}
\item the blob is compact and connected (because $\emb(X,Y,h)$ is);
\item the blob is semi-smooth;
\item every semi-smooth subset of $\RR^2$ is a $2$-manifold with boundary;
\item there are no holes in the blob.
\end{itemize} 

We refer the reader to Sections \ref{sec:blob} to \ref{sec:noholes} for the definition of a semi-smooth set and the breakdown of these steps. Suffice it to say that Ioffe's theorem implies that points on the boundary of $\blob(x,X,Y,h)$ come from conformal embeddings $f$ whose restriction $f^\star :X \setminus \{x\} \to Y\setminus \{f(x)\}$ is a slit mapping. We then use a variational formula for extremal length to deduce information about the shape of $\blob(x,X,Y,h)$ near such boundary points.

\subsection{Moving one point at a time}

We now explain how to deduce that $\emb(X,Y,h)$ is contractible from the previous results, still assuming hypothesis \myhyperlink{h2}. Pick a countable dense set $\{x_1,x_2,...\}$ in $X$ and let $F \in \emb(X,Y,h)$ be any conformal embedding. We define a deformation retraction of the space $\emb(X,Y,h)$ into $\{F\}$ by moving one point at a time. Given a map $f$ in $\emb(X,Y,h)$, we join $\lift_{x_1}(f)$ to $\lift_{x_1}(F)$ by a path $\gamma_1$ in $\blob(x_1,X,Y,h)$. Such a path exists since $\blob(x_1,X,Y,h)$ is homeomorphic to a closed disk.

For every $t \in [0,1]$, we then look at where $x_2$ can go under maps $g$ in $\emb(X,Y,h)$ which satisfy $\lift_{x_1}(g)=\gamma_1(t)$. This defines a new kind of blob, call it $\blob_t(x_2)$. We show that $\blob_t(x_2)$ moves continuously with $t$, which allows us to construct a second path $\gamma_2$ from $\lift_{x_2}(f)$ to $\lift_{x_2}(F)$ with the property that for every $t \in [0,1]$, the point $\gamma_2(t)$ belongs to $\blob_t(x_2)$.

Proceeding by induction, we obtain a sequence of paths $\gamma_1,\gamma_2,\ldots$ such that for every $n \in \NN$ and every $t\in [0,1]$, there exists at least one map $f_t^n$ in $\emb(X,Y,h)$ such that $\lift_{x_j}(f_t^n)=\gamma_j(t)$ for every $j\in\{1,\ldots,n\}$. If we fix $t$ and pass to a subsequence, we get some limit $f_t \in \emb(X,Y,h)$ for which $\lift_{x_j}(f_t)=\gamma_j(t)$ for every $j \in \NN$. Since any two limits agree on the dense set $\{x_1,x_2,...\}$, we actually have convergence $f_t^n \to f_t$ as $n \to \infty$ (without passing to a subsequence). By a similar argument, $f_t$ depends continuously on $t$. We thus found a path from $f$ to $F$ in $\emb(X,Y,h)$. 

We construct the paths $\gamma_1, \gamma_2, \ldots$ carefully enough so that they depend continuously on the initial map $f$, hence the path $t \mapsto f_t$ also depends conti\-nuous\-ly on $f$. The end result is a deformation retraction of $\emb(X,Y,h)$ into $\{F\}$. See Section \ref{sec:movingpoints} for details.

\subsection{Notes and references}

The holomorphic couch problem arose in the context of renormalization in complex dynamics. Although Theorem \ref{thm:CouchThm} does not have any direct application to dynamics, some of the tools used here do (see \cite{Dylan}).

The space $\emb(\DD,\CC)$, or rather its subspace $\mathcal S$ of conformal embeddings $f:\DD \to \CC$ satisfying $f(0)=0$ and $f'(0)=1$, was the subject of much interest until the solution of the Bieberbach conjecture by de Branges in 1984. It is easy to see that $\mathcal S$ is contractible in the compact-open topology. On the other hand, $\mathcal S$ has isolated points when equipped with the topology of uniform convergence of the Schwarzian derivative \cite{Zippers}. The literature on the class $\mathcal{S}$ is quite vast. In comparison, not much has been written about conformal embeddings between general Riemann surfaces. Exceptions include \cite{EarleMarden}, \cite{Ioffe}, \cite{Ioffe2}, \cite{JenkinsBook}, \cite{SchifferSpencer} and \cite{Shiba}.

The holomorphic couch problem for embeddings of a multiply-punctured disk in a multiply-punctured sphere was first considered in \cite{Royden2}. How\-ever, the solution presented there relies in part on a ri\-gi\-di\-ty claim \cite{Royden} which is known to be false in ge\-ne\-ral \cite{Jenkins}. 

Theorems \ref{thm:AlmostRigidity} and \ref{thm:localmaxintro} generalize the fact that in the homotopy class of an essential simple closed curve $\alpha$ on a finite Riemann surface $Y$, there exists a unique embedded annulus $X \subset Y$ of largest modulus, and $X$ is the horizontal cylinder of a quadratic differential $\psi$ on $Y$ (known as the Jenkins--Strebel differential corresponding to $\alpha$).

In \cite{EarleMarden}, Earle and Marden consider a functional similar to our modulus of extension $\m$ where they keep the annuli disjoint from $X$. Their construction is more natural than ours as it does not require any choice of parametrization of $\partial X$. We decided to glue the annuli to $X$ in order to deal with a connected surface. The main theorem in \cite{EarleMarden} is analogous to Theorem \ref{thm:localmaxintro}, but its uniqueness statement is false for the same reason that Ioffe's uniqueness statement is. 

There is a plethora of other extremal pro\-blems on Riemann surfaces who\-se solutions involve quadratic differentials (see e.g. \cite{Krushkal} and the references there\-in). These are all examples of ``Teichm\"uller's principle'' \cite[p.48]{JenkinsBook}.

For the class $\mathcal{S}$ of normalized univalent functions from $\DD$ to $\CC$, a suitable version of the blob is actually a round disk. More precisely, for every $z \in \DD$ and every $f \in \mathcal S$, the quantity $w=\log(f(z)/z)$ satisfies
$$
\left| w - \log \frac{1}{1-|z|^2}\right| \leq \log \frac{1+|z|}{1-|z|}
$$
and every value $w$ satisfying the inequality is achieved for some $f\in \mathcal S$. This was proved by Grunsky in 1932 (see \cite[p.323]{Duren}). The blob for $K$-quasiconformal homeomorphisms of the disk with prescribed boundary values was studied in \cite{Strebel2} (see also \cite{EarleLakic} for a generalization to arbitrary hyperbolic surfaces). Our approach for proving that the blob is homeomorphic to a closed disk seems similar to Strebel's, but the context is different.

The idea of moving one point at a time to get an isotopy is reminiscent of the finite ``holomorphic axiom of choice'' used by Slodkowski to extend holomorphic motions \cite{Slodkowski}. Our isotopies are holomorphic in the space variable and continuous in the time variable, whereas holomorphic motions are the other way around.

\begin{acknowledgements}
The results of this paper are adapted from the author's Ph.D. thesis at the Gra\-dua\-te Center of the City University of New York. I thank my advisor Jeremy Kahn for sugges\-ting this problem and for his continued interest and support. I also thank Dylan Thurston, Kevin Pilgrim, Frederick Gardiner and Patrick Hooper for useful conversations, Fran\c{c}ois Labourie for pointing out the connection with the $h$-principle, and the anonymous referees for their helpful comments and suggestions.
\end{acknowledgements}

\section{Preliminaries} \label{sec:prelim}

\subsection{Ideal boundary and punctures}

A Riemann surface is \emph{hyperbolic} if its universal covering space is conformally isomorphic to the unit disk $\DD$. The only non-hyperbolic Riemann surfaces are the Riemann sphere $\CHAT$, the complex plane $\CC$, the once-punctured plane $\CC \setminus \{0\}$ and complex tori. A hyperbolic surface $X$ can be regarded as the quotient of its universal covering space $\DD$ by its group of deck transformations $\Gamma$. The \emph{limit set} $\Lambda_\Gamma$ is the set of accumulation points in $\partial \DD$ of the $\Gamma$-orbit of any point $z \in \DD$, and the \emph{set of discontinuity} is $\Omega_\Gamma = \partial \DD \setminus \Lambda_\Gamma$. The \textit{ideal boundary} of $X$ is $\partial X = \Omega_\Gamma / \Gamma$. The union $X \cup\partial X =  (\DD \cup \Omega_\Gamma) / \Gamma$ is naturally a bordered Riemann surface, since $\Gamma$ acts properly discontinuously and analytically on $\DD \cup \Omega_\Gamma$. If $X$ is a finite hyperbolic surface, then $\partial X$ has finitely many connected components, each homeomorphic to a circle.

A \emph{puncture} in a Riemann surface $X$ is an end corresponding to a proper (preimages of compact sets are compact) conformal embedding $\overline \DD \setminus \{0\}\to X$. For example, $\CC$ has one puncture at infinity and $\CC \setminus \{0\}$ has two punctures. For hyperbolic surfaces, punctures are the same as a cusps, or ends with parabolic monodromy. Every puncture can be \emph{filled}, meaning that one can add the missing point and extend the complex structure there. The set of punctures of $X$ is denoted by $\dot{X}$.

Given a finite Riemann surface $X$, we write $\what X =X \cup \partial X \cup \dot X$ for the compact bordered Riemann surface obtained after adding the ideal boundary and filling the punctures (by definition, a non-hyperbolic Riemann surface has empty ideal boundary). Suppose that $\partial X$ is non-empty. Then if we take two copies of $\what X$---the second with reversed orientation---and glue them along $\partial X$ with the identity, we get a closed Riemann surface called the \emph{double} of $\what X$. Because of this construction, many theorems about Riemann surfaces of finite type are also true for finite Riemann surfaces.

\subsection{Automorphisms} \label{subsec:aut}

Given a Riemann surface $X$, let $\aut_0(X)$ be its group of conformal automorphisms homotopic to the identity. It is well-known that:
\begin{itemize}
\item $\aut_0(\CHAT)$ acts simply transitively on ordered triples of distinct points in $\CHAT$;
\item $\aut_0(\CC)$ acts simply transitively on ordered pairs of distinct points in $\CC$;
\item $\aut_0(\CC\setminus\{0\})$ acts simply transitively on $\CC\setminus\{0\}$;
\item if $X$ is a torus, then $\aut_0(X)$ acts simply transitively on $X$; 
\item $\aut_0(\DD)$ acts simply transitively on the unit tangent bundle of $\DD$ (with respect to the hyperbolic metric);
\item  if $r \in [0,1)$, then $\aut_0(\DD \setminus r\overline{\DD})$ is isomorphic to the circle $S^1$.
\end{itemize}
In all other cases, $\aut_0(X)$ only contains the identity.

\begin{lemma} \label{lem:selfembisidentity}
Let $X$ be a hyperbolic Riemann surface not isomorphic to $\DD$ nor $\DD \setminus \{0\}$, and let $h: X \to X$ be a holomorphic map homotopic to the identity. Then $h$ is equal to the identity unless $X$ is an annulus and $h$ is a rotation.
\end{lemma}

In this lemma, $X$ is not assumed to be finite and $h$ is not assumed to be either injective or surjective.

\begin{proof}
By the Schwarz lemma, $h$ is $1$-Lipschitz with respect to the hyperbolic metric. Therefore, if $\alpha$ is a closed geodesic in $X$, then $h(\alpha)$ is at most as long as $\alpha$. But geodesics minimize length in their homotopy class, so that $h(\alpha)=\alpha$. In particular, $h$ is an isometry along $\alpha$. 

If $X$ is an annulus, then it contains a unique simple closed geodesic $\alpha$. We can post-compose $h$ by a rotation $r$ of $X$ so that $r \circ h$ is equal to the identity on $\alpha$ and hence on all of $X$ by the identity principle. 

If $X$ is not an annulus, then it contains a closed geodesic $\alpha$ which self-intersects exactly once. Then $h$ fixes this self-intersection point, thus all of $\alpha$ pointwise, and hence all of $X$ pointwise by the identity principle.
 \end{proof}

It follows that $\aut_0(X)$ is path-connected for any Riemann surface $X$. More precisely, $\aut_0(X)$ is
\begin{itemize}
\item homotopy equivalent to the unit tangent bundle of $X$ if $\pi_1(X) = \{ 0 \}$;
\item homotopy equivalent to $S^1$ if $\pi_1(X) \cong \ZZ$;
\item homeomorphic to $S^1 \times S^1$ if $\pi_1(X) \cong \ZZ^2$;
\item a point if $\pi_1(X)$ is non-abelian.
\end{itemize}
Theo\-rems \ref{thm:CouchThm} and \ref{thm:defret} can be viewed as generalizations of this. When the domain has finite type, Theo\-rem \ref{thm:defret} actually follows from the above.

\begin{lemma} \label{lem:emptyidealbdry}
Suppose that $h:X \to Y$ is a topological embedding bet\-ween finite Riemann surfaces, where $\partial X$ is empty and $Y$ is not the sphere with at most $2$ punctures nor a torus. Then $\emb(X,Y,h)$ contains at most one element.
\end{lemma}
\begin{proof}
Since conformal embeddings send punctures to punctures or regu\-lar points, every $f$ in $\emb(X,Y,h)$ extends to a conformal embedding $\what f$ from $X \cup \dot{X}$ to $Y \cup \dot{Y}$. As $\what X = X \cup \dot{X}$ is compact and $\what{f}$ is open, the map $\what{f}$ is surjective. In particular, $Y \cup \dot{Y}=\what Y$ is compact. If $f,g \in \emb(X,Y,h)$, then the inverses of the extensions $\what f,\what g : \what X \to \what Y$ are homotopic rel $\dot{Y}$. The composition $\what g \circ \what f^{-1}$ is thus a conformal automorphism of $\what Y$ homotopic to the identity rel $\dot{Y}$. If $\what Y$ has genus $0$, then $\dot Y$ must contain at least $3$ points by hypothesis, so that $\what f = \what g$. If $\what Y$ has genus $1$, then $\dot Y$ contains at least $1$ point, so that $\what f = \what g$. If $\what Y$ has genus at least $2$, then $\what f = \what g$ by Lemma \ref{lem:selfembisidentity}.
 \end{proof}

By the same argument, if $\partial X = \varnothing$ and $Y$ is the sphere with at most $2$ punctures or a torus, then $\emb(X,Y,h)$ is either empty or homeomorphic to $\aut_0(Y)$. Theorem \ref{thm:defret} for domains of finite type follows easily.

\begin{proof}[Proof of Theorem \ref{thm:defret} in the case where $\partial X = \varnothing$.]

The image of any element in $\emb(X,Y,h)$ has finite complement in $Y$ so that $\pi_1(h)$ surjects onto $\pi_1(Y)$. If $\pi_1(Y)$ is non-abelian, then $h$ is generic and $\emb(X,Y,h)$ is a singleton by Lemma \ref{lem:emptyidealbdry}. Otherwise $\emb(X,Y,h)$ is homeomorphic to $\aut_0(Y)$, whose homotopy type was described above. 
\end{proof}

\subsection{Montel's theorem}

The simplest version of Montel's theorem says that the set of all holomorphic maps from $\DD$ to $\overline{\DD}$ is compact. This implies a similar result for holomorphic maps between arbitrary hyperbolic surfaces by lifting to the universal covers. 

A sequence of maps $f_n:X \to Y$ between Riemann surfaces \emph{diverges locally uniformly} if for every compact sets $K \subset X$ and $L\subset Y$, the sets $f_n(K)$ and $L$ are disjoint for all large enough $n$. A set $\mathcal{F}$ of maps between two Riemann surfaces $X$ and $Y$ is \emph{normal} if every sequence in $\mathcal{F}$ admits either a locally uniformly convergent subsequence or a locally uniformly divergent subsequence.

\begin{theorem}[Montel's theorem]
If $X$ and $Y$ are hyperbolic surfaces, then every set of holomorphic maps from $X$ to $Y$ is normal.
\end{theorem}

See \cite[p.34]{Milnor}. Note that the limit of a convergent sequence of holomorphic maps is holomorphic. If every map in the sequence is injective, then the limit is either injective or constant. If every map in the sequence is locally injective, then the limit is either locally injective or constant.

\subsection{Quasiconformal maps}

Let $K\geq 1$. A \emph{$K$-quasiconformal map}\footnote{This notion is usually called ``quasiregular map'', and the expression ``quasiconformal map'' is usually reserved for homeomorphisms. It seems more convenient to modify the noun ``map'' instead, to indicate further attributes. For example we will use ``quasiconformal immersion'', ``quasiconformal embedding'', or ``quasiconformal homeomorphism'' for a quasiconformal map which is a local embedding, an embedding, or a homeomorphism respectively.} between Riemann surfaces is a map $f$ such that in charts, its first partial derivatives in the distributional sense are locally in $L^2$ and the formal matrix $\mathrm{d}f$ of partial derivatives satisfies the inequality $\|\mathrm{d}f\|^2 \leq K \det(\mathrm{d}f)$ almost everywhere. For almost every point $z$, the real linear map $\mathrm{d}_z f$ sends circles in the tangent plane at $z$ to ellipses of eccentricity $\|\mathrm{d}_zf\|^2/\det(\mathrm{d}_zf)$ in the tangent plane at $f(z)$, and this ratio is called the \emph{pointwise dilatation} of $f$ at $z$. The \emph{dilatation} of $f$, denoted $\Dil(f)$, is the smal\-lest $K\geq 1$ for which $f$ is $K$-quasiconformal. This is the same as the essential supremum of the pointwise dilatation of $f$. 

A \emph{Beltrami form} on a Riemann surface $X$ is a map $\mu : TX\setminus \vec{0} \to \CC$ such that $\mu(\lambda \mathbf{v})= (\overline \lambda / \lambda)\mu(\mathbf{v})$ for every $\mathbf{v} \in TX\setminus \vec{0}$ and every $\lambda \in \CC \setminus \{0\}$ where $\vec{0}$ denotes the zero section. In charts, the \emph{Wirtinger derivatives} of a quasiconformal map $f$ are   
$$\partial f =\frac{1}{2}\left(\frac{\partial f}{\partial x}-i\frac{\partial f}{\partial y}\right) \quad \text{and} \quad \dbar f =\frac{1}{2}\left(\frac{\partial f}{\partial x}+i\frac{\partial f}{\partial y}\right).$$
The ratio $\dbar f / \partial f$ is naturally a Beltrami form, and is sometimes called the \emph{Beltrami coefficient} of $f$. The Beltrami coefficient of $f$ encodes the field of ellipses in $TX$ which $\mathrm{d}f$ sends to circles. 

The measurable Riemann mapping theorem of Morrey and Ahlfors--Bers says that every measurable ellipse field with bounded eccentricity is the Beltrami coefficient of a quasiconformal homeomorphism.

\begin{theorem}[The measurable Riemann mapping theorem] \label{thm:mrmt}
Let $X$ be a Riemann surface and let $\mu$ be a measurable Beltrami form on $X$ such that $\|\mu\|_\infty < 1$. Then there exists a Riemann surface $Y$ and a quasiconformal homeo\-morphism $f:X \to Y$ such that $\dbar f / \partial f = \mu$ almost everywhere. The surface and the quasiconformal homeomorphism are unique up to post-composition with a conformal homeomorphism.
\end{theorem}

An important consequence is the following factorization principle. Suppose that $f:X \to Y$ is a quasiconformal map. Then $f=g \circ h$ where $h: X \to Z$ is a quasiconformal homeomorphism and $g:Z \to Y$ is holomorphic. Indeed, we can take $h$ to be the solution of the Beltrami equation with $\mu = \dbar f / \partial f$ and let $g = f \circ h^{-1}$. 

Another fundamental property of $K$-quasiconformal homeomorphisms is compactness under appropriate normalization \cite[p.32]{AhlforsLectures}.

\begin{theorem}
For every $K\geq 1$, the space of all $K$-quasiconformal homeomorphisms from $\DD$ to $\DD$ fixing $0$ is compact.
\end{theorem}

Lastly, we will use the fact that quasiconformal homeomorphisms send punctures to punctures. As a consequence, any quasiconformal embedding bet\-ween finite Riemann surfaces extends to a quasiconformal embedding between the surfaces with their punctures filled.  

The reader may consult \cite{AhlforsLectures} for more background on quasiconformal homeo\-morphisms.

\subsection{Quadratic differentials}

A \emph{quadratic differential} on a Riemann surface $X$ is a map $\vphi : TX \to \CHAT$ such that $\vphi(\lambda \mathbf{v}) = \lambda^2 \vphi(\mathbf{v})$ for every $\mathbf{v} \in TX$ and every $\lambda \in \CC$. A quadratic differential on $X$ is \emph{holomorphic} (resp. \emph{meromorphic}) if for every open set $U \subset X$, and every holomorphic vector field $\mathbf{v}:U \to TU$, the function $\vphi \circ \mathbf{v} : U \to \CC$ is holomorphic (resp. meromorphic). All quadratic differentials in this paper will be holomorphic or meromorphic. The \emph{pull-back} $f^* \vphi$ of a quadratic differential $\vphi$ by a holomorphic map $f$ is defined in the usual way by the formula $f^* \vphi (\mathbf{v}) = \varphi(\mathrm{d}f(\mathbf{v}))$.

A vector $\mathbf{v}\in TX$ is \emph{horizontal} (resp. \emph{vertical}) for $\varphi$ if $\vphi(\mathbf{v})>0$ (resp. $\vphi(\mathbf{v})<0$). Let $I \subset \RR$ be an interval. A piecewise smooth arc $\gamma : I \to X$ is \emph{horizontal} (resp. \emph{vertical}) if $\gamma'(t)$ is horizontal (resp. vertical) wherever it is defined. Such a trajectory is called \emph{regular} if it does not contain any zero or pole of $\vphi$. The absolute value $|\varphi|$ is an area form, and its integral $\|\vphi\|=\int_X |\vphi|$ is the \emph{norm} of $\vphi$. For a finite Riemann surface $X$, we denote by $\qd(X)$ the set of all integrable holomorphic quadratic differentials $\vphi$ on $X$ which extend analytically to the ideal boundary of $X$, and such that $\vphi(\mathbf{v})\in \RR$ for every vector $\mathbf{v}$ tangent to $\partial X$. Every $\vphi \in \qd(X)$ extends to a meromorphic quadratic differential on $\what X$ with at most simple poles on $\dot{X}$. The set $\qd^+(X)$ is similarly defined, but with the additional requirements that $\vphi \geq 0$ along $\partial X$ and that $\vphi$ is not identically zero. The set $\qd(X)$ is a real vector space inside of which $\qd^+(X)$ forms a convex cone. 

For every simply connected open set $U \subset \what X$ where a quadratic differential $\varphi$ does not have any zero or pole, there exists a locally injective holomorphic map $z: U \to \CC$ such that $\vphi=\mathrm{d}z^2$. The map $z$ is unique up to translation and sign and is called a \emph{natural coordinate} when it is injective. If $\vphi \in \qd(X)$, then the atlas of natural coordinates for $\vphi$ is a \emph{half-translation structure} on $\what X$ minus the singularities of $\vphi$. The transition maps for this atlas have the form $z \mapsto \pm z + c$. Such a structure induces a flat geometry with cone points on $\what X$. We return to this geometry at the end of the present section. The standard reference for this material is \cite{Strebel}.

\subsection{Teichm\"uller's theorem}

A \emph{Teichm\"uller homeomorphism} (usually called Teichm\"uller map) between finite Riemann surfaces $X$ and $Y$ is a homeomorphism $f:X \to Y$ such that there exists a constant $K>1$ and non-zero $\vphi \in \qd(X)$ and $\psi \in \qd(Y)$ such that $f$ is locally of the form $x+iy \mapsto Kx + iy$ in natural coordinates for $\vphi$ and $\psi$, up to sign and translation. Such a homeomorphism is $K$-quasiconformal with constant pointwise dilatation.

The following theorems of Teichm\"uller are of central importance.

\begin{theorem}[Teichm\"uller's existence theorem]  \label{thm:TeichExistence}
Let $h$ be a quasiconformal homeomorphism between finite Riemann surfaces. If there is no conformal homeomorphism homotopic to $h$, then there is a Teichm\"uller homeomorphism homotopic to $h$. 
\end{theorem}

\begin{theorem}[Teichm\"uller's uniqueness theorem] \label{thm:TeichUniqueness}
Let $f:X \to Y$ be a Teichm\"uller homeomorphism of dilatation $K$ between finite Riemann surfaces. If $g$ is a $K$-quasiconformal homeomorphism homotopic to $f$, then $g \circ f^{-1}$ is a conformal automorphism of $Y$ homotopic to the identity. If $Y$ is not an annulus nor a torus, then $g=f$.  
\end{theorem}

Teich\-m\"ul\-ler's theorem is usually stated and proved for closed Riemann surfaces, but the general case follows from the closed case by doubling across the ideal boundary and by taking a branched cover of degree $2$ or $4$ ramified at the punctures \cite{Ahlfors}.

\subsection{Teichm\"uller spaces}

Let $S$ be a finite Riemann surface. The \emph{Teichm\"uller space} $\teich(S)$ is defined as the set of pairs $(X,f)$ where $X$ is a finite Riemann surface and $f:S\to X$ is a quasiconformal homeomorphism, modulo the equivalence relation $(X,f) \sim (Y,g)$ if and only if $g \circ f^{-1}$ is homotopic to a conformal homeomorphism. The equivalence class of $(X,f)$ is denoted $[X,f]$, or just $X$ when the \emph{marking} $f$ is implicit. The \emph{Teichm\"uller distance} between two points of $\teich(S)$ is defined as
$$
d([X,f],[Y,g])= \frac{1}{2} \inf  \log \Dil(h)
$$
where the infimum is taken over all quasiconformal homeomorphisms $h$ homotopic to $g \circ f^{-1}$. By Teichm\"uller's theorem, the infimum is realized by a (usually unique) quasiconformal homeomorphism $h$ which is either conformal or a Teichm\"uller homeomorphism. 

The space $\teich(S)$ is a contractible real-analytic manifold of finite dimension. Let $\mathcal{M}(X)$ denote the space of essentially bounded Beltrami forms on $X \in \teich(S)$. By the measurable Riemann mapping theorem, the tangent space to $\teich(S)$ at $X$ can be identified with the quotient of $\mathcal{M}(X)$ by its subspace $\mathcal{M}_0(X)$ of infinitesimally trivial deformations. There is a na\-tural pairing bet\-ween $\mathcal{M}(X)$ and $\qd(X)$ given by
$$
\langle \mu , \vphi \rangle = \re \int_X \mu \vphi,
$$
and it turns out that $\mathcal{M}_0(X) = \qd(X)^\perp$ with respect to this pairing. The tangent and cotangent spaces to $\teich(S)$ at $X$ are thus isomorphic to $\mathcal{M}(X) / \qd(X)^\perp$ and $\qd(X)$ respectively. See \cite{Earle2} and \cite{Earle} for more details.

\begin{remark}
In the literature, $\teich(S)$ is called the \emph{reduced} Teich\-m\"uller space of $S$. The (unreduced) Teichm\"uller space of $S$ is defined similarly, except that pairs $(X,f)$ and $(Y,g)$ are declared equivalent if $g\circ f^{-1}$ is homotopic to a conformal homeomorphism rel $\partial X$. This unreduced Teichm\"uller space is always a contractible complex manifold, but has infinite dimension whenever the ideal boundary of $S$ is non-empty. We will only consider reduced Teichm\"uller spaces in this paper, and will therefore omit the adjective ``reduced''.  
\end{remark}

\subsection{Homotopies}

If two maps are homotopic, then they induce the same homomorphism between fundamental groups, up to conjugation. The converse also holds under appropriate conditions \cite[p.60]{AhlforsLectures} \cite[\S 6]{Bers}. 

\begin{lemma}\label{lem:ahlfors}
Let $X$ be a space which has a universal cover, let $Y$ be a metric space whose universal cover is a uniquely geodesic space in which geodesics depend continuously on endpoints, and let $f_0,f_1:X \to Y$ be continuous maps. Suppose that for some $x \in X$ the induced homo\-mor\-phisms $\pi_1(f_j) : \pi_1(X,x) \to \pi_1(Y,f_j(x))$ agree up to conjugation by a path between $f_0(x)$ and $f_1(x)$. Then $f_0$ and $f_1$ are homotopic.
\end{lemma}
\begin{proof}
Let $\wtilde X$ and $\wtilde Y$ be the universal covers of $X$ and $Y$, and let $\alpha$ be a path connecting $f_0(x)$ to $f_1(x)$ which conjugates the homomorphisms $\pi_1(f_0)$ and $\pi_1(f_1)$. Given a lift $\wtilde f_0 : \wtilde X \to \wtilde Y$, the path $\alpha$ allows us to lift $f_1$ in such a way that $\wtilde f_0$ and $\wtilde f_1$ are equivariant with respect to the same homomorphism of deck groups. The homotopy from $\wtilde f_0$ to $\wtilde f_1$ sending $(x,t) \in \wtilde X \times [0,1]$ to the point at proportion $t$ along the geodesic from $\wtilde f_0(x)$ to $\wtilde f_1(x)$ in $\wtilde Y$ is continuous and equivariant, so it descends to a homotopy from $f_0$ to $f_1$. 
\end{proof}

This is also true if $X$ is a CW-complex and $Y$ is a $K(\pi,1)$ \cite[p.90]{Hatcher}. If $X$ and $Y$ are finite Riemann surfaces and $Y$ is not the sphere, then either hypotheses are satisfied. The most useful consequence for us is that homotopy classes of maps between finite Riemann surfaces are closed.

\begin{corollary} \label{cor:homotopyclassclosed}
Let $X$ and $Y$ be finite Riemann surfaces and let $f_n,f:X \to Y$ be continuous maps such that $f_n \to f$ as $n\to \infty$. Then $f_n$ is homotopic to $f$ for all large enough $n$. 
\end{corollary}
\begin{proof}
First assume that $Y$ is the Riemann sphere. By the Hopf theorem \cite[p.50]{Milnor2}, two maps $X \to Y$ are homotopic if and only if they have the same degree. Since topological degree depends continuously on the map, $f_n$ eventually has the same degree as $f$ and is therefore homotopic to it.

Suppose that $Y$ is not the Riemann sphere. Then the universal cover of $Y$ supports a metric which is proper (whose closed balls are compact) and is uniquely geodesic. In such a metric, geodesics depend continuously on endpoints. Thus the hypotheses of Lemma \ref{lem:ahlfors} are satisfied.

Let $\beta_1, \ldots, \beta_k$ be loops based at $x \in X$ which generate $\pi_1(X,x)$ and let $V$ be a simply connected neighborhood of $f(x)$ in $Y$. Let $n$ be large enough so that $f_n(x) \in V$ and so that there is a homotopy between the restrictions $f_n|\beta_j$ and $f|\beta_j$ which keeps the image of the basepoint $x$ inside $V$ for every $j \in \{1, \ldots, k \}$. Let $\alpha$ be any path from $f(x)$ to $f_n(x)$ in $V$. Then $\alpha * f_n(\beta_j)*\overline{\alpha}$ is homotopic to $f(\beta_j)$ as loops based at $f(x)$, for every $j$. In other words, $\alpha$ conjugates the induced homo\-mor\-phisms $\pi_1(f)$ and $\pi_1(f_n)$.  By Lemma \ref{lem:ahlfors}, $f_n$ is homotopic to $f$. 
\end{proof}

In Teichm\"uller theory, one often goes back and forth between punctures and marked points as convenient. This passage is justified by the fact that quasiconformal homeomorphisms send punctures to punctures. Moreover, homotopies defined in the complement of punctures can be modified as to extend to the punctures. 

\begin{definition}
Let $h:X \to Y$ be a quasiconformal embedding between finite Riemann surfaces. Recall that $h$ extends to a quasiconformal embedding $\overline{h}$ from $X \cup \dot X$ to $Y \cup \dot{Y}.$ We say that a puncture $p \in \dot{X}$ is \emph{essential} if $\overline{h}(p) \in \dot{Y}$ and is \emph{inessential} if $\overline{h}(p) \in Y$.  
\end{definition}

\begin{lemma}\label{lem:punctures}
Let $X$ and $Y$ be finite Riemann surfaces, let $f_0, f_1: X \to Y$ be quasiconformal embeddings, let $\overline f_0$ and $\overline f_1$ be their extensions to $X \cup \dot{X}$, and let $H : X \times [0,1] \to Y$ be a homotopy from $f_0$ to $f_1$. Let $E \subset \dot{X}$ be the set of punctures which are essential with respect to $f_0$. Then there exists a homotopy $\overline{H}:X\cup \dot{X} \times [0,1]\to Y \cup \dot{Y}$ from $\overline f_0$ to $\overline f_1$ which is constant on $E\times [0,1]$, maps $X \times [0,1]$ into $Y$, and whose restriction to $X \times [0,1]$ is homotopic to $H$.
\end{lemma}
\begin{proof}
For each $p\in \dot{X}$, let $D_p \subset X\cup \dot{X}$ be an embedded closed disk such that $D_p \cap \dot{X}  = \{p\}$. Further assume that the disks $D_p$ are all disjoint. The idea is to define $\overline H = H$ on $X \setminus \bigcup_{p\in \dot{X}} D_p \times [0,1]$ and then extend this to a continuous map $ \overline H: D_p \times [0,1] \to Y \cup \dot{Y}$ sending $D_p \setminus \{p \}$ into $Y$ for each $p \in \dot{X}$. If the puncture $p\in \dot{X}$ is inessential, then the annulus $H(\partial D_p \times [0,1])$ is contractible in $Y$ so that $H$ extends to a continuous map $\overline H: D_p \times [0,1] \to Y$. If $p\in \dot{X}$ is essential, then the annulus $H(\partial D_p \times [0,1])$ is contractible in $Y \cup \{f_0(p)\}$ so that $H$ extends to a continuous map $\overline H: D_p \times [0,1] \to Y \cup \{f_0(p)\}$. We can further choose the extension to satisfy $\overline{H}^{-1}(f_0(p))=\{p\} \times [0,1]$.
\end{proof}

\subsection{Geometry of quadratic differentials} \label{subsec:geom}

Let $X$ be a finite Riemann surface and let $\varphi \in \qd(X)\setminus \{0\}$. Then $\varphi$ induces a notion of direction as well as a Riemannian metric with singularities $|\varphi|$ on $X$. This metric is Euclidean except at the zeros of $\varphi$, where it has cone points. At a zero of order $k$, the metric looks like a cone with angle $(k+2)\pi$. Because the total angle at cone points is at least $2\pi$, the induced metric is locally $\cat(0)$. However, the metric $|\varphi|$ is not complete whenever $X$ is not closed. Nevertheless, given an arc $\gamma : [0,1] \to X$ there exists a unique ``geodesic representative'' $\gamma^\dagger:[0,1] \to \what X$ which is a limit of arcs homotopic to $\gamma$ in $X$ and has minimal $|\varphi|$-length among such arcs. To see this, it is convenient to pass to the universal cover first.

Let $\pi: \wtilde{X} \to X$ be the universal covering map and let $\wtilde{\varphi} = \pi^* \varphi$. The Riemannian metric with singularities $|\wtilde{\varphi}|$ induces a distance $d$ on $\wtilde{X}$ in the usual way: for any $x,y\in \wtilde{X}$ the distance $d(x,y)$ is defined as the infimum of $\int_\gamma \sqrt{|\wtilde{\varphi} |}$ over all piecewise smooth paths $\gamma$ between $x$ and $y$. The resulting metric space $(\wtilde{X}, d)$ is $\cat(0)$. It follows that its metric completion $(\wtilde{X}', d')$  is also $\cat(0)$ \cite[Corollary 3.11]{Bridson}. When $X$ is hyperbolic, $\wtilde{X}'$ is the universal cover $\wtilde X$ together with the the set of discontinuity and the parabolic limit points. Observe that $\what X$ is the metric completion of $X$ in the metric $|\vphi|$ and that $\pi$ extends to a continuous map $\pi' : \wtilde{X}' \to \what X$. 

Since $(\wtilde{X}',d')$ is a complete $\cat(0)$ space, there is a unique geodesic segment between any two points of $\wtilde{X}'$. Thus given an arc $\gamma : [0,1] \to X$, we can lift it to an arc $\wtilde{\gamma} : [0,1] \to \wtilde{X}$, find the geodesic $\wtilde\gamma^\dagger: [0,1] \to \wtilde{X}'$ bet\-ween the endpoints of $\wtilde{\gamma}$, then define the ``geodesic representative''\footnote{The quotation marks are because $\gamma^\dagger$ is not necessarily geodesic in $\what X$. For example, it may pass through a pole of $\vphi$.} of $\gamma$ to be $\gamma^\dagger = \pi'\circ \wtilde\gamma^\dagger$. 

By hypothesis, $\vphi$ extends analytically to the ideal boundary of $X$ (if any). Hence $\wtilde \vphi$ extends analytically to the set of discontinuity in $\wtilde{X}'$. Say that a \emph{singularity of $\wtilde \vphi$} is a point in $\wtilde{X}'$ where $\wtilde \vphi$ has a zero or is not defined (those are the parabolic limit points). Since $(\wtilde{X}',d')$ is locally Euclidean, geodesics are straight lines in natural coordinates away from the singularities of $\wtilde \vphi$. If a geodesic passes though a singularity, then there should not be any shortcut on either side, which translates into an angle condition at the singularity. More precisely, an arc $\gamma: [0,1] \to \wtilde{X}'$ is geodesic if and only if 
\begin{itemize}
\item $\gamma$ is smooth except at singularities of $\wtilde{\varphi}$;
\item the argument of $\wtilde{\varphi}(\gamma'(t))$ is locally constant where $\gamma$ is smooth;
\item at a singularity of $\wtilde{\varphi}$, the cone angle on either side\footnote{At a boundary point of $\wtilde X'$, the cone angle only makes sense on one side of $\gamma$, but we can define the cone angle on the ``other side'' to be $+\infty$ by convention.} of $\gamma$ is at least $\pi$.
\end{itemize}

For example, an arc which is horizontal and does not backtrack is geodesic. Actually, any horizontal arc which does not backtrack minimizes horizontal travel (defined as the integral of $|\re \sqrt{\wtilde{\varphi}} |$) between its endpoints, since it is quasi-transverse to the vertical foliation \cite[Proposition 2.5]{HubbardMasur}.

A \emph{geodesic polygon} in a $\cat(0)$ space is a closed curve which is piecewise geodesic. Self-intersecting polygons are allowed. Given a geodesic polygon $P$ in a $\cat(0)$ space, and three consecutive vertices $a$, $b$ and $c$ along $P$, let $\angle abc$ be the Alexandrov angle between the geodesics $[a,b]$ and $[b,c]$ at $b$ \cite[p.184]{Bridson}. If $P$ is a geodesic polygon in $\wtilde{X}'$, then $\angle abc$ is the minimum between $\pi$ and the two cone angles on either side of $[a,b]\cup[b,c]$ at $b$, as measured in the metric induced by $\wtilde \vphi$. 

The next Proposition is a generalization of Strebel's ``divergence principle'' \cite[p.77]{Strebel}. It says that if two geodesic rays $\gamma_0$ and $\gamma_1$ in $\wtilde{X}'$ are such that the angles they form with the geodesic $[\gamma_0(0),\gamma_1(0)]$ sum to at least $\pi$, then the distance between $\gamma_0(t)$ and $\gamma_1(t)$ is non-decreasing as a function of $t \geq 0$. In the original statement, each angle is assumed to be at least $\pi/2$. The result actually holds in any $\cat(0)$ space.

\begin{lemma} \label{lem:quad}
Let $Q$ be a geodesic quadrilateral with vertices $a,b,c,d$ in a $\cat(0)$ space $(M,\mu)$. Suppose that $\mu(a,d)=\mu(b,c)$ and $\angle dab + \angle abc \geq \pi$. Then $\mu(c,d) \geq \mu(a,b)$ with equality only if $Q$ is isometric to a (possibly degenerate) Euclidean paral\-le\-logram.  
\end{lemma}
\begin{proof}
Consider the geodesic triangles $\Delta abc$ and $\Delta acd$. Let $\Delta a'b'c'$ and $\Delta a'c'd'$ be corresponding comparison triangles in $\RR^2$. We may assume that these two triangles are on opposite sides of the segment $[a', c']$ so that the quadrilateral $Q'$ they form is simple. We need to show that $|c' - d'| \geq |a' - b'|$. By the law of cosines, it suffices to prove that $\angle d' a'c' \geq \angle b'c'a'$. Indeed, the triangles $\Delta a'b'c'$ and $\Delta a'c'd'$ share the side $[a',c']$ and the sides $[b',c']$ and $[d',a']$ are congruent. There are two cases to consider depending on whether $Q'$ is convex or not. 

Suppose first that $\angle d'a'c' + \angle c'a'b' \leq \pi$. Note that angles in $\Delta a'b'c'$ and $\Delta a'c'd'$ are are at least as large as corresponding angles in $\Delta abc$ and $\Delta acd$. Moreover, angles are subadditive in the sense that $\angle dab \leq \angle dac + \angle cab$ for instance. Since the sum of $\angle d'a'c'$ and $\angle c'a'b'$ is at most $\pi$, we have
$$\angle d'a'b' = \angle d'a'c' + \angle c'a'b' \geq  \angle dac + \angle cab \geq \angle dab$$
so that 
$$
\angle d'a'b' + \angle a'b'c' \geq \angle dab + \angle abc \geq \pi.
$$
We deduce that
\begin{align*}
0 & \leq \angle d'a'b' + \angle a'b'c' - \pi \\
& = (\angle d'a'c' +\angle c'a'b') + \angle a'b'c' - (\angle c'a'b' + \angle a'b'c'+ \angle b'c'a') \\
& =  \angle d'a'c' - \angle b'c'a'
\end{align*}
as required. 

If $\angle d' a'c' = \angle b'c'a'$, then $Q'$ is a parallelogram (which is possibly contained in a line) and all the above inequalities are equations. Thus the angles in $\Delta abc$ and $\Delta acd$ are the same as in their comparison triangles, which implies that $\Delta abc$ and $\Delta acd$ are isometric to their comparison triangles. By considering the other diagonal of $Q$, we get that any two adjacent sides of $Q$ lie in a flat triangle. The equality $\angle dab = \angle dac + \angle cab $ means that the triangles  $\Delta abc$ and $\Delta acd$ line up, in the sense that the union of their convex hulls is convex. Hence $Q$ spans a Euclidean parallelogram. 

Now suppose that $\angle d'a'c'+ \angle c'a'b' > \pi$. Then 
\begin{align*}
0 &\leq \pi + \angle a'b'c' - \pi \\
&< (\angle d'a'c' +\angle c'a'b') + \angle a'b'c' - (\angle c'a'b' + \angle a'b'c'+ \angle b'c'a') \\
& = \angle d'a'c' - \angle b'c'a'.
\end{align*}
Equality cannot hold in this case. 
 \end{proof}

One can deduce from this the well-known fact that if two simple closed geodesics in $X$ (with respect to $|\vphi|$) are freely homotopic, then they span a Euclidean cylinder.

\section{Ioffe's theorem} \label{sec:Ioffe}

The goal of this section is to reprove Theorem \ref{thm:IoffeIntro} from Ioffe, that is, to characterize quasiconformal embeddings that have minimal dilatation in their homotopy class. The main motivation for doing so is that Ioffe's original statement contains a mistake (we pointed this out in the introduction and will return to this issue in Subsection \ref{subsec:nonunique}). Thus we wanted to work out the proof in detail to make sure that the remaining part was correct. 

The results are stated for quasiconformal embeddings from one finite Riemann surface to another, but they extend to quasiconformal embeddings from a finite union of finite Riemann surfaces to a finite Riemann surface, which is is the context considered in Ioffe's paper \cite{Ioffe}. 

\subsection{Compactness}

Recall that an embedding $h$ between finite Riemann surfaces is \emph{trivial} if the image of $\pi_1(h)$ is trivial, is \emph{cyclic} if the image of $\pi_1(h)$ is infinite cyclic, and is \emph{generic} if the image of $\pi_1(h)$ is non-abelian. An embedding $h$ is called \emph{parabolic} if the image of $\pi_1(h)$ is infinite cyclic and generated by a loop around a puncture. 

The following compactness lemma (cf. \cite[Theorem 1.1]{Ioffe}) guarantees the existence of extremal quasiconformal embeddings, that is, with minimal dilatation in their homotopy class.

\begin{lemma} \label{lem:compactness}
Let $K \geq 1$ and let $h : X \to Y$ be a $K$-quasiconformal embedding bet\-ween finite Riemann surfaces. The space of all $K$-quasiconformal embeddings homotopic to $h$ is compact if and only if $h$ is neither trivial nor parabolic.
\end{lemma}
\begin{proof}
Suppose that $Y$ is equal to $\CHAT$, $\CC$ or $\CC\setminus \{0\}$. Then $h$ is either trivial or parabolic. Moreover, the space of $K$-quasiconformal emdeddings homotopic to $h$ is closed under post-composition with elements of $\aut_0(Y)$, which is non-compact. It follows that the space of $K$-quasiconformal emdeddings homotopic to $h$ is also non-compact. 

Suppose that $Y$ is a torus. The group $\aut_0(Y)$ is homeomorphic to $Y$ and is in particular compact. Let $W$ be the set of all $K$-quasiconformal embeddings homotopic to $h$ and let $N$ be the space of of $K$-quasiconformal embeddings homotopic to $h$ rel $x_0$, where $x_0$ is any point in $X$. Then $W$ is homeomorphic to $\aut_0(Y) \times N$, hence is compact if and only if $N$ is. Now $N$ is the same as the space of $K$-quasiconformal embeddings homotopic to the restriction $h^\star: X \setminus \{x_0\} \to Y \setminus \{h(x_0)\}$. Moreover $h$ non-trivial if and only if $h^\star$ is non-parabolic. The case where $Y$ is a torus thus reduces to the case where $Y$ is a once-punctured torus, hence hyperbolic.

For the rest of the proof, we may assume that $Y$ is hyperbolic. This implies that $X$ is hyperbolic as well, since there is no non-constant quasiconformal map from a non-hyperbolic surface to a hyperbolic surface. 

If $h$ is trivial then its image is contained in a disk in $Y$. In other words, $h$ can be written as $h=F \circ g$ where $F: \DD \to Y$ is a conformal embedding and $g : X \to \DD$ is a $K$-quasiconformal embedding. Consider the sequence $h_n=F_n \circ g$ where $F_n(z)=F(z/n)$. Each $h_n$ is a $K$-quasiconformal embedding homotopic to $h$, but the sequence converges to a constant map. Similarly, if $h$ is parabolic, then we can form a sequence of homotopic $K$-quasiconformal embeddings which diverges to the corresponding puncture. 

If $h$ is neither trivial nor parabolic then it is is either generic, or cyclic but not parabolic. We treat the generic case first and return to the cyclic case at the end. Suppose that $h$ is generic. Pick any basepoint $b_0 \in X$ and let $\pi_X:(\DD,0) \to (X,b_0)$ and $\pi_Y: (\DD,0) \to (Y,h(b_0))$ be universal covering maps with respective deck groups $\deck(X)$ and $\deck(Y)$. Let $\wtilde h$ be the unique lift of $h$ such that $\wtilde{h}(0) = 0$. For any $\alpha \in \deck(X)$, let $\Theta(\alpha)$ be the unique element in $\deck(Y)$ sending $0$ to $\wtilde{h}(\alpha\cdot 0)$. Then $\Theta$ is a homomorphism.

  If $f: X \to Y$ is a $K$-quasiconformal embedding homotopic to $h$, then it lifts to a $K$-quasiconformal immersion $\wtilde f : \DD \to \DD$ satisfying 
$$
\wtilde f (\alpha \cdot z)=\Theta(\alpha)\cdot \wtilde f (z)
$$ 
for every $z \in \DD$ and every $\alpha \in \deck(X)$. The lift $\wtilde f$ is uniquely determined by this equation, since there is only one homotopy class of homotopy from $h$ to $f$ according to Lemma \ref{lem:lift}. We do not really need uniqueness of the lift here (any lift would do), which is why the lemma is postponed to a later section. Observe that $f$ is injective if and only if 
$$
\pi_Y(\wtilde f (z)) = \pi_Y(\wtilde f (w)) \Rightarrow \pi_X(z) = \pi_X(w).
$$

Let $f$ be a $K$-quasiconformal embedding homotopic to $h$ and $\wtilde f$ be its lift as above. Write $\wtilde f = F \circ g$ where $g:\DD \to \DD$ is a $K$-quasiconformal homeomorphism fixing $0$ and $F: \DD \to \DD$ is a holomorphic immersion. As mentioned in Section \ref{sec:prelim}, the space of $K$-quasiconformal homeomorphisms of $\DD$ fixing the origin and the space of holomorphic maps from $\DD$ to $\overline \DD$ are both compact. Thus, given a sequence of $K$-quasiconformal embeddings $f_n$ homotopic to $h$ and corresponding lifts $\wtilde f_n = F_n \circ g_n$ factored as above, we can pass to a subsequence such that $F_n \to F$ and $g_n \to g$ and hence $\wtilde f_n \to \wtilde f := F \circ g.$

We claim that the limit $\wtilde f$ is not constant. Indeed, if there is a constant $w_0 \in \overline \DD$ such that $\wtilde f(z)=w_0$ for every $z \in \DD$, then 
$$w_0 =\wtilde f(\alpha\cdot z)=\lim_{n\to \infty}\wtilde f_n(\alpha\cdot z)=\lim_{n\to \infty}\Theta(\alpha)\cdot \wtilde f_n(z)=\Theta(\alpha)\cdot \wtilde f(z)=\Theta(\alpha)\cdot w_0$$ for every $\alpha \in \deck(X)$. However, a Fuchsian group which fixes a point in $\overline \DD$ is cyclic, contradicting the assumption that $h$ is generic. In particular, the holomorphic map $g$ is not constant, hence has image in $\DD$. By Hurwitz's theorem in complex analysis, $g$ is locally injective. Therefore $\wtilde f=F \circ g$ is a $K$-quasiconformal immersion. Moreover, the equality $\wtilde f (\alpha \cdot z)=\Theta(\alpha)\cdot \wtilde f (z)$ for every $z \in \DD$ and every $\alpha \in \deck(X)$ implies that $\wtilde f$ descends to a $K$-quasiconformal immersion $f : X \to Y$. 

It remains to show that $f$ is injective. If $f(\pi_X(z))=f(\pi_X(w))$, then $\pi_Y(\wtilde f(z))=\pi_Y(\wtilde f(w))$. Since $\wtilde f_n \to \wtilde f$ and since these maps are open, we can find a sequence $z_n$ converging to $z$ and a sequence $w_n$ converging to $w$ such that $\wtilde f_n(z_n) = \wtilde f (z)$  and $\wtilde f_n (w_n)= \wtilde f (w)$ for all $n$ large enough. Then $\pi_Y(\wtilde f_n(z_n))=\pi_Y(\wtilde f_n (w_n))$, which implies that $\pi_X(z_n)=\pi_X(w_n)$ since $f_n$ is injective. Taking the limit as $n \to \infty$, we obtain $\pi_X(z)=\pi_X(w)$. 

Lastly, $f$ is homotopic to $h$ because it is a limit of maps which are (see Corollary \ref{cor:homotopyclassclosed}). We have shown that if $Y$ is hyperbolic and $h$ is generic, then the space of $K$-quasiconformal embeddings homotopic to $h$ is sequentially compact.

Suppose now that $h$ is cyclic but not parabolic. The same construction as above still applies, but the image of $\Theta$ is a cyclic group generated by a hyperbolic element $\beta \in \deck(Y)$. Morover, given an embedding $f$ homotopic to $h$, its lift $\wtilde f$ is only defined up to post-composition with powers of $\beta$. Let $D \subset \DD$ be a fundamental domain for $\beta$ whose closure is disjoint from the fixed points of $\beta$. Let $\{f_n\}_{n\geq 1}$ be a sequence of $K$-quasiconformal embeddings homotopic to $h$. By applying an appropriate power of $\beta$, we can choose a lift $\wtilde f_n$ of $f_n$ such that $\wtilde f_n(0) \in D$. As above, we can extract a subsequence of $\{\wtilde{f}_n\}_{n\geq 1}$ converging to some map $\wtilde{f}$. By our normalization, the limit $\wtilde f(0)$ belongs to the closure $\overline D$, so it is not one of the fixed points of $\beta$. It follows that $\wtilde f$ is not constant. The rest of the argument applies verbatim to conclude that $\wtilde f$ descends to a $K$-quasiconformal embedding $f$ homotopic to $h$.  
\end{proof}

The following corollary is immediate.

\begin{corollary} \label{cor:mindil}
Let $h : X \to Y$ be a non-trivial and non-parabolic quasiconformal embedding between finite Riemann surfaces. Among the quasiconformal embeddings homotopic to $h$, there is one with mini\-mal dilatation.
\end{corollary}
\begin{proof}
Dilatation is lower semi-continuous, hence achieves its minimum on any compact set.
 \end{proof}

\begin{remark}
If $h:X \to Y$ is a quasiconformal embedding between finite Riemann surfaces which is trivial or parabolic, then it is homotopic to a conformal embedding. This is because $h(X)$ is contained in a subset $Z \subset Y$ isomorphic to the disk with at most one puncture or the sphere with at most two punctures, whose Teichm\"uller space is a point. Using Theorem \ref{thm:mrmt}, we can apply a quasiconformal deformation on $Z$ to make $h$ conformal. The hypothesis that $h$ be non-trivial and non-peripheral is thus superfluous in Corollary \ref{cor:mindil}.
\end{remark}

\subsection{Teichm\"uller embeddings}

Let us recall the definition of a Teichm\"uller embedding.

\begin{definition}
A \emph{Teichm\"uller embedding} of dilatation $K\geq 1$ between finite Riemann surfaces $X$ and $Y$ is an injective continuous map $f:X \to Y$ for which there exist quadratic differentials $\vphi \in \qd^+(X)$ and $\psi \in \qd^+(Y)$ such that $f$ has the form $x+iy \mapsto Kx +iy$ in natural coordinates and such that $Y \setminus f(X)$ is a finite union of  points and horizontal arcs for $\psi$. We say that $\vphi$ and $\psi$ are \emph{initial} and \emph{terminal} quadratic differentials for $f$. A \emph{slit mapping} is a conformal Teichm\"uller embedding, i.e. one with $K=1$. 
\end{definition}

\begin{remark}
The horizontal arcs in the definition are allowed to overlap. In ge\-ne\-ral, the complement of the image of a Teichm\"uller embedding is an analytic graph (see Figure \ref{fig:backtrack} for example).
\end{remark}

\begin{remark} \label{rem:homeonotemb}
Despite the appellation, a Teichm\"uller homeomorphism between surfaces with non-empty ideal boundary is not necessarily a Teichm\"uller embedding, as its initial and terminal quadratic differentials are allowed be ne\-ga\-tive along the ideal boun\-da\-ry. 
\end{remark}

\begin{figure}[htp]
\centering
\includegraphics[scale=.8]{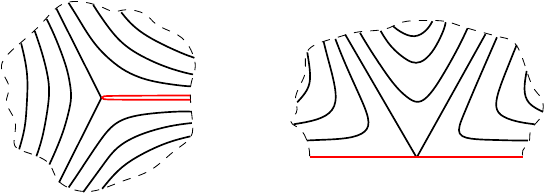}
\caption{Unfolding a slit ending at a simple zero (a $3$-prong singularity) of a quadratic differential yields half of a zero of order $4$ (a $6$-prong singularity).}
\label{fig:unfold}
\end{figure}

\begin{remark} If $Y$ is a finite Riemann surface, $\psi \in \qd^+(Y)$, and $X$ is the complement of a finite union of points and of horizontal arcs in $Y$, then $\iota^*\psi$ belongs to $\qd^+(X)$, where $\iota:X \to Y$ is the inclusion map. This means that $\iota^*\psi$ extends analytically to the ideal boundary of $X$ and that the latter is horizontal. Near the endpoint of a horizontal slit, one needs to take a square root in order to unfold the slit to an ideal boundary component. If we pull-back the quadratic differential $z^k \mathrm{d}z^2$ in $\CC$ by the squaring map $s(z)=z^2$ from $\HH$ to $\CC \setminus [0,\infty)$, we get the quadratic differential $4 z^{2k+2} \mathrm{d}z^2$. In other words, unfolding a slit ending at a singularity of order $k\geq -1$ yields a singularity of order $2k+2 \geq 0$ on the boundary. In particular, if the slit ends at a simple pole then the unfolded quadratic differential is regular at the corresponding point. It is perhaps more natural to count the number of prongs: an $n$-prong singularity transforms into half of a $2n$-prong singularity (see Figure \ref{fig:unfold}). The ideal boundary remains horizontal in the process of unfolding so that $\iota^*\psi \in \qd^+(X)$.

Every slit mapping $f:X \to Y$ arises in this way, in the sense that $f$ gives a conformal homeomorphism from $X$ to $f(X)$ and $f(X)\subset Y$ is obtained by removing finitely many points and horizontal arcs for some $\psi \in \qd^+(Y)$. The quadratic differential $\vphi = f^* \psi$ is redundant data.
\end{remark}

\begin{remark}
Every Teichm\"uller embedding can be factored as a Teichm\"uller homeomorphism $G:X \to Z$ (which happens to also be a Teichm\"uller embedding) followed by a slit mapping $F:Z \to Y$, where the terminal quadratic differential of $G$ agrees with the initial quadratic differential of $F$. 
\end{remark}

\begin{remark} \label{token}
If there is a Teichm\"uller embedding $X \to Y$, then $\qd^+(Y)$ is non-empty. This imposes some restrictions on $Y$, namely, it cannot be a sphere with at most $3$ punctures or a disk with at most $1$ puncture. By the same token, a Teichm\"uller embedding is always non-trivial and non-parabolic. 
\end{remark}

\begin{remark}
If $f$ is a Teichm\"uller embedding of dilatation $K$ with initial and terminal quadratic differentials $\vphi$ and $\psi$, then
$\dbar f / \partial f = k \overline{\vphi}/|\vphi|$ on $X$ and $\dbar (f^{-1}) / \partial (f^{-1}) = -k\overline{\psi} /|\psi|$ on $f(X)$, where $k=\frac{K-1}{K+1}$. This is one way of descri\-bing $f$ without using coordinates.
\end{remark}

We present Ioffe's theorem in two parts. The first part says that every extremal quasiconformal embedding is either conformal or a Teichm\"uller embedding.

\begin{theorem}[Extremal embeddings are Teichm\"uller] \label{thm:IoffeExistence}
Let $f : X \to Y$ be a quasiconformal embedding between finite Riemann surfaces with mini\-mal dilatation in its homotopy class. If $f$ is not conformal, then it is a Teichm\"uller embedding.
\end{theorem}

\begin{proof}
Suppose that $f$ is not conformal and let
$$
\mu := \begin{cases} \overline\partial (f^{-1}) / \partial (f^{-1}) & \text{on }f(X) \\ 0 & \text{on }Y \setminus f(X). \end{cases}
$$
Let $F : Y \to Y_\mu$ be the solution to the Beltrami equation $\overline \partial F / \partial F = \mu$ provided by Theorem \ref{thm:mrmt}.
By construction we have $\overline \partial (F\circ f) / \partial (F\circ f) = 0$ so that $F \circ f$ is a conformal embedding. By Teichm\"uller's theorem (Theorems \ref{thm:TeichExistence} and \ref{thm:TeichUniqueness}), there exists a quasiconformal homeomorphism $G: Y_\mu \to Y$ homotopic to $F^{-1}$ with mi\-ni\-mal dilatation. Moreover, $G$ is either conformal or a Teichm\"uller homeomorphism. The composition $G \circ F \circ f : X \to Y$ is a quasiconformal embedding homotopic to $f$, so that
$$
\Dil(f) \leq \Dil(G \circ F \circ f) \leq \Dil(G) \leq \Dil(F^{-1}) = \Dil(f).
$$
Thus all the terms in this chain are equal. The equality $\Dil(G) = \Dil(F^{-1})$ implies that $F$ has minimal dilatation in its homotopy class. Since $F$ is not conformal, it is a Teichm\"uller homeomorphism. This means that there is a non-zero $\psi\in \qd(Y)$ and a constant $k\in (0,1)$ such that $\mu = -k \overline \psi / |\psi|$ almost everywhere. In particular, $Y \setminus f(X)$ has measure zero and $f$ has constant pointwise dilatation.

Thought of as a homeomorphism from $X$ to $f(X)$, the map $f$ has minimal dilatation in its homotopy class and is thus a Teichm\"uller homeomorphism. Let $\vphi$ and $\omega$ be its initial and terminal quadratic differentials. Since $F \circ f$ is conformal, the directions of maximal stretching for $F$ and $f^{-1}$ must be perpendicular, which means that $\psi = c \omega$ on $f(X)$ for some $c>0$, which we may assume is equal to $1$ by rescaling.

We have to show that $\vphi \in \qd^+(X)$. If not, then $\vphi < 0$ along some segment $I \subset \partial X$. We will explicitly construct a quasiconformal embedding $\wtilde f$ from $X$ to $f(X)$ with pointwise dilatation smaller than $f$ near $I$. We may work in a natural coordinate chart for $\vphi$ in which $I$ is equal to the vertical segment $[-i,i]$ in the plane and $X$ is to the right of $I$. There is also a natural chart for $\omega$ in which $f$ takes the form $x+iy \mapsto Kx + i y$. Let $\Delta$ be the isoceles triangle with base $[-i,i]$ and apex $\delta > 0$. 
\begin{figure}[htp]  \centering
\includegraphics[scale=.7]{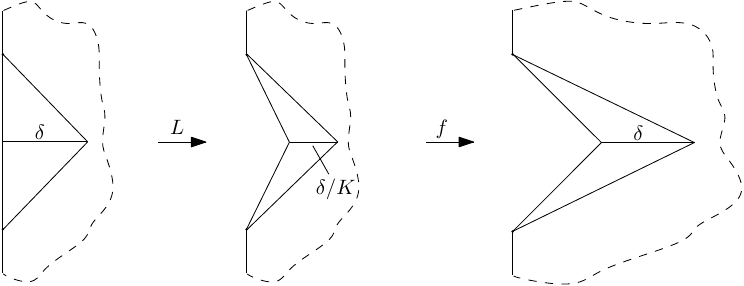}
\caption{Reducing the pointwise dilatation of an embedding near a vertical boundary arc.}
\label{fig:punchinwards}
\end{figure}
Consider the map $L : \Delta \to \Delta$ which is affine on the upper and lower halves of $\Delta$, fixes all three vertices of $\Delta$, and sends the midpoint of $I$ to $(1-1/K)\delta$. Extend $L$ to be the identity on $X \setminus \Delta$ and let $\wtilde f = f \circ L$. The linear part of $\wtilde f$ on the lower half of $\Delta$ is equal to $$\begin{pmatrix} 1 & (K-1)\delta \\ 0 & 1\end{pmatrix}$$ and the dilatation of this matrix tends to $1$ as $\delta \to 0$. A similar statement holds in the upper half of $\Delta$. Therefore, if $\delta$ is small enough, then the embedding $\wtilde f$ has strictly smaller pointwise dilatation than $f$ on $\Delta$. Moreover, the global dilatation of $\wtilde f$ is the same as $f$, so that $\wtilde f$ also has minimal dilatation in its homotopy class. By the first paragraph of the proof, the pointwise dilatation of $\wtilde f$ must be constant. This is a contradiction, and hence $\vphi \in \qd^+(X)$.

It remains to show that $f(X)$ is the complement of a graph which is horizontal with respect to $\psi$. Recall that $\what X$ is the compactification of $X$ obtained by adding its ideal boundary and filling its punctures. The metrics induced by $|\vphi|$ and $|\psi|$ extend to complete metrics on $\what X$ and $\what Y$. Since $f : X \to Y$ is $K$-Lipschitz with respect to these metrics, it extends to a $K$-Lipschitz map $\what f : \what X \to \what Y$. Moreover, $\what f$ is surjective since $Y \setminus f(X)$ has measure zero and hence empty interior. Let $I$ be the closure of a connected component of $\partial X \setminus \{\text{zeros of }\vphi\}$. There is a sequence $\{I_n\}$ of arcs in $X$ which are horizontal for $\vphi$ and converge uniformly to $I$. Since the image arcs $f(I_n)$ are all horizontal for $\psi$, they can only accumulate onto horizontal arcs, and thus $\what f(I)$ is horizontal. Therefore, the complement $\what Y \setminus f(X)=\what{f}(\partial X \cup \dot{X})$ is a union of finitely many points and horizontal arcs for $\psi$. In particular, the ideal boundary $\partial Y$ is horizontal for $\psi$ so that $\psi \in \qd^+(Y)$.
 \end{proof}

In the last paragraph of the proof we actually showed the following useful criterion.

\begin{lemma} \label{lem:fullmeasure}
Let $X$ and $Y$ be finite Riemann surfaces, let $\vphi \in \qd^+(X)$ and $\psi \in \qd(Y)\setminus \{0 \}$, and let $f: X \to Y$ be an embedding which is locally of the form $x+iy \mapsto Kx +iy$ in natural coordinates. If $f(X)$ is dense in $Y$, then $\psi \in \qd^+(Y)$ and $f$ is a Teichm\"uller embedding with respect to $\vphi$ and $\psi$.
\end{lemma}

The second part of Ioffe's theorem says that every Teichm\"uller embedding is extremal. The proof is very similar to the proofs of Teichm\"uller's uniqueness theorem given in \cite{Bers} and \cite[Chapter 11]{FarbMargalit}, only with additional technicalities due to the lack of compactness. 

\begin{theorem}[Teichm\"uller embeddings are extremal] \label{thm:IoffeUniqueness}
Let $f : X \to Y$ be a Teichm\"uller embedding of dilatation $K$ with initial and terminal quadratic differential $\vphi$ and $\psi$, and let $g: X \to Y$ be a $K$-quasiconformal embedding homotopic to $f$. Then $g$ is a Teichm\"uller embedding of dilatation $K$ with initial and terminal quadratic differentials $\vphi$ and $\psi$, and $g \circ f^{-1}: f(X) \to g(X)$ is locally a translation in natural coordinates for $\psi$.
\end{theorem}

\begin{proof}
We may assume that all the punctures of $X$ are essential with respect to $f$. Otherwise, the set $I$ of inessential punctures can be filled and $f$ extended to a Teichm\"uller embedding $\wtilde f : X \cup I \to Y$. Similarly, $g$ extends to a $K$-quasiconformal embedding $\wtilde g : X \cup I \to Y$ homotopic to $\wtilde f$. If $\wtilde g$ is a Teich\-m\"ul\-ler embedding with respect to $\psi$, then so is its restriction $g$.

By rescaling, we may assume that $\|\vphi\|=1$. This implies that $\|\psi\|=K$ since $f$ multiplies area by a factor $K$ and $f(X)$ has full measure in $Y$. Recall that $\what X$ and $\what Y$ are the metric completions of $X$ and $Y$ with respect to the distance induced by $|\vphi|$ and $|\psi|$. Since $f$ is $K$-Lipschitz with respect to these metrics, it extends to a $K$-Lipschitz map $\what f:\what X \to \what Y$. Since $f$ is assumed to have only essential punctures, $\what f$ maps $\dot{X}$ into $\dot{Y}$.

Let $M>1$. For every $n \in \NN$, let $ G_n : \what X \to \what X$ be a smooth $M$-quasicon\-for\-mal and $M$-Lipschitz embedding such that 
\begin{itemize}
\item $G_n(\dot{X})=\dot{X}$;
\item $G_n$ is homotopic to the identity rel $\dot{X}$;
\item $G_n(\what X)$ is contained in $X \cup \dot{X}$;
\item $G_n \to \id$ uniformly in the $C^1$ norm as $n \to \infty$.
\end{itemize}
Here is one way to construct such maps. Given a smooth vector field on $\what X$ pointing inwards on $\partial X$ and vanishing on $\dot{X}$, we can let $G_n$ be the corres\-pon\-ding flow at small enough time $t_n$. Then let $g_n = g\circ G_n : \what X \to \what Y$, which is a $KM$-quasiconformal embedding. The purpose of this construction is to circumvent the following difficulty: when $X \cup \dot{X}$ is non-compact $g$ has no reason to extend continuously to $\partial X$ while $g_n$ does.

Fix $n$ for a little while. By construction, $g_n|_X$ is homotopic to $g$ and hence to $f$. By Lemma \ref{lem:punctures}, there is a homotopy $H : X \times [0,1] \to Y$ from $f$ to $g_n|_X$ that extends to be constant at the punctures (which are assumed to be essential). By flowing $\partial X$ inward at the beginning and end of the homotopy as above, we may further assume that $H$ extends continuously to $\partial X \times [0,1]$. Let $\what H$ be the extension of $H$ to $\what X \times [0,1]$.

 For every $x \in X$, let $\ell(x)= \inf \int_\alpha \sqrt{|\psi|}$ where the infimum is taken over all piecewise smooth paths $\alpha : [0,1] \to Y$ that are homotopic to $t \mapsto H(x,t)$ rel endpoints. In general, the infimum need not be realized since the restriction of $|\psi|$ to $Y$ is not complete. However, $\ell(x)$ is equal to the length of the ``geodesic representative'' $\gamma_x$ of $t \mapsto H(x,t)$ in the completion $\what{Y}$ as explained in Subsection \ref{subsec:geom}. 
 
Since $H$ is continuous, the maps $x \mapsto \gamma_x$ and $x \mapsto \ell(x)$ are continuous. Moreover, they extend continuously to $\what X$. Indeed, for $x \in \what X \setminus X$ we can define $\gamma_x$ as the limit of $\gamma_{x_n}$ as $n \to \infty$, where $\{x_n\}_{n\geq 1}$ is a sequence in $X$ converging to $x$. This limit exists and does not depend on the sequence $\{x_n\}_{n \geq 1}$ since the path $t \mapsto \what H(x,t)$ is well-defined even for $x \in \what X \setminus X$. Moreover, its length $\ell(x)$ is the limit of the lengths $\ell(x_n)$. As $\what X$ is compact, there exists a constant $B$ such that $\ell(x)<B$ for every $x \in \what X$. This constant $B$ depends on $n$, but we will make it disappear before changing $n$.

Let $\eta$ be a horizontal arc of length $L>0$ in $X$. Since $f$ is a Teichm\"uller embedding of dilatation $K$, it sends $\eta$ to a horizontal arc of length $KL$ in $Y$. Let $x_0$ and $x_1$ be the endpoints of $\eta$. We can obtain a path homotopic to $f(\eta)$ in $Y$ by taking the concatenation of a piecewise smooth path $\alpha_0$ of length at most $B$ homotopic to $t \mapsto H(x_0,t)$, the image $g_n(\eta)$, and a piecewise smooth path $\alpha_1$ of length at most $B$ homotopic to $t \mapsto H(x_1,1-t)$. Since horizontal arcs minimize horizontal travel among all homotopic paths, we have
\begin{align*}
KL  = \int_{f(\eta)} |\re \sqrt{\psi}| & \leq \int_{\alpha_0} |\re \sqrt{\psi}|+\int_{g_n(\eta)} |\re \sqrt{\psi}|+\int_{\alpha_1} |\re \sqrt{\psi}| \\
& \leq  2B + \int_{g_n(\eta)} |\re \sqrt{\psi}|.
\end{align*}

Let $\mathrm{d}g_n$ denote the matrix of partial derivatives of $g_n$ with respect to na\-tu\-ral coordinates\footnote{The matrix is only defined up to sign, but no matter.} and $(\mathrm{d}g_n)_{1,1}$ its first entry. If $z=x+iy$ and $\zeta$ are natural coordinates for $\vphi$ and $\psi$, then $(\mathrm{d}g_n)_{1,1}=\re(\partial(\zeta \circ g_n\circ z^{-1}) / \partial x)$. If $g_n$ is absolutely continuous on $\eta$, then we have 
$$\int_\eta |(\mathrm{d}g_n)_{1,1}|\cdot \sqrt{|\vphi|} = \int_{g_n(\eta)} |\re \sqrt{\psi}| \geq KL - 2B.$$

Remove from $X$ all trajectories that go through a puncture of $X$ or a zero of $\vphi$ and denote the resulting full measure subset by $U$. For every $x \in U$, there is a unique (possibly closed) bi-infinite horizontal trajectory through $x$. For every $L>0$ and every $x \in U$, let $\eta(x,L)$ be the horizontal arc of length $L$ centered at $x$. Since $g_n$ is quasiconformal, it is absolutely continuous on almost every horizontal trajectory. Upon applying Fubini's theorem, we find
\begin{align*}
\int_U |(\mathrm{d}g_n)_{1,1}|\cdot|\vphi| &= \int_{x\in U} \left(\frac{1}{L} \int_{\eta(x,L)} |(\mathrm{d}g_n)_{1,1}|\cdot \sqrt{|\vphi|}  \right) \cdot |\vphi | \\
& \geq  \left(K  - \frac{2B}{L}\right) \int_U |\varphi|.
\end{align*}
Letting $L \to \infty$, we obtain $\int_U |(\mathrm{d}g_n)_{1,1}|\cdot|\vphi| \geq K \int_U |\varphi|$ and hence 
$$\int_X |(\mathrm{d}g_n)_{1,1}|\cdot|\vphi| \geq K \int_X |\varphi|=K,$$ since $U$ has full measure in $X$. 

Now that we got rid of the constant $B=B(n)$, we can vary the index $n$. We claim that $\int_X |(\mathrm{d}g)_{1,1}|\cdot|\vphi|=\lim_{n\to \infty} \int_X |(\mathrm{d}g_n)_{1,1}|\cdot|\vphi|$ and hence that $\int_X |(\mathrm{d}g)_{1,1}|\cdot|\vphi|\geq K$. This is a consequence of the Vitali convergence theorem \cite[p.94]{RoydenBook}. In order to apply the theorem, we need to check that the functions $|(\mathrm{d}g_n)_{1,1}|$ are uniformly integrable. First observe that $$\int_X \det(\mathrm{d} g)\cdot |\varphi|  = \int_{g(X)}|\psi| \leq K,$$ so that $\det(\mathrm{d} g)$ is integrable. It follows that for every $\eps>0$, there exists a $\delta>0$ such that if $A \subset X$ is measurable and $\int_A |\vphi| < \delta$, then $\int_A \det(\mathrm{d} g)\cdot |\varphi| < \eps$. Now if $\int_A |\vphi| < \delta / M^2$, then $\int_{G_n(A)} |\vphi|< \delta$ since $G_n$ is $M$-Lipschitz. By the Cauchy-Schwarz inequality we have
\begin{align*}
\left(\int_A |(\mathrm{d}g_n)_{1,1}|\cdot|\vphi|\right)^2 & \leq \int_A |(\mathrm{d}g_n)_{1,1}|^2\cdot |\vphi| \leq \int_A \|\mathrm{d}g_n\|^2 \cdot |\vphi| \\
& \leq KM \int_A \det(\mathrm{d}g_n) \cdot |\vphi| \\
& = KM \int_{A}  \det(\mathrm{d}_{G_n(z)}g)\det(\mathrm{d}_z G_n)\cdot |\vphi| \\
& = KM \int_{G_n(A)} \det(\mathrm{d} g)\cdot |\varphi| < KM \eps,
\end{align*}
which shows uniform integrability and proves the claim.

Applying Cauchy-Schwarz to the inequality $K \leq \int_X |(\mathrm{d}g)_{1,1}|\cdot|\vphi|$ yields
\begin{align*}
K^2 & \stackrel{\text(a)}{\leq} \left(\int_X |(\mathrm{d}g)_{1,1}|\cdot|\vphi|\right)^2  \stackrel{\text(b)}{\leq}   \int_X |(\mathrm{d}g)_{1,1}|^2\cdot|\vphi| \\
& \stackrel{\text(c)}{\leq}  \int_X \|\mathrm{d}g\|^2 \cdot |\varphi| \leq   K \int_X \det(\mathrm{d}g) \cdot |\varphi| \\
& =  K \int_{g(X)} |\psi| \stackrel{\text(d)}{\leq}  K \int_{Y} |\psi| =  K^2.
\end{align*}
Since the two ends of this chain of inequalities agree, each intermediate inequality is in fact an equation.
Equality in (b) implies that $|(\mathrm{d}g)_{1,1}|$ is equal to a constant almost everywhere on $X$, and that constant is equal to $K$ by (a). The inequality (c) is based on
$$
|(\mathrm{d}g)_{1,1}|^2 \leq |(\mathrm{d}g)_{1,1}|^2 + |(\mathrm{d}g)_{2,1}|^2 = |\mathrm{d}g\bigl(\begin{smallmatrix} 1 \\ 0 \end{smallmatrix}\bigr)|^2 \leq \sup_{\|\mathbf{v}\|=1} |\mathrm{d}g(\mathbf{v})| ^2 = \|\mathrm{d}g\|^2.
$$
Equality implies that $\mathrm{d}g\bigl(\begin{smallmatrix} 1 \\ 0 \end{smallmatrix}\bigr)=\pm \bigl(\begin{smallmatrix} K \\ 0 \end{smallmatrix}\bigr)$. Moreover, since $\mathrm{d}g$ stretches maximally in the horizontal direction which is preserved, $\mathrm{d}g$ must be diagonal, i.e. $\mathrm{d}g=\pm \bigl(\begin{smallmatrix} K & 0 \\ 0 & \ast \end{smallmatrix}\bigr)$ with $0< \ast \leq K$. Then the equality $K^2 = \|\mathrm{d}g\|^2 = K \det(\mathrm{d}g)$ determines that $\mathrm{d} g=\pm \bigl(\begin{smallmatrix} K & 0 \\ 0 & 1 \end{smallmatrix}\bigr)$ almost everywhere on $X$. 

Since $\mathrm{d} f=\bigl(\begin{smallmatrix} K & 0 \\ 0 & 1 \end{smallmatrix}\bigr)$ up to sign as well, we have $\mathrm{d}(g \circ f^{-1})=\pm \bigl(\begin{smallmatrix} 1 & 0 \\ 0 & 1 \end{smallmatrix}\bigr)$ almost everywhere on $f(X)$. The Beltrami coefficient of $g \circ f^{-1}$ is thus equal to $0$ almost everywhere on $f(X)$, so that $g \circ f^{-1}$ is conformal and in particular smooth. Since $f$ is smooth except at the zeros of $\varphi$, the same holds for $g$. Therefore the equality $\mathrm{d}g=\pm \bigl(\begin{smallmatrix} K & 0 \\ 0 & 1 \end{smallmatrix}\bigr)$ holds everywhere except at the zeros of $\vphi$, and $g$ takes the form $x + i y \mapsto \pm(Kx+iy)+c$ in natural coordinates. Equality in (d) means that $g(X)$ has full measure in $Y$. By Lemma \ref{lem:fullmeasure}, $g$ is a Teichm\"uller embedding with respect to $\vphi$ and $\psi$. Finally, the equality $\mathrm{d}(g \circ f^{-1})=\pm \bigl(\begin{smallmatrix} 1 & 0 \\ 0 & 1 \end{smallmatrix}\bigr)$ holds everywhere, so that $g\circ f^{-1}$ is a local translation. 
 \end{proof}

Theorems \ref{thm:IoffeExistence} and \ref{thm:IoffeUniqueness} together imply Theorem \ref{thm:IoffeIntro} from the introduction.

\begin{remark} \label{rem:termindiff}
Theorem \ref{thm:IoffeUniqueness} does not say that the quadratic differentials $\vphi$ and $\psi$ are unique up to scale because that is not the case when $K=1$. A slit mapping may be so with respect to a large-dimensional family of quadratic differentials. 

For example, suppose that $Y$ admits an anti-conformal involution $\sigma$ and let $X$ be the complement in $Y$ of finitely many arcs contained in the fixed locus of $\sigma$. Any element of $\qd^+(Y/\sigma)$ can be doubled to a quadratic differential in $\qd^+(Y)$ which is non-negative along the fixed locus of $\sigma$. The inclusion map $X \hookrightarrow Y$ is a slit mapping with respect to any quadratic differential obtained in this way.

However, if $K>1$ then $\vphi$ and $\psi$ are unique up to a positive scalar, since the Beltrami coefficients $\dbar f / \partial f$ and $\dbar (f^{-1}) / \partial (f^{-1})$ on $X$ and $f(X)$ encode the directions of maximal stretching, which are the horizontal and vertical directions of $\vphi$ and $\psi$ respectively. 
\end{remark}

\subsection{Lack of uniqueness} \label{subsec:nonunique}

If $f$ and $g$ are homotopic Teichm\"uller embeddings between finite Riemann surfaces $X$ and $Y$, then the inclusion map $f(X) \hookrightarrow Y$ and the composition $g\circ f^{-1} : f(X) \to Y$ are homotopic slit mappings with respect to the same quadratic differential by Theorem \ref{thm:IoffeUniqueness}. If $X$ and $Y$ are closed, then  $g\circ f^{-1}$ is a conformal automorphism of $Y$ homotopic to the identity. If $Y$ has genus at least $2$, then $g=f$ by Lemma \ref{lem:selfembisidentity}. 

However, the map $g\circ f^{-1}$ does not extend to all of $Y$ in general and indeed, Teichm\"uller embeddings are not necessaily unique in their homotopy class. There are two obvious ways for uniqueness to fail:
\begin{itemize}
\item if $Y$ is a torus, then we can post-compose $f$ with any automorphism of $Y$ isotopic to the identity;
\item if $f(X)$ is contained in an annulus $A\subset Y$, then we can post-compose $f$ with rotations of $A$. 
\end{itemize}

In \cite{Ioffe}, Ioffe claims these are the only exceptions, but this is wrong\footnote{The source of the mistake is \cite[Lemma 3.2]{Ioffe}. Similarly, \cite{Ioffe2}, \cite{EarleMarden}, and \cite{Golubev} contain minor errors as they build up on the false claim.}. The next simplest example is as follows. Let $Y$ be  a round annulus in the plane with a concentric circular arc removed, and let $X$ be the same annulus but with a slightly longer arc removed. Then we can obviously rotate $X$ inside of $Y$ by some amount. This gives a $1$-parameter family of slit mappings between triply connected domains. If the annulus is centered at the origin, then the quadratic differential in play is $-\mathrm{d} z^2 / z^2$.

One might think that every counterexample comes from a torus or an annulus with slits, but this is not the case. Here is a general method for constructing examples of slit mappings which are not unique in their homotopy class, with essentially any codomain $Y$. Let $\omega$ be a holomorphic $1$-form on $Y$ such that $\omega^2 \in \qd^+(Y)$. Then we can find finitely many horizontal arcs for $\omega^2$ such that their complement $X \subset Y$ is not rigid as follows. For every point $y \in \what Y$ which is either a zero of $\omega$ or a puncture of $Y$, and every trajectory $\gamma$ of $\omega$ ending at $y$ in forward time, remove a neighborhood of $y$ in $\gamma$ from $Y$ to obtain $X$. For all small enough $t > 0$, the time-$t$ flow for the vector field $1/\omega$ is well-defined on $X$, and is a slit mapping homotopic to the inclusion map yet different from it. The slit annulus example is a special case of this.

\begin{figure}[htp] 
\centering
\includegraphics[scale=.7]{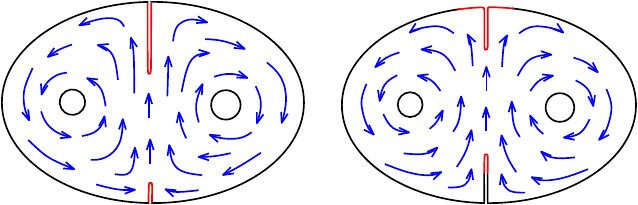}
\caption{A $1$-parameter family of slit mappings between pairs of pants.}
\label{fig:sliding}
\end{figure}

In the next section, we show that this kind of phenomenon is the worse that can happen: we can always get from any slit mapping to any homotopic one by flowing horizontally (except in the torus case), although the quadratic differential does not have to be globally the square of an abelian differential in general.

\section{Slit mappings are almost rigid} \label{sec:rigidity}

We prove that slit mappings are almost rigid in their isotopy class, in the sense that we can get from any one to any other by flowing horizontally (except in the torus case).

\begin{theorem} \label{thm:isotopythruslits}
Let $f_0,f_1: X \to Y$ be distinct homotopic slit mappings with terminal quadratic differential $\psi$. There is an isotopy $f_t$ from $f_0$ to $f_1$ through slit mappings with terminal quadratic differential $\psi$ such that for every $x \in X$, the path $t \mapsto f_t(x)$ is a regular geodesic for $\psi$ whose length and slope are independent of $x$. If $Y$ is not a torus, then the paths $t \mapsto f_t(x)$ are horizontal.
\end{theorem}

\begin{remark}
The statement is still true if we replace each occurrence of ``slit mappings'' with ``Teichm\"uller embeddings of dilatation $K$'', since this is really a statement about $f_1 \circ f_0^{-1}$, which is a slit mapping either way (Theorem \ref{thm:IoffeUniqueness}).
\end{remark}

\begin{proof}
The setup is the same as in Theorem \ref{thm:IoffeUniqueness}. We can assume that $f_0$ and $f_1$ have only essential punctures, since inessential punctures can be filled and play no special role whatsoever. As explained earlier, the maps $f_0$ and $f_1$ extend to continuous maps $\what f_0$ and $\what f_1$ from $\what X$ to $\what Y$. By Lemma \ref{lem:punctures}, there is a homotopy $H : X \times [0,1] \to Y$ from $f_0$ to $f_1$ which extends to be constant at the punctures. By pushing $\partial X$ inward at the beginning and end of the homotopy, we may further assume that $H$ extends continuously to $\partial X \times [0,1]$. To summarize, we have a homotopy $\what H : \what X \times [0,1] \to \what Y$ from $\what f_0$ to $\what f_1$ which sends $\dot{X} \times [0,1]$ into $\dot{Y}$ and whose restriction $H$ to $X \times [0,1]$ has image in $Y$.

Given $x \in X$, let $\ell(x)$ be the infimum $|\psi|$-length of piecewise smooth arcs homotopic rel endpoints to the path $t \mapsto H(x,t)$ in $Y$. This is equal to the length of the ``geodesic representative'' $\gamma_x$ of $t \mapsto H(x,t)$, which has image in $\what Y$ a priori. Note that $\gamma_x$ has endpoints $f_0(x)$ and $f_1(x)$ by construction. As explained in the proof of Theorem \ref{thm:IoffeUniqueness}, the maps $x \mapsto \gamma_x$ and $x \mapsto \ell(x)$ are continuous and extend continuously to $\what X$. 

Recall that the initial quadratic differential $\vphi = f_0^*\psi = f_1^*\psi$ is the same for both slit mappings and that $f_1 \circ f_0^{-1}$ is a local translation with respect to $\psi$ by Theorem \ref{thm:IoffeUniqueness}. In what follows, a \emph{singularity of $\psi$} is either a zero of $\psi$ or a puncture of $Y$. The rationale behind this is that even if $\psi$ is regular at some puncture, its pull-back to the universal cover has infinite angle singularities at the corresponding boundary points.

\begin{claim}
Suppose that $\ell$ has a local maximum at $x \in \what X$. Then $x \notin \dot{X}$ and $\ell$ is constant in a neighborhood of $x$. If $x \in X$, then $\gamma_x$ does not contain any singularity of $\psi$, and $\vphi$ does not have a zero at $x$. Moreover, for every $y$ near $x$ there is an isometric immersion from a Euclidean parallelogram $P$ to $Y$ which sends two opposite sides of $P$ to $\gamma_x$ and $\gamma_y$.
\end{claim}
\begin{proof}[Proof of Claim]
Suppose that $x \in \dot{X}$. Then $\ell(x)=0$ since $H$ is constant on $\dot{X}$. Since $x$ is a local maximum, $\ell(y) = 0$ for all $y$ in a neighborhood $U$ of $x$. This means that $\what f_0= \what f_1$ on $U$ and hence on all of $\what X$ by the identity principle. This contradicts the hypothesis that $f_0\neq f_1$.

We need to set up some notation. Let $\pi_X: \wtilde X \to X$ and $\pi_Y: \wtilde Y \to Y$ be the universal covers and let $\wtilde X'$ and $\wtilde Y '$ be their metric completions with respect to $\wtilde \vphi = \pi_X^* \vphi$ and $\wtilde \psi = \pi_Y^* \psi$. Recall that $\pi_X$ and $\pi_Y$ extend to continuous maps $ \pi_X': \wtilde X' \to \what X$ and $\pi_Y' : \wtilde Y' \to \what Y$. The homotopy $\what H:\what X \times [0,1] \to \what Y$ lifts to a continuous map $\wtilde H: \wtilde X' \times [0,1] \to \wtilde Y'$ under these. For any $(u,j) \in \wtilde X' \times \{0,1\}$, let us write $\wtilde H_j(u) := \wtilde H(u,j)$. Note that $\wtilde H_j$ is a lift of $\what f_j$. In particular, the restriction of $\wtilde H_j$ to $\wtilde X$ is a local translation with respect to $\wtilde\vphi$ and $\wtilde \psi$. For any $u \in \wtilde X'$, let $\wtilde \gamma_u$ be the geodesic between $\wtilde H_0(u)$ and $\wtilde H_1(u)$ in $\wtilde Y'$ and let $\wtilde \ell(u)$ be its length. By definition, $\gamma_{\pi_X'(u)} =  \pi_X' ( \wtilde \gamma_u)$ and $\ell(\pi_X'(u))=\wtilde \ell(u)$. Let $u \in \wtilde X'$ be any lift of $x$ under $\pi_X'$. Then $\wtilde \ell$ has a local maximum at $u$.  Let $B$ be a ball centered at $u$ in $\wtilde X'$ such that $\wtilde \ell(v) \leq \wtilde \ell(u)$ for every $v \in B$. Let $\eta : [-1,1] \to B$ be a geodesic such that $\eta(0)=u$. Let $\eta^-$ and $\eta^+$ be its restrictions to $[-1,0]$ and $[0,1]$ respectively, and let $v^\pm =\eta(\pm 1)$. Given two geodesic rays $\alpha$ and $\beta$ sharing an endpoint in $\wtilde X'$ or $\wtilde Y'$, we will denote their Alexandrov angle by $\angle(\alpha,\beta)$ instead of the three point notation from Subsection \ref{subsec:geom} (this angle takes values in $[0,\pi]$).

After possibly shrinking $B$ and $\eta$, the images $\wtilde H_0(\eta^+)$ and $\wtilde H_1(\eta^+)$ are geodesic segments of the same length in $\wtilde Y '$ (and similarly for $\eta^-$). This is because $\wtilde H_0$ and $\wtilde H_1$ are local isometries in $\wtilde X$. Their extensions to $(\pi_X')^{-1}(\partial X)$ are also local isometries except perhaps at zeros of $\vphi$ where they can fold the boundary in two. In any case, $\wtilde H_0$ and $\wtilde H_1$ are isometries along any sufficiently short geodesic rays at $u$.

Here is the \myhypertarget{main observation}{main observation}. Let $\sigma \in \{+, -\}$. By Lemma \ref{lem:quad}, if $$\angle (\wtilde H_0(\eta^\sigma),\wtilde \gamma_u) + \angle (\wtilde \gamma_u, \wtilde H_1(\eta^\sigma))\geq \pi$$ then $\wtilde \ell(v^\sigma) \geq \wtilde \ell(u)$. But $\wtilde \ell(v^\sigma) \leq \wtilde \ell(u)$ by hypothesis. Since equality holds, the geodesic quadrilateral $Q^\sigma$ with sides $\wtilde \gamma_u$, $\wtilde\gamma_{v^\sigma}$, $\wtilde H_0(\eta^\sigma)$ and $\wtilde H_1(\eta^\sigma)$ is isometric to a Euclidean parallelogram and $\angle (\wtilde H_0(\eta^\sigma),\wtilde \gamma_u) + \angle (\wtilde \gamma_u, \wtilde H_1(\eta^\sigma))= \pi$.

\paragraph{Case 1:} Suppose that $x \in X$. 

Then after possibly shrinking the ball $B$ centered at $u\in \wtilde X$, the restrictions of $\wtilde H_0$ and $\wtilde H_1$ to $B$ are isometries. This means that $\wtilde H_0(\eta)$ and $\wtilde H_1(\eta)$ are geodesics and hence that
\begin{equation*}
\pi = \angle(\wtilde H_j(\eta^+),\wtilde H_j(\eta^-)) \leq \angle (\wtilde H_j(\eta^+),\wtilde \gamma_u) + \angle( \wtilde \gamma_u,\wtilde H_j(\eta^-))
\end{equation*}
for $j=0$ and $j=1$ by subadditivity of angles. Thus if $$\angle (\wtilde H_0(\eta^+),\wtilde \gamma_u) + \angle (\wtilde \gamma_u, \wtilde H_1(\eta^+)) < \pi$$ then $$\angle (\wtilde H_0(\eta^-),\wtilde \gamma_u) + \angle (\wtilde \gamma_u, \wtilde H_1(\eta^-))> \pi.$$ Lemma \ref{lem:quad} implies that $\wtilde \ell(v^-)> \wtilde \ell (u)$, which is a contradiction. We conclude that $\angle (\wtilde H_0(\eta^+),\wtilde \gamma_u) + \angle (\wtilde \gamma_u, \wtilde H_1(\eta^+))\geq \pi$. By the \myhyperlink{main observation}, $Q^+$ is isometric to a Euclidean parallelogram and $\wtilde \ell (v^+) = \wtilde \ell(u)$. Since the geodesic $\eta$ through $u$ was arbitrary, $\wtilde \ell$ is constant in a neighborhood of $u$.

Choose $\eta$ in such a way that each of $H_0(\eta^+)$, $H_0(\eta^-)$, $H_1(\eta^+)$ and $H_1(\eta^-)$ forms a positive angle with $\wtilde \gamma_u$. Then the corresponding parallelograms $Q^+$ and $Q^-$ are non-degenerate and lie on opposite sides of $\wtilde \gamma_u$. In particular, there is no excess angle on either side of $\wtilde \gamma_u$ in the metric $|\wtilde \psi|$. In other words, the interior of $\wtilde \gamma_u$ does not contain any singularity of $\wtilde \psi$. Suppose however that $\wtilde \psi$ has a zero at an endpoint of $\wtilde \gamma_u$. Then $\wtilde \vphi$ has a zero at $u$ as both $\wtilde H_0$ and $\wtilde H_1$ are local translations in a neigborhood of $u$. Since $u$ is in the interior, the cone angle at $u$ in the metric $|\wtilde \vphi|$ is at least $3\pi$. Thus there is a whole sector of points $v^+ \in B \setminus \{u\}$ such that $\angle (\wtilde H_0(\eta^+) , \wtilde \gamma_u) = \pi$ where $\eta^+$ is the geodesic $[u,v^+]$. For each such $v^+$, we then have $\angle (\wtilde H_1(\eta^+) , \wtilde \gamma_u) = 0$, which means that $\wtilde H_1(\eta^+)$ and $\wtilde \gamma_u$ share a segment. In other words, $\wtilde H_1$ collapses a whole sector of $B$ into $\wtilde \gamma_u$. But this is impossible since $\wtilde H_1$ is locally injective. This shows that $\wtilde \gamma_u$ is completely free of singularities. In particular, $\wtilde \vphi=(\wtilde H_0)^* \wtilde \psi$ does not have a zero at $u$.
 
Since $\wtilde \gamma_u$ is regular and has endpoints in $\wtilde Y$, it is contained in $\wtilde Y$. Similarly, for any short geodesic ray $\eta^+$ at $u$, the resulting parallelogram $Q^+$ is contained in $\wtilde Y$. The covering map $\pi_Y : \wtilde Y \to Y$ is an isometric immersion which sends two parallel sides of $Q^+$ to $\gamma_x$ and $\gamma_y$ where $y=\pi_Y(v^+)$.
 
\paragraph{Case 2:} Suppose that $x \in \partial X$ and that $\vphi$ has a zero at $x$.

Then the total cone angle at $u$ in the metric $|\wtilde \vphi|$ is at least $2\pi$. For any geodesic $\eta$ through $u$ in $B$, the \myhyperlink{main observation} still holds.  Now for every geodesic ray $\eta^+$ from $u$, there exists a ray $\eta^-$ such that $\wtilde H_0(\eta)$ and $\wtilde H_1(\eta)$ are geodesic, where $\eta = \eta^+ \cup \eta^-$. This is because there is enough angle at $u$ to make sure that the cone angle on either side of $\wtilde H_j(\eta)$ at $\wtilde H_j(u)$ is at least $\pi$.

The argument from Case 1 applies for any such geodesic $\eta$: if the angles at the endpoints of $\wtilde \gamma_u$ in $Q^+$ sum to less than $\pi$, then $\wtilde \ell(v^-)> \wtilde \ell (u)$. Hence the angles sum to at least $\pi$ and $Q^+$ is a Euclidean parallelogram by the \myhyperlink{main observation}. This shows that $\wtilde \ell(v^+) = \wtilde \ell(u)$. This holds for every $v^+ \in B \setminus \{u\}$, so that $\wtilde \ell$ is constant there.

\paragraph{Case 3:} Suppose that $x \in \partial X$ and that $\vphi$ does not have a zero at $x$.

Then the cone angle at $u$ in the metric $|\wtilde \vphi|$ is $\pi$ and there is only one geodesic $\eta$ through $u$, namely the horizontal one running along the boundary. In this case, neither $\wtilde H_0$ nor $\wtilde H_1$ folds the boundary of $\wtilde X'$ in two at $u$. Indeed, if $\wtilde H_j$ did fold in two at $u$, then $\wtilde \psi$ would have a simple pole at $\wtilde H_j(u)$. But all the simple poles of $\psi$ in $\what Y$ have been unwrapped to infinite angle singularities in $\wtilde Y'$. 

Therefore $\wtilde H_0(\eta)$ and $\wtilde H_1(\eta)$ are geodesic. Thus the same argument as above applies to show that the corresponding quadrilaterals $Q\pm$ are parallelograms. If one of $Q^+$ or $Q^-$ is degenerate, then so is the other one. 

Suppose first that $Q^\pm$ are degenerate so that the beginning and end of $\wtilde \gamma_u$ are contained in $\wtilde H_0(\eta)$ and $\wtilde H_1(\eta)$ respectively. Moreover if the beginning is contained in $\wtilde H_0(\eta^+)$ then the end is contained in $\wtilde H_1(\eta^-)$, and vice versa. It follows that for any $w \in B \setminus \{u\}$ the geodesic $\tau =[u,w]$ satifies 
\begin{equation} \label{eq:sumpi}
\angle (\wtilde H_0(\tau), \wtilde \gamma_u)+\angle (\wtilde \gamma_u, \wtilde H_1(\tau))= \pi.
\end{equation}
By the \myhyperlink{main observation}, $\wtilde \gamma_w$ and $\wtilde \gamma_u$ are the opposite sides of a Euclidean parallelogram and $\wtilde \ell(w) = \wtilde\ell(u)$.

Now suppose that $Q^+$ and $Q^-$ are non-degenerate. Then the interior of $\wtilde \gamma_u$ is free of singularities, as $Q^+$ and $Q^-$ lie on opposite sides of it. This implies that equation \eqref{eq:sumpi} holds for any geodesic ray $\tau$ at $u$ in $B$, as parallel transport is well-defined along $\wtilde \gamma_u$. Therefore $\wtilde \ell$ is constant on $B$.
 \end{proof}

It follows from the above claim that $\ell$ is constant. Indeed, the subset of $\what X$ where $\ell$ attains its maximum is open by the claim and is closed by continui\-ty, hence equal to $\what X$. In particular, every point of $\what X$ is a local maximum of $\ell$, so the additional conclusions of the claim hold everywhere. Specifically,
\begin{enumerate}[label=(C\arabic*)]
 \item \label{c1} $\dot{X}$ is empty\footnote{Recall that punctures are assumed to be essential in the proof. This proves that if two slit mappings are homotopic and have an essential puncture, then they are equal. \label{nopunctures}};
 \item \label{c5} $\vphi$ does not have any zero in $X$; 
 \item \label{c3} for every $x \in X$ and $y$ near $x$, the geodesics $\gamma_x$ and $\gamma_y$ form opposite sides of an isometrically immersed parallelogram in $Y$.
 
\end{enumerate} 

Let $x \in X$. Parametrize $\gamma_x : [0,1] \to Y$ proportionally to arc length in such a way that $\gamma_x(0) = f_0(x)$ and $\gamma_x(1)=f_1(x)$. For every $(x,t) \in X \times [0,1]$, define $f_t(x) = \gamma_x(t)$. Then $f_t$ has image in $Y$ and depends continuously on $t$. By \ref{c3}, $f_t\circ f_0^{-1}$ is a local translation with respect to $\psi$ for every $t \in [0,1]$. Thus $f_t$ is a local translation with respect to $\vphi$ and $\psi$. It remains to prove that $f_t$ is injective and that its image is the complement of finitely many horizontal arcs and points.

Suppose that $Y$ is a torus. Then $\psi = \omega^2$ for an abelian diffe\-ren\-tial $\omega$. For any $x \in X$, the tangent vector $v(x,t)=\omega(\gamma_x'(t)) \in \CC$ does not depend on $t\in [0,1]$ since $\gamma_x$ is a smooth geodesic. As nearby geodesics $\gamma_x$ and $\gamma_y$ are parallel and of the same length, $v$ is locally constant as a function of $x \in X$, hence constant on $X$. This means that $f_t$ differs from $f_0$ by the translation $z \mapsto z + t v$ with respect to $\omega$, and is therefore a slit mapping.

Suppose that $Y$ is an annulus. Then $\psi$ is again the square of an abelian differential $\omega$. By the same argument as in the previous paragraph, $v=\omega(\gamma_x'(t))$ is independent of $(x,t) \in X \times [0,1]$. If the geodesics $\gamma_x$ are not horizontal, then they point away from $Y$ on one of the two boundary components, which is absurd. Hence $v\in \RR$ and $f_t$ differs from $f_0$ by the horizontal translation $z \mapsto z + t v$ with respect to $\omega$ (i.e., by a rotation of $Y$).

Our proof that $f_t$ is a slit mapping in the general case is rather indirect. By Lemma \ref{lem:fullmeasure}, it suffices to show that $f_t$ is injective and that $f_t(X)$ is dense in $Y$. To prove this, we first show that the geodesics $\gamma_x$ are horizontal.

\begin{claim} 
If $Y$ is not a torus, then $\gamma_x$ is horizontal for every $x\in \what X$.
\end{claim}
\begin{proof}
By the same argument as in the torus and annulus case, the unoriented slope $\psi(\gamma_x'(t))$ is constant on $X \times [0,1]$. Since the map $x \mapsto \gamma_x$ is continuous on $\what X$ and since the limit of a sequence of geodesics of slope $s$ has slope $s$ wherever it is smooth, the slope function is constant on all $\what X$. Thus if we find a single horizontal geodesic, then all geodesics are horizontal.

If $Y$ is an annulus, then we have already shown that $\gamma_x$ is horizontal for every $x\in \what X$. We can leverage the same idea to prove that if $Y$ contains a horizontal cylinder, then $\gamma_x$ is horizontal for some (hence all) $x \in \what X$.

Let $A \subset Y$ be a maximal horizontal cylinder for $\psi$.  Let $\alpha$ be a closed horizontal trajectory in $A \cap f_0(X)$. Then $f_t\circ f_0^{-1}(\alpha)$ is a closed horizontal trajectory homotopic to $\alpha$---hence is contained in $A$---for every $t\in [0,1]$. Since all but finitely many horizontal trajectories in $A$ are contained in $f_0(X)$, we have that $f_t\circ f_0^{-1}(A \cap f_0(X)) \subset A$ for every $t \in [0,1]$. In particular, the geodesic $\gamma_x$ is contained in $A$ for every $x \in f_0^{-1}(A)$.

Observe that $f_0(X) \cap A$ has either one or two connected components. If $f_0(X) \cap A$ is connected, then the same argument as in the annulus case applies: if the geodesics $\gamma_x$ are not horizontal, then they point away from $A$ on one of its two boundary components, contradicting the fact that they are contained in $A$. Suppose that $f_0(X) \cap A$ has two components. Since $X$ is connected, $f_0(X)$ intersects each boundary component of $A$. If $\gamma_x$ points strictly inward at such a boundary point $f_0(x)$, then for small $t>0$ the map $f_t\circ f_0^{-1}$ drags points from the complement of $A$ into $A$, which is impossible. Said differently, if $t>0$ and $\alpha \subset A \cap f_t(X)$ is a closed horizontal trajectory sufficiently close to $f_0(x)$, then $f_0 \circ f_t^{-1}(\alpha)$ is a closed horizontal trajectory outside of $A$ yet homotopic to $\alpha$, contradicting the maximality of $A$. In either case, the geodesics $\gamma_x$ have to be horizontal.  

If $Y$ has a horizontal cylinder then we are done, so suppose that it does not have any. Let $t \in [0,1]$. Since $f_t:X \to Y$ is $1$-Lipschitz with respect to $|\vphi|$ and $|\psi|$, it extends to a $1$-Lipschitz map $\what f_t : \what X \to \what Y$ between the completions. Let $C$ be a connected component of $\partial X$, let $N$ be neighborhood of $C$ in $\what X$ such that $N\cap X$ is an annulus, and let $\alpha \subset N\cap X$ be homotopic to $C$ in $N$. If $f_0(\alpha)$ is not homotopic to a point or a puncture in $Y$, then $\what f_t(C)$ becomes the ``geodesic representative'' of $f_0(\alpha)$ after erasing its backtracks (see Figure \ref{fig:backtrack}). Since $Y$ contains no horizontal cylinders, the ``geodesic representative'' of $f_0(\alpha)$ is unique. This means that $\what f_t(C)$ is contained in the same horizontal leaf for every $t \in [0,1]$, hence that $\gamma_x$ is horizontal for every $x \in C$.

\begin{figure}[htp] 
\centering
\includegraphics[scale=.7]{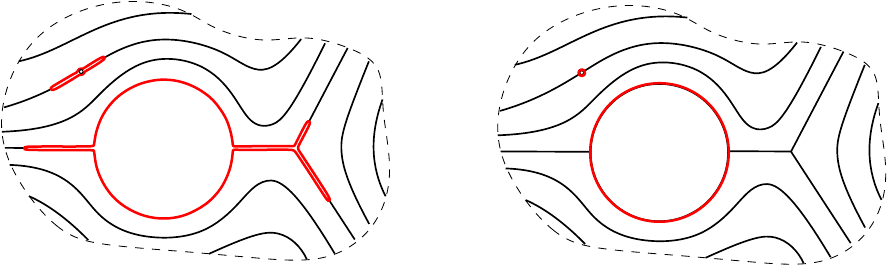}
\caption{Pulling a horizontal curve tight shrinks it to a point, a puncture, or a geodesic.}
\label{fig:backtrack}
\end{figure}

We may thus assume that $f_0(\alpha)$ is homotopic to a point or a puncture in $Y$ for every $C$, $N$ and $\alpha$ as above. This means that $\what f_t(C)$ is completely folded onto itself, as it reduces to a point after erasing its backtracks. Thus $\what f_t(N)$ forms an open neighborhood of $\what f_t(C)$ in $\what Y$. This holds for every component $C$ of $\partial X$ and every $t \in [0,1]$. Moreover, $\partial Y$ is empty for otherwise there is some component $C$ of $\partial X$ such that $\what f_0(C)$ intersects $\partial Y$, forcing $f_0(\alpha)$ to be essential.

Since $Y$ is not a torus nor an annulus, $\psi$ has at least one singularity, say at $y \in \what Y$. As $\what f_0 : \what X \to \what Y$ is surjective, there is some $x\in \what X$ such that $\what f_0 (x) = y$. We have that $x \in \partial X$ since $\dot{X}$ is empty \ref{c1} and $\vphi$ does not have any zeros in $X$ \ref{c5}. Let $C$ be the component of $\partial X$ containing $x$. If $\gamma_x$ is not horizontal, then $\what f_t(C)$ moves off of its horizontal leaf for small enough $t>0$. On the other hand, $y$ is contained in the open set $\what f_t(N)$ for all sufficiently small $t>0$. But then $y \in f_t(N \setminus C) \subset f_t(X) \subset Y$ so that $\psi$ has a zero at $y$ and hence $\vphi$ has zeros on $f_t^{-1}(y)$, contradicting \ref{c5}. Therefore $\gamma_x$ is horizontal.
 \end{proof}

Since $f_t$ is a local translation with respect to $\vphi$ and $\psi$, the area of its image is at most the area of $X$, with equality only if $f_t$ is injective. Recall that the set $U\subset X$ of points through which passes a bi-infinite regular horizontal trajectory of $\varphi$ has full measure in $X$. For every bi-infinite regular horizontal trajectory $\eta \subset U$, its image $f_0(\eta)$ is a bi-infinite regular horizontal trajectory of $\psi$. For any $x \in \eta$, the horizontal geodesic $\gamma_x$ is thus contained in $f_0(\eta)$. Hence $f_t\circ f_0^{-1}$ acts as a translation along $f_0(\eta)$ for every $t \in [0,1]$, so that $f_t(\eta) = f_0(\eta)$. It follows that $f_t(U)=f_0(U)$. But $f_0(U)$ has full measure in $Y$ since $U$ has full measure in $X$ and $f_0(X)$ has full measure in $Y$. Therefore $f_t(X)$ has full measure in $Y$, whose area is equal to that of $X$. Hence $f_t$ is injective. 

It follows that $f_t$ is a slit mapping by Lemma \ref{lem:fullmeasure}. Thus $(x,t) \mapsto f_t(x)$ is an isotopy from $f_0$ to $f_1$ through slit mappings. If $Y$ is not a torus, then the path $t \mapsto f_t(x)=\gamma_x(t)$ is horizontal for every $x\in X$ by the previous claim.
 
\end{proof}

\begin{remark} \label{rem:nopunctures2}
As observed in footnote \ref{nopunctures}, if $f_0,f_1:X \to Y$ are homotopic slit mappings with an essential puncture, then $f_0 = f_1$. Similarly, if their initial quadratic differential has a zero in $X$ then $f_0=f_1$. As another example, if their terminal quadratic differential has a horizontal saddle connection whose interior is contained in $f_j(X)$ for $j=0$ or $j=1$ then $f_0 = f_1$. This is because the saddle connection has nowhere to go under the horizontal isotopy. The same statements apply to Teichm\"uller embeddings.
\end{remark}

\begin{remark}
If $Y$ is not the torus and $f_0,f_1:X \to Y$ are distinct homotopic slit mappings, then their terminal quadratic differential $\psi$ is unique up to a positive scalar (cf. Remark \ref{rem:termindiff}). Indeed, for every $x\in X$ and every $\psi \in \qd^+(Y)$ with respect to which $f_0$ and $f_1$ are slit mappings, $f_0(x)$ and $f_1(x)$ are on the same horizontal leaf of $\psi$. Thus if $f_t$ is \emph{any} isotopy between $f_0$ and $f_1$ through slit mappings, then $t \mapsto f_t(x)$ is horizontal for every $x\in X$ and every $\psi\in \qd^+(Y)$ with respect to which $f_0$ and $f_1$ are slit mappings. But if two quadratic differentials have the same horizontal trajectories, then they are positive multiples of each other. The contrapositive is worth mentioning explicitly: if the terminal quadratic differential of a slit mapping $f$ is not unique up to scale, then $f$ is the unique conformal embedding in its homotopy class. For example, if $\dim \qd^+(Y)>1$ and $X \subset Y$ has finite complement, then the inclusion map $X \hookrightarrow Y$ is the only conformal embedding in its homotopy class (cf. Lemma \ref{lem:emptyidealbdry}).
\end{remark}

Let $h:X \to Y$ be a topological embedding between finite Riemann surfaces. Recall that the space $\emb(X,Y,h)$ is the set of all conformal embeddings homotopic to $h$ equipped with the compact-open topology. We now deduce Theorem \ref{thm:AlmostRigidity} from the introduction, which describes the homeomorphism type of $\emb(X,Y,h)$ whenever the latter contains a slit mapping.

\begin{proof}[Proof of Theorem \ref{thm:AlmostRigidity}]
If $\emb(X,Y,h)$ contains a slit mapping with terminal differential $\psi$, then every $f \in \emb(X,Y,h)$ is a slit mapping with respect to $\psi$ by Theorem \ref{thm:IoffeUniqueness}. Let $x \in X$. The evaluation map $\emb(X,Y,h) \to Y$ sending $f$ to $f(x)$ is continuous by definition of the compact-open topology. We now show that it is injective. 

Suppose that $f_0(x)=f_1(x)$ for distinct $f_0,f_1 \in \emb(X,Y,h)$ and some $x \in X$. Our goal is to show that $f_0 = f_1$. If not, let $f_t$ be the isotopy from $f_0$ to $f_1$ provided by Theorem \ref{thm:isotopythruslits}. Then $\gamma: S^1 \to Y$ defined by $\gamma(t)=f_t(x)$ is a closed regular geodesic for $\psi$.

Suppose that $Y$ is a torus. For each $t \in [0,1]$, there is a unique translation $g_t \in \aut_0(Y)$ such that $g_t(\gamma(t))=\gamma(0)$. Then $t \mapsto g_{t} \circ f_t$ is an isotopy from $f_0$ to $f_1$ which is constant at $x$. Thus the restrictions $f_0^\star, f_1^\star : X \setminus \{x\} \to Y \setminus \{f_0(x)\}$ are homotopic slit mappings with an essential puncture. It follows that $f_0=f_1$ by Remark \ref{rem:nopunctures2}. Alternatively, the proof of Theorem \ref{thm:isotopythruslits} show that $f_0$ and $f_1$ differ by an element of $\aut_0(Y)$, which has to be the identity since it fixes the point $f_0(x)=f_1(x)$.

Suppose that $Y$ is not a torus. Then $\gamma$ is horizontal. Let $A \subset Y$ be a cylinder foliated by closed horizontal trajectories homotopic to $\gamma$ which is symmetric about $\gamma$. Since $f_0(x)$ is a fixed point of the local translation $f_1 \circ f_0^{-1}$, we have that $f_1 \circ f_0^{-1}$ is either equal to the identity near $f_0(x)$ or a half-turn around $f_0(x)$. If $f_1 \circ f_0^{-1}$ is equal to the identity near $f_0(x)$, then $f_0=f_1$ by the identity principle. Otherwise, let $\sigma$ be the conformal involution of $A$ which fixes $f_0(x)$ and permutes its two boundary components. By the identity principle $f_1 \circ f_0^{-1} = \sigma$ on the connected component $U$ of $A\cap f_0(X)$ containing $f_0(x)$. Let $\alpha \subset U$ be a closed horizontal trajectory and let $\beta = f_0^{-1}(\alpha)$. Then $f_0(\beta)=\alpha$ and $f_1(\beta)=\sigma(\alpha)$ are not homotopic in $Y$, for they have reverse orientations. This is a contradiction, which proves that the first situation prevails. 

Thus the evaluation map is injective. Since $\emb(X,Y,h)$ contains a slit mapping, $h$ is non-trivial and non-parabolic by Remark \ref{token}. It follows that $\emb(X,Y,h)$ is compact by Lemma \ref{lem:compactness}. Any injective continuous map from a compact space to a Hausdorff space is a homeomorphism onto its image. Let $V(x)$ be the image of evaluation map at $x\in X$. If $Y$ is a torus, then $V(x)=Y$ since $\aut_0(Y)$ acts on $\emb(X,Y,h)$ by composition on the left and acts transitively on $Y$. 

Suppose that $Y$ is not a torus. Then $V(x)$ is path-connected and is contained in a regular horizontal trajectory by Theorem \ref{thm:isotopythruslits}. As $V(x)$ is compact, it is either a closed trajectory, a geodesic segment, or a point. If $h$ is cyclic, then the image of any $f \in \emb(X,Y,h)$ is contained in some annulus $A \subset Y$. Since $f$ can be post-composed with rotations of $A$, the set $V(x)$ contains a circle, hence is equal to a circle. 

Conversely, if $V(x)$ is a horizontal circle then its free homotopy class in $Y$ does not depend on $x\in X$ since $X$ is path-connected and $V(x)$ depends continuously on $x$. Therefore $\bigcup_{x\in X} V(x)$ is contained in a horizontal cylinder $A\subset Y$, which means that $f(X) \subset A$ for every  $f \in \emb(X,Y,h)$ so that $h$ is cyclic. This shows that if $h$ is generic (so that $Y$ is not a torus), then $V(x)$ is either a segment or a point.
 \end{proof}

As pointed out in the introduction, if $\emb(X,Y,h)$ contains a slit mapping  then the conclusions of Theorem \ref{thm:defret} follow immediately from Theorem \ref{thm:AlmostRigidity}. In that case $h$ is necessarily non-trivial. If $Y$ is a torus then $\emb(X,Y,h)$ is homeomorphic to a torus, hence homotopy equivalent to one. If $Y$ is not a torus and $h$ is cyclic, then $\emb(X,Y,h)$ is homeomorphic to a circle. If $h$ is generic, then $\emb(X,Y,h)$ is homeomorphic to a point or a compact interval, hence contractible.

\section{The modulus of extension} \label{sec:localmax}

In this section, we characterize local maxima of the modulus of extension (Theorem \ref{thm:localmaxintro}) and use that to prove that the space of all conformal embeddings in a given homotopy class is connected under some conditions (Theorem \ref{thm:connected}).

Let $h:X \to Y$ be a topological embedding between finite Riemann surfaces where $X$ has non-empty ideal boundary. We will impose further conditions on $h$ in a moment, but for now the only hypothesis is that $\partial X \neq \varnothing$.

\begin{figure}[htp]  \centering
\includegraphics[scale=.8]{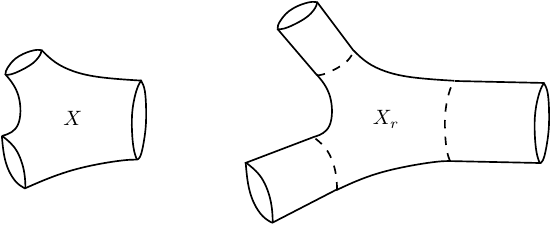}
\caption{The surface $X_r$ is obtained by gluing a cylinder of modulus $r$ to each ideal boundary component of $X$. }
\label{fig:enlargement}
\end{figure}

 For each connected component $C$ of $\partial X$, choose an analytic para\-me\-tri\-zation $\zeta_C :S^1 \to C$. For every $r \in (0,\infty]$ and every component $C$ of $\partial X$, glue a copy of the cylinder $S^1 \times [0,r)$ to $X\cup \partial X$ along $S^1 \times \{0\}$ using the map $\zeta_C$ (see Figure \ref{fig:enlargement}). We denote the resulting surface by $X_r$ and also let $X_0=X$. If $\rho \leq r$, then the inclusion $[0,\rho) \subset [0,r)$ induces a conformal embedding $X_{\rho} \subset X_{r}$. We call the directed family of Riemann surfaces $\{X_r\}$ a \emph{$1$-parameter family of enlargements of $X$}. Note that for every $r\in[0,\infty]$, there is a homeomorphism $H_r : X_r \to X$ which when followed by the inclusion $X \subset X_r$ is homotopic to the identity. We will abuse notation and write $\emb(X_r,Y,h)$ instead of $\emb(X_r,Y,h\circ H_r)$.

 Given a conformal embedding $f \in \emb(X,Y,h)$, we define the \emph{modulus of extension of $f$} as $$\m(f)=\sup \{\, r\in[0,\infty] \mid f \text{ extends to a conformal embedding }X_r \hookrightarrow Y \,\}.$$
This number depends on the $1$-parameter family of enlargements $\{X_r\}_{r \in [0,\infty]}$ which is fixed once and for all. It is easy to see that the supremum is realized and that the conformal extension of $f$ to $X_{\m(f)}$ is unique.

\begin{lemma}
For every $f \in \emb(X,Y,h)$, there is a unique conformal embedding $\what f : X_{\m(f)} \to Y$ extending $f$.
\end{lemma}
\begin{proof}
Let $r=\m(f)$, let $r_n \nearrow r$, and let $f_n : X_{r_n} \to Y$ be a conformal embedding extending $f$. Every $x \in X_r$ is contained in $X_{r_n}$ when $n$ is large enough since $X_{r_n} \nearrow X_r$. Define $\what f (x) = f_n(x)$. This does not depend on $n$ since $f_j = f_{j+1}$ on $X_{r_j}$ by the identity principle. The function $\what f$ is holomorphic and injective on $X_r$ because each $f_n$ is holomorphic and injective. The uniqueness of $\what f$ follows from the identity principle.
 \end{proof}

Similarly, $\m$ is upper semi-continuous on $\emb(X,Y,h)$. Again, this holds without any extra assumptions on $h$.

\begin{lemma} \label{lem:uppersemi}
The modulus of extension $\m$ is upper semi-continuous.
\end{lemma}

\begin{proof}
Suppose that $f_n \to f$ in $\emb(X,Y,h)$. Pass to a subsequence so that $r_n=\m(f_n)$ converges to some $r \in [0,\infty]$. We have to show that $\m(f) \geq r$. If $r=0$, then there is nothing to prove so we assume that $r>0$.

If $Y$ is not hyperbolic, then its group $\aut_0(Y)$ of conformal automorphisms homotopic to the identity acts simply transitively on $k$-tuples of distinct points for some $k \in \{1,2,3\}$. Pick some $k$-tuple $P\subset X$ of distinct points. Then there exists a sequence of automorphisms $T_n \in \aut_0(Y)$ converging to the identity such that $T_n\circ f_n(P)=f(P)$ for every $n$. Observe that $\m(f_n) = \m(T_n \circ f_n)$, so we might as well assume from the start that $f_n(P)=f(P)$ for every $n$. Under this assumption, we can consider $f_n$ and $f$ as maps from $X\setminus P$ to $Y \setminus f(P)$. Now $Y \setminus f(P)$ is hyperbolic, so the problem is reduced to that case.

Assume that $Y$ is hyperbolic. Let $\rho \in (0,r)$ and let $\what{f}_n$ be the conformal extension of $f_n$ to $X_{r_n}$. If $n$ is large enough, then $\rho\leq r_n$, and we let $g_n$ be the restriction of $\what{f}_n$ to $X_\rho$. Since $g_n$ cannot diverge locally uniformly (its restriction to $X$ converges to $f$), we may pass to a subsequence such that $g_n \to g$ for some $g \in \emb(X_\rho,Y,h)$ by Montel's theorem. The restriction of $g$ to $X$ is necessarily equal to $f$. In other words, $g$ is a conformal embedding extending $f$, so that $\m(f) \geq \rho$. Since $\rho \in (0,r)$ was arbitrary we have $\m(f) \geq r$. 
 \end{proof}

Since $\m$ is upper semi-continuous, it attains its maximum whenever the space $\emb(X,Y,h)$ is compact. This happens precisely when $h$ is non-trivial and non-parabolic by Lemma \ref{lem:compactness}. Let us recall hypothesis \myhypertarget{h3}{(\textbf{H})} from the introduction:
\begin{itemize}
\item $h:X\to Y$ is a non-trivial and non-parabolic embedding between finite Riemann surfaces;
\item $X$ has non-empty ideal boundary;
\item $\{X_r\}_{r\in[0,\infty]}$ is a $1$-parameter family of enlargements of $X$;
\item $\m$ is the associated modulus of extension;
\item $\emb(X,Y,h)$ is non-empty.
\end{itemize}
Note that this hypothesis prevents $Y$ from being the sphere with at most $2$ punctures or the disk with at most $1$ puncture. Under hypothesis \myhyperlink{h3}, $\m$ a\-chie\-ves its maximum at some $f\in\emb(X,Y,h)$. We will see that any such maxi\-mum is the restriction of a slit mapping. We first need to show that any limit of Teichm\"uller embeddings is itself a Teichm\"uller embedding.

\begin{definition}
Let $X_n \in \teich(X)$ and $Y_n \in \teich(Y)$ be such that $X_n \to X$ and $Y_n \to Y$ as $n \to \infty$, and let $\sigma_n: X_n \to X$ and $\tau_n : Y_n \to Y$ be the Teichm\"uller homeomorphisms that respect the markings. Let $h_n:X_n \to Y_n$ and $h:X \to Y$ be any maps. We say that $h_n \to h$ as $n \to \infty$ if $\tau_n \circ  h_n \circ \sigma_n^{-1} \to h$ locally uniformly on $X$ as $n \to \infty$. 
\end{definition}

\begin{lemma} \label{lem:limteich}
Let $h:X \to Y$ be a non-trivial and non-parabolic embedding between finite Riemann surfaces. Suppose that $X_n \to X$ in $\teich(X)$ and $Y_n \to Y$ in $\teich(Y)$, and let $f_n : X_n \to Y_n$ be a sequence of Teichm\"uller embeddings homotopic to $h$. Then there exists a Teichm\"uller embedding $f:X \to Y$ homotopic to $h$ such that $f_n \to f$ after passing to a subsequence.
\end{lemma}

Note that passing to a subsequence might be necessary due to the non-uniqueness of Teichm\"uller embeddings.

\begin{proof}
By Lemma \ref{lem:compactness}, we may pass to a subsequence such that $f_n \to f$ for some quasiconformal embedding $f: X \to Y$ homotopic to $h$. By Theorem \ref{thm:IoffeUniqueness}, each $f_n$ has minimal dilatation in its homotopy class. It follows that $f$ has minimal dilatation in its homotopy class. If $f$ is not conformal, then it is a Teichm\"uller embedding by Theorem \ref{thm:IoffeExistence}, and we are done.

Suppose that $f$ is conformal. We need to show that $f$ is a slit mapping. Let $\psi_n$ be the terminal quadratic differential of $f_n$, and let $g_n : Y_n \to Y_n'$ be a Teichm\"uller homeomorphism of dilatation $e^2$ and initial quadratic differential $\psi_n$. Since $d(Y_n,Y_n')=1$ for every $n$, $Y_n \to Y$ as $n \to \infty$, and closed balls in $\teich(Y)$ are compact, we may pass to a subsequence such that $Y_n'\to Y'$ as $n\to \infty$, for some $Y' \in \teich(Y)$ with $d(Y,Y')=1$. Let $g: Y \to Y'$ be a Teichm\"uller homeomorphism that respects the markings. By a standard argument similar to the one above, $g_n\to g$ after passing to yet another subsequence (this is only necessary if $Y$ is an annulus or a torus).

By construction, $g_n \circ f_n : X_n \to Y_n'$ is a Teichm\"uller embedding. Moreover, $g_n \circ f_n \to g \circ f$. As in the first paragraph of this proof, $g \circ f$ has minimal dilatation in its homotopy class. This dilatation is equal to $e^2>1$, so that $g \circ f$ is a Teichm\"uller embedding. Since $f$ is conformal, the terminal quadratic differential of $g \circ f$ is equal to the terminal quadratic differential of $g$. It follows that $f= g^{-1} \circ (g \circ f)$ is a slit mapping with respect to the initial quadratic differential of $g$. 
 \end{proof}

 We obtain the following characterization of the global maxima of $\m$.

\begin{proposition}\label{prop:globalmax}
Assume hypothesis \myhyperlink{h3}. Let $f$ be a global maximum of $\m$ with $\m(f)<\infty$ and let $\what f$ be the conformal extension of $f$ to $X_{\m(f)}$. Then $\what f$ is a slit mapping.
\end{proposition}

\begin{proof}
Let $R=\m(f)$. For every $r>R$, there is no conformal embedding $g: X_r \to Y$ whose restriction to $X$ is homotopic to $h$, for otherwise we would have $\m(g|_X)\geq r > R = \m(f)$. By Corollary \ref{cor:mindil} and Theorem \ref{thm:IoffeExistence}, there exists a Teichm\"uller embedding $g_r : X_r \to Y$ whose restriction to $X$ is homotopic to $h$. It is easy to see that $X_r$ moves continuously in $\teich(X)$ as a function of $r \in [0,\infty)$. By Lemma \ref{lem:limteich}, we can extract a limiting Teichm\"uller embedding $g: X_R \to Y$ from some subsequence of $g_r$ as $r \to R$. Since $\what f$ is conformal and homotopic to $g$, Theorem \ref{thm:IoffeUniqueness} implies that $\what f$ is itself a Teichm\"uller embedding. A conformal Teichm\"uller embedding is a slit mapping.
 \end{proof}

\begin{remark} \label{rem:triplypunc}
 Observe that every end of the surface $X_\infty$ is a puncture since the cylinder $S^1 \times [0,\infty)$ is conformally isomorphic to $\overline \DD \setminus \{0\}$. Thus if $\m(f)=\infty$, then the extension $\what f$ extends further to a conformal homeomorphism between $\what{X_\infty}$ and $\what{Y}$. In particular, $Y \setminus \what{f} (X_\infty)$ is finite, so that $\what f$ is a slit mapping with respect to any $\psi \in \qd^+(Y)$. This is unless $Y$ is the triply-punctured sphere in which case $\qd^+(Y)$ is empty. Thus the hypothesis $\m(f)<\infty$ in the above theorem is superfluous unless $Y$ is the triply-punctured sphere. In the latter case $\what f$ is nevertheless unique in its homotopy class by Lemma \ref{lem:emptyidealbdry}.
\end{remark}

We can in fact strengthen Proposition \ref{prop:globalmax} by replacing the word ``global'' with ``local''. The proof requires another lemma similar to Lemma \ref{lem:limteich}.

\begin{lemma} \label{lem:exhaustion} Let $h:X \to Y$ be a non-trivial and non-parabolic embedding between finite Riemann surfaces. Suppose that $Y_n \nearrow Y$, where the inclusion $\iota_n:Y_n \hookrightarrow Y$ is homotopic to a homeomorphism, and let $f_n : X \to Y_n$ be a sequence of Teichm\"uller embeddings such that $\iota_n \circ f_n $ is homotopic to $h$. Then there exists a Teichm\"uller embedding $f:X \to Y$ homotopic to $h$ such that $f_n \to f$ after passing to a subsequence.
\end{lemma}

The difference with Lemma \ref{lem:limteich} is that $Y_n$ is allowed to live in a diffe\-rent Teichm\"uller space than $Y$. For example, a surface with punctures can be exhausted by surfaces with holes. A priori, the terminal quadratic differentials could converge to a quadratic differential with double poles at the punctures of $Y$ and then the limiting map would not be a Teichm\"uller embedding. What prevents this from happening is that the domain $X$ is fixed.  

\begin{proof}
Let $K_n$ be the dilatation of $f_n$. The sequence $K_n$ is non-increasing and thus converges to a limit $K\geq 1$. After passing to a subsequence, $f_n$ converges to a $K$-quasiconformal embedding $f:X \to Y$.

We claim that $f$ has minimal dilatation in its homotopy class. If not, let $g:X \to Y$ be a quasiconformal embedding homotopic to $f$ such that $\Dil(g) < \Dil(f)$. If $r>0$ is small enough, then there is a quasiconformal ho\-meo\-morphism $\sigma: X_r \to X$ of dilatation strictly smaller than $\Dil(f)/\Dil(g)$. Thus the quasiconformal embedding $F$ consisting of the inclusion $X \hookrightarrow X_r$ followed by $g \circ \sigma$ has dilatation strictly less than $\Dil(f)$. Then $F(X) \subset Y_n$ whenever $n$ is large enough. Indeed, the ends of $F(X)$ which are not punctures are contained in a compact subset of $Y$ since they are surrounded by the union of collars $g \circ \sigma(X_r \setminus X)$. Thus the only way to go to infinity in the closure $\overline{F(X)}$ relative to $Y$ is via punctures of $F(X)$ that are also punctures of $Y$. For every such puncture $p$ and every $n\in \NN$, a neighborhood of $p$ in $Y$ is contained in $Y_n$. This is because $f_n$ is quasiconformal and as such it cannot map punctures of $X$ to ends of $Y_n$ with finite modulus. Since $Y_n \nearrow Y$, the inclusion $F(X) \subset Y_n$ holds for $n$ large enough. But the inequality $\Dil(F)< \Dil(f) \leq \Dil(f_n)$ contradicts the hypothesis that $f_n$ is a Teichm\"uller embedding and hence has minimal dilatation in its homotopy class.

If $f$ is not conformal, then it is a Teichm\"uller embedding by Theorem \ref{thm:IoffeExistence}, and we are done. Thus suppose that $f$ is conformal but is not a slit mapping. By Proposition \ref{prop:globalmax}, there exists an $r>0$ and a conformal embedding $g : X_r \to Y$ whose restriction to $X$ is homotopic to $f$. By the above argument, we have $g(X_{r/2}) \subset Y_n$ if $n$ is large enough, and thus $g(X) \subset Y_n$ with complement having non-empty interior. On the other hand, the restriction $g|_X : X \to Y_n$ is conformal and homotopic to the Teichm\"uller embedding $f_n:X \to Y_n$. By Theorem \ref{thm:IoffeUniqueness}, $g|_X$ is a slit mapping so that $Y_n \setminus g(X)$ has empty interior, contradiction.
 \end{proof}

We come to the main result of this section, which is that every local maxi\-mum of $\m$ is the restriction of a slit mapping. 

\begin{theorem}\label{thm:localmax}
Assume hypothesis \myhyperlink{h3}. Let $f$ be a local maximum of $\m$ with $\m(f)<\infty$, and let $\what f$ be the conformal extension of $f$ to $X_{\m(f)}$. Then $\what f$ is a slit mapping.
\end{theorem}

\begin{proof}
Let $R = \m(f)$. We first show that the complement $Y \setminus \what{f}(X_R)$ is horizontal for some meromorphic quadratic differential on $Y$, and is in particular an analytic graph. Let $\{x_1,x_2,\ldots\}$ be a dense subset of $X$. 

\begin{claim}
There exists an $n \in \NN$ such that if $g \in \emb(X,Y,h)$ satisfies $g(x_j)=f(x_j)$ for every $j \in \{1,...,n\}$, then $\m(g) \leq \m(f)$.
\end{claim}
\begin{proof}[Proof of Claim]
Suppose on the contrary that for every $n \in \NN$ there exists an element $g_n$ of $\emb(X,Y,h)$ satisfying $g_n(x_j)=f(x_j)$ for every $j \in \{1,...,n\}$ such that $\m(g_n) > \m(f)$. As $\emb(X,Y,h)$ is compact, every subsequence of $\{g_n\}_{n=1}^\infty$ has a subsequence converging to some $g \in \emb(X,Y,h)$. Any limit $g$ agrees with $f$ on the dense set $\{x_1,x_2,....\}$, and hence is equal to $f$. Thus $g_n \to f$ with $\m(g_n) > \m(f)$. This contradicts the hypothesis that $f$ is a local maximum of $\m$.
 \end{proof}

Let $n$ be as in the claim, and let $P = \{x_1, x_2, \ldots, x_n \}$. Then for every $r>R$, there is no conformal embedding $g: X_r \to Y$ homotopic to $f$ rel $P$. By Corollary \ref{cor:mindil} and Theorem \ref{thm:IoffeExistence}, there exists a Teichm\"uller embedding $g_r : X_r\setminus P \to Y\setminus {f}(P)$ homotopic to $f$ rel $P$.

 Let $g$ be any limit of any subsequence of $g_r$ as $r\searrow R$. Then $g: X_R\setminus P \to Y\setminus f(P)$ is a Teichm\"uller embedding by Lemma \ref{lem:limteich}. Since $\what f$ is conformal and homotopic to $g$ rel $P$, Theorem \ref{thm:IoffeUniqueness} implies that $\what f$ is itself a slit mapping, considered as a map from $X_R\setminus P$ to $Y\setminus {f}(P)$. Therefore the complement $Y \setminus \what f (X_R)$ is a finite union of horizontal arcs for some meromorphic quadratic differential on $Y$, possibly with simple poles on the set $f(P)$.
 
Let $\Gamma = Y \setminus \what f (X_R)$, let $\{y_1,y_2,...\}$ be a dense subset of the graph $\
\Gamma$ minus its vertices, and fix a Riemannian metric on $Y$. 

\begin{claim}
There exists a $k\in \NN$ such that for every $r>R$ and every $\eps > 0$, there is no conformal embedding $g : X_r \to Y$ homotopic to $h$ whose image is disjoint from the balls $\overline{B}(y_1,\eps),\ldots,\overline{B}_\eps(y_k,\eps)$.
\end{claim}
\begin{proof}[Proof of Claim]
Suppose that for every $k\in \NN$ there exist an $r_k>R$, an $\eps_k > 0$, and a con\-for\-mal em\-bed\-ding $g_k : X_{r_k} \to Y$ whose res\-tric\-tion to $X$ is homotopic to $h$ such that $g_k(X_{r_k})$ is disjoint from the balls $\overline{B}(y_1,\eps_k),\ldots,\overline{B}(y_k,\eps_k)$. We may assume that $r_k \to R$ and $\eps_k \to 0$. Let $g$ be any limit of any subsequence of the sequence $\{g_k\}$. Then $g(X_R)$ is disjoint from the set $\{y_1,y_2,...\}$ and hence from its closure $\Gamma$, so that $\what{f}^{-1} \circ g : X_R \to X_R$ is a conformal embedding homotopic to the identity. If $X_R$ is not an annulus, then Lemma \ref{lem:selfembisidentity} implies that $g=\what{f}$ and hence $g_k \to \what{f}$. If $X_R$ is an annulus, then we may pre-compose each $g_k$ by a rotation so that we still get $g_k \to \what{f}$. Since $\m(g_k|_X) \geq r_k > R = \m(f)$, this contradicts the hypothesis that $f$ is a local maximum of $\m$. 
 \end{proof}

Let $k$ be as in the last claim, and let $Q=\{y_1,\ldots, y_k\}$. For each $\eps>0$, let $Y_\eps = Y \setminus (\overline{B}(y_1,\eps)\cup \cdots \cup \overline{B}(y_k,\eps))$. Let $\iota_\eps : Y_\eps \to Y\setminus Q$ be a homeomorphism homotopic to the inclusion map, and let $h_\eps = \iota_\eps^{-1} \circ \what{f}$. The embedding $h_\eps : X_R \to Y_\eps$ followed by the inclusion $Y_\eps \hookrightarrow Y$ is homotopic to $h$. By the claim, for every $r>R$, there is no conformal embedding $g : X_r \to Y_\eps$ homotopic to $h_\eps$. Therefore, there is a Teichm\"uller embedding $g_\eps^r: X_r \to Y_\eps$ homotopic to $h_\eps$. Letting $r\searrow R$, we can extract a limiting Teichm\"uller embedding $g_\eps: X_R \to Y_\eps$ by Lemma \ref{lem:limteich}.

Since $Y_\eps \nearrow (Y \setminus Q)$ as $\eps \searrow 0$, we can apply Lemma \ref{lem:exhaustion} and obtain a Teichm\"uller embedding $g : X_R \to Y \setminus Q$ as a limit of a subsequence of $\{g_\eps\}$. Since $\what{f} : X_R \to Y\setminus Q$ is homotopic to $g$, it is a slit mapping with respect to some $\psi \in \qd^+(Y \setminus Q)$. Thus $\psi$ is meromorphic on $Y$ with at most simple poles on $Q$. Moreover, the graph $\Gamma = Y \setminus \what{f}(X_R)$ is horizontal for $\psi$. Since every point of $Q$ is contained in the interior of an edge of $\Gamma$, the quadratic differential $\psi$ cannot have simple poles on $Q$. Indeed, there is only one horizontal ray emanating from any simple pole. Therefore, $\psi$ is holomorphic on $Y$ and the map $\what f : X_R \to Y$ is an honest slit mapping.

 \end{proof}

Once again, the hypothesis $\m(f) < \infty$ is not necessary unless $Y$ is the triply-punctured sphere.

As a consequence of Theorem \ref{thm:localmax}, every conformal embedding which is not a slit mapping can be approximated by a sequence of conformal embeddings each of which extends by some amount.

\begin{corollary} \label{cor:tameapprox}
Assume hypothesis \myhyperlink{h3} and let $g \in \emb(X,Y,h)$. If $g$ is not a slit mapping, then there is a sequence $\{g_n\}$ converging to $g$ in $\emb(X,Y,h)$ such that $g_n$ extends to a conformal embedding of $X_{r_n}$ into $Y$ for some $r_n>0$. 
\end{corollary}
\begin{proof}
If $\m(g)>0$, then we can take $g_n=g$. If $\m(g)=0$ but $g$ is not a local maximum of $\m$, then there exists a sequence $g_n \to g$ with $\m(g_n)>0$. If $\m(g)=0$ and $g$ is a local maximum of $\m$, then $g$ is a slit mapping by the previous theorem.
 \end{proof}

A strong converse to Theorem \ref{thm:localmax} holds due to Theorem \ref{thm:IoffeUniqueness}. 

\begin{lemma} \label{lem:localisglobal} 
Assume hypothesis \myhyperlink{h3}. Suppose that $g: X_r \to Y$ is a slit mapping such that $g|_X$ is homotopic to $h$. Then $g|_X$ is a global maximum of $\m$.
\end{lemma}
\begin{proof}
First observe that $\m(g|_X) \geq r$. Suppose that $\m(f) \geq \m(g|_X)$ for some element $f$ of $\emb(X,Y,h)$ and let $\what f$ be the maximal extension of $f$. Then $\what{f}|_{X_r}$ is homotopic to $g$. By Theorem \ref{thm:IoffeUniqueness}, $\what{f}|_{X_r}$ is a slit mapping. In particular, the complement of $\what{f} (X_r)$ has empty interior in $Y$. Therefore $X_{\m(f)} \setminus X_r$ is empty so that $\m(f)\leq r \leq \m(g|_X)$.
 \end{proof}

Furthermore, the almost rigidity of slit mappings implies that the set of local maxima of $\m$ is path-connected. 

\begin{lemma} \label{lem:setoflocalmax}
Assume hypothesis \myhyperlink{h3}. The set $M$ of all local maxima of $\m$ is homeomorphic to either a point, a compact interval, a circle, or a torus, and $\m$ is constant on $M$.
\end{lemma}
\begin{proof}
Suppose first that there is some $f \in M$ such that $\m(f) < \infty$. Then by Theorem \ref{thm:localmax}, the maximal extension $\what f$ is a slit mapping. By Lemma \ref{lem:localisglobal}, $f$ is a global maximum of $\m$. In particular, $\m(g) < \infty$ for every $g \in M$ and thus every $g \in M$ is a global maximum of $\m$. In particular, $\m$ is constant on $M$, say equal to $R$. The map $M \to \emb(X_R, Y , h)$ defined by $f\mapsto \what f$ is a homeomorphism with inverse $g \mapsto g|_X$. By Theorem \ref{thm:AlmostRigidity}, the evaluation map $\emb(X_R,Y,h) \to Y$ sending $f$ to $f(x)$ is a homeomorphism onto its image for every $x \in X_R$, and its image is either a point, a compact interval, a circle, or a torus.

Otherwise, $\m$ is constant equal to $\infty$ on $M$. In this case $M$ is homeomorphic to $\emb(X_\infty,Y,h)$, which is the same as $\aut_0(Y)$. This is either a point or a torus (see Subsection \ref{subsec:aut}).   
 \end{proof}

Theorem \ref{thm:localmax}, Lemma \ref{lem:localisglobal} and Lemma \ref{lem:setoflocalmax} together imply Theorem \ref{thm:localmaxintro} from the introduction. The fact that $\m$ has a connected plateau of local maxima easily implies that the space $\emb(X,Y,h)$ is connected.

\begin{proof}[Proof of Theorem \ref{thm:connected}]
Suppose that $\emb(X,Y,h)=E_0 \cup E_1$ where $E_0$ and $E_1$ are disjoint non-empty closed sets. Then each of $E_0$ and $E_1$ is both compact and open. Since $\m$ is upper semi-continuous, the restriction $\m|_{E_j}$ attains its maximum at some $f_j\in E_j$. Then $f_j$ is a local maximum of $\m$ since $E_j$ is open. By Lemma \ref{lem:setoflocalmax}, $f_0$ and $f_1$ are both contained in a connected subset $M$ of $\emb(X,Y,h)$. But $M = (M \cap E_0) \cup (M \cap E_1)$ is disconnected. Contradiction.
 \end{proof}

\section{The blob and its boundary} \label{sec:blob}

We say that a continuous map $h:X \to Y$ between any Riemann surfaces is \emph{generic} if the induced homomorphism $\pi_1(h):\pi_1(X,x) \to \pi_1(Y,h(x))$ has non-abelian image. This implies that $X$ and $Y$ are hyperbolic. Given such a generic map $h$, let $\maps(X,Y,h)$ be the space of all continuous maps $f:X \to Y$ homotopic to $h$. The following lemma shows that for every $x \in X$ and every $f \in \maps(X,Y,h)$, there is a well-defined way to lift the image point $f(x)$ to the universal cover of $Y$.

\begin{lemma} \label{lem:lift}
Let $h:X \to Y$ be a generic map between Riemann surfaces. Let $H: X \times [0,1] \to Y$ be a homotopy from $h$ to some $f\in \maps(X,Y,h)$. Then for every $x \in X$, the homotopy class rel endpoints of the path $t \mapsto H(x,t)$  does not depend on the choice of $H$.
\end{lemma}
\begin{proof}
We use two standard facts about hyperbolic surfaces:
\begin{itemize}
\item every abelian subgroup of $\pi_1(Y,h(x))$ is cyclic;
\item every non-trivial element in $\pi_1(Y,h(x))$ is the positive power of a unique pri\-mi\-tive element.  
\end{itemize}

Let $G$ be any other homotopy from $h$ to $f$. By composing $H$ with $G$ run backwards, we get a homotopy from $h$ to itself, hence a map $F : X \times S^1 \to Y$. Suppose that the loop $\gamma(t)= F(x,t)$ is not trivial in $\pi_1(Y,h(x))$. Then it is equal to $\beta^k$ for some primitive element $\beta$ and some $k> 0$.

Let $\alpha$ be any loop in $X$ based at $x$. Then the map $S^1 \times S^1 \to Y$ given by $(s,t) \mapsto F(\alpha(s),t)$ induces a homomorphism of $\ZZ^2$ into $\pi_1(Y,h(x))$. The image of this homomorphism is cyclic, and contains both $[h\circ \alpha]$ and $[\gamma]$. From the existence and uniqueness of primitive roots in $\pi_1(Y,h(x))$, it follows that $[h\circ \alpha]=\beta^j$ for some $j \in \ZZ$. Since $\alpha$ was arbitrary, we deduce that the image of the homomorphism $\pi_1(h) : \pi_1(X,x) \to \pi_1(Y,h(x))$ is contained in the cyclic group $\langle \beta \rangle$. This contradicts the hypothesis that $h$ is generic. We conclude that the loop $\gamma(t) = F(x,t)$ is null-homotopic. Equivalently, the paths $t \mapsto H(x,t)$ and $t \mapsto G(x,t)$ are homotopic rel endpoints.

 \end{proof}

\begin{definition} \label{def:lift}
Let $h:X \to Y$ be a generic map between Riemann surfaces, let $f \in \maps(X,Y,h)$ and let $x \in X$. We define $\lift_x(f)$ to be the homotopy class rel endpoints of the path $t \mapsto H(x,t)$ in $Y$ where $H$ is any homotopy from $h$ to $f$. This is well-defined by Lemma \ref{lem:lift}. By definition, $\lift_x(f)$ is an element of the universal cover $\wtilde Y$ of $Y$ based at $h(x)$. This universal cover is isomorphic to $\DD$, but for our purposes it will be better to think of it as the set of homotopy classes of paths in $Y$ starting at $h(x)$. 
\end{definition}

\begin{remark}
If $h:X \to Y$ is a cyclic map between Riemann surfaces, i.e., such that the image of $\pi_1(h)$ is an infinite cyclic subgroup $C\leq \pi_1(Y,h(x))$, and $Y$ is not a torus then a similar construction defines a lift from $\maps(X,Y,h)$ to the annulus cover $A \to Y$ associated to $C$.  
\end{remark}

\begin{lemma} \label{lem:cont}
Let $h:X \to Y$ be a generic map between Riemann surfaces, where $X$ is finite. Then $\lift_x: \maps(X,Y,h) \to \wtilde Y$ is continuous.
\end{lemma}

\begin{proof}
Since $\maps(X,Y,h)$ is metrizable, it suffices to prove sequential continuity. Let $f_n,f \in \maps(X,Y,h)$ be such that $f_n \to f$ as $n \to \infty$. Let $K\subset X$ be a compact deformation retract of $X$ containing $x$. Let $\eps>0$ be smaller than the minimum of the injectivity radius of $Y$ over $f(K)$ with respect to the hyperbolic metric and let $n$ be large enough so that $|f_n-f|<\eps$ on $K$. For every $(\xi,t) \in K \times [0,1]$, let $F_n(\xi,t)$ be the point at proportion $t$ along the unique shortest length geodesic between $f(\xi)$ and $f_n(\xi)$ in $Y$. This gives a continuous homotopy from $f|_K$ to $f_n|K$. By composing the deformation retraction $X \to K$ with $F_n$, we get a homotopy $G_n$ from $f$ to $f_n$ moving points of $K$ by distance at most $\eps$.

Given any homotopy $H$ from $h$ to $f$, the concatenation $H\ast G_n$ (that is, $H$ followed by $G_n$) is a homotopy from $h$ to $f_n$. Thus $\lift_x(f_n)$ can be represented as $\alpha \ast \beta_n$ where $\alpha$ is any representative of $\lift_x(f)$ and $\beta_n(t)=G_n(x,t)$. Since the geodesic $\beta_n$ is contained in the ball of radius $\eps$ centered at $\alpha(1)$ in $Y$, and $\eps$ can be taken arbitrarily small, we have $\lift_x(f_n) \to \lift_x(f)$ as $n \to \infty$.
 \end{proof}

We now get back to conformal embeddings and look at where a given point can go under all conformal embeddings in a given homotopy class.

\begin{definition} Let $h: X \to Y$ be a generic embedding between finite Riemann surfaces and let $x \in X$. We define $\blob(x,X,Y,h)$ to be the image of $\emb(X,Y,h)$ by the map $\lift_x$ from Definition \ref{def:lift}. This is a subset of the universal cover $\wtilde Y$ of $Y$ based at $h(x)$. Given $c \in \wtilde Y$, we will write $\pi_Y(c)$ for its projection to $Y$, that is, the endpoint $\gamma(1)$ of any representative $\gamma \in c$. 
\end{definition}

Our previous results imply that $\blob(x,X,Y,h)$ is at most $1$-dimensional whenever $\emb(X,Y,h)$ contains a slit mapping.

\begin{proposition} \label{prop:blobslit}
Let $h: X \to Y$ be a generic embedding between finite Riemann surfaces and let $x \in X$. If $\emb(X,Y,h)$ contains a slit mapping, then $\blob(x,X,Y,h)$ is homeomorphic to a point or a compact interval.
\end{proposition}
\begin{proof}
By Theorem \ref{thm:AlmostRigidity}, the evaluation map $\ev_x:\emb(X,Y,h) \to Y$ sending $f$ to $f(x)$ is  a homeomorphism onto its image $V(x)$. Moreover, $V(x)$ is either a point or a compact interval since $Y$ is not a torus and $h$ is not cyclic. The restriction of the universal covering map $\pi_Y : \wtilde Y \to Y$ to $\blob(x,X,Y,h)$ is a homeomorphism onto its image $V(x)$ with inverse $\lift_x \circ \ev_x^{-1}$.
 \end{proof}

We will see that $\blob(x,X,Y,h)$ is not much more complicated when $\emb(X,Y,h)$ does not contain any slit map\-ping. Let us re\-call hypo\-the\-sis \myhypertarget{h4}{(\textbf{H'})} from the introduction:
\begin{itemize}
\item $h:X \to Y$ is a generic embedding between finite Riemann surfaces;
\item $\partial X \neq \varnothing$;
\item $\emb(X,Y,h)$ is non-empty and does not contain any slit mapping.
\end{itemize}

Assume hypothesis \myhyperlink{h4}. Since $\emb(X,Y,h)$ is compact (Lemma \ref{lem:compactness}) and connected (Theorem \ref{thm:connected}), and since $\lift_x$ is conti\-nuous (Lemma \ref{lem:cont}), the set $\blob(x,X,Y,h)$ is compact and connected for any $x \in X$. Our goal is to show that the blob is homeomorphic to a closed disk (Theorem \ref{thm:blobisdisk}). The strategy of the proof is to analyze the boundary of the blob. We will show that every point in $\partial\blob(x,X,Y,h)$ is attained by a special kind of map in $\emb(X,Y,h)$ which we call a slit mapping rel $x$. 

\begin{definition}
Let $h:X \to Y$ be a topological embedding between finite Riemann surfaces and let $x \in X$. We say that a map $f \in \maps(X,Y,h)$ is a \emph{Teichm\"uller embedding rel $x$} if the restriction $f^\star:X \setminus \{x\} \to Y \setminus \{f(x)\}$ is a Teichm\"uller embedding. A \emph{slit mapping rel $x$} is a Teichm\"uller embedding rel $x$ which is conformal.
\end{definition}

The distinction to make here is that the initial and terminal quadratic differentials of $f$ are allowed to have simple poles at $x$ and $f(x)$ respectively. For example, there are no slit mappings into the triply-punctured sphere, but plenty of slit mappings relative to a point. 

In order to characterize the boundary points of the blob, we first need a standard construction for pushing a point around on a Riemann surface. To paraphrase \cite[p.97]{FarbMargalit}: imagine placing your finger on a surface and pushing it along a smooth path, dragging the rest of the surface along as you go. The diffeomorphism obtained at the end is called a \emph{point-pushing diffeomorphism}.

\begin{lemma} \label{lem:push}
Let $Y$ be any Riemann surface and let $\gamma:[0,1] \to Y$ be a smooth immersion. There exists an isotopy $H: Y \times [0,1] \to Y$ through quasiconformal diffeomorphisms such that $H(y,0)=y$ for every $y \in Y$ and $H(\gamma(0),t)=\gamma(t)$ for every $t \in [0,1]$.
\end{lemma}
\begin{proof}
First assume that $\gamma$ embedded. Extend the vector field $\gamma_*(\partial/\partial t)$ to a smooth vector field $V$ supported in a tubular neighborhood of $\gamma([0,1])$, and define $H$ to be the vector flow along $V$. For each $t \in [0,1]$, the map $y \mapsto H(y,t)$ is a diffeomorphism which is the identity outside a compact set, hence is quasiconformal.

If $\gamma$ is not embedded, break it up into finitely many embedded subarcs $\gamma_j$, then construct an isotopy $H_j$ on each corresponding subinterval $[a_{j-1},a_j]$ using the above method. The isotopy $H(y,t)$ is defined as $H_1(y,t)$ for $t \in [a_0,a_1]$, then $H_2(H_1(y,a_1),t)$ for $t \in [a_1, a_2]$ and so on, picking up where we left off at each step.
 \end{proof}

Given a smooth immersion $\gamma:[0,1] \to Y$ into a Riemann surface and an isotopy $H$ as in Lemma \ref{lem:push}, we say that $y \mapsto \push_\gamma(y):=H(y,1)$ is a \emph{point-pushing diffeomorphism} along $\gamma$. By definition, $\push_\gamma$ is isotopic to the identity and satisfies $\push_\gamma(\gamma(0))=\gamma(1)$. Of course, $\push_\gamma$ depends on the specific choice of $H$ but its isotopy class rel $\gamma(0)$ only depends on the homotopy class of $\gamma$ by the Birman exact sequence \cite[Theorem 4.6]{FarbMargalit}.

Here are two elementary observations. If $h:X \to Y$ is a generic map between Riemann surfaces and $\gamma : [0,1] \to Y$ is an immersed curve such that $\gamma(0)=h(x)$, then $\lift_x(\push_\gamma \circ h)=[\gamma]$ by definition. Secondly, if two maps $f$ and $g$  in $\maps(X,Y,h)$ are homotopic rel $x$, then obviously $\lift_x(f)=\lift_x(g)$. The converse is also true, as can be seen using point-pushing on $Y$. Thus $\lift_x$ detects when two homotopic maps are homotopic rel $x$.

We say that a map $f \in \maps(X,Y,h)$ \emph{realizes} a point $c \in \wtilde Y$ if $\lift_x(f)=c$. We can now prove that points outside the blob are realized by Teichm\"uller embeddings rel $x$, provided that there is some quasiconformal embedding in the homotopy class of $h$.

\begin{lemma} \label{lem:bdryblobteich}
Let $h:X \to Y$ be a generic quasiconformal embedding between finite Riemann surfaces, let $x \in X$, and let $c \in \wtilde Y \setminus \blob(x,X,Y,h)$. Then there exists a unique Teichm\"uller embedding $f$ rel $x$ homotopic to $h$ such that $\lift_x(f)=c$. 
\end{lemma}
\begin{proof}
Let $\gamma: [0,1]\to Y$ be a smooth immersed curve in $c$, for example the hyperbolic geodesic, and let $\push_\gamma$ be a quasiconformal point-pushing diffeomorphism along $\gamma$. Then $\push_\gamma \circ h$ is a quasiconformal embedding homotopic to $h$ such that $\lift_x(\push_\gamma \circ h) = [\gamma] = c$. Let $f^\star:X \setminus \{x\} \to Y \setminus \{ \gamma(1) \}$ be a quasiconformal embedding homotopic to the restriction of $\push_\gamma \circ h$ with minimal dilatation (Corollary \ref{cor:mindil}) and let $f:X \to Y$ be its extension. Since $f$ is ho\-mo\-to\-pic to $\push_\gamma \circ h$ rel $x$, we have $\lift_x(f) = \lift_x(\push_\gamma \circ h)=c$.

Suppose that $f^\star$ is conformal. Then $f\in\emb(X,Y,h)$ so that $c=\lift_x( f)$ belongs to $\blob(x,X,Y,h)$, contrary to our assumption. By Theo\-rem \ref{thm:IoffeExistence}, $f^\star$ is a Teichm\"uller embedding, so that $f$ is a Teichm\"uller embedding rel $x$. 

If $g$ is another Teichm\"uller embedding rel $x$ homotopic to $h$ such that $\lift_x(g)=c$, then $g$ is homotopic to $f$ rel $x$. Thus the restrictions $g^\star$ and $f^\star$ to $X \setminus \{x\} \to Y \setminus \{\gamma(1)\}$ are homotopic Teichm\"uller embeddings with an essential puncture, and we conclude that $g=f$ by Remark \ref{rem:nopunctures2}.
 \end{proof}

We deduce that points on the boundary of the blob are realized by slit mappings rel $x$.

\begin{proposition} \label{prop:bdryblob}
Assume \myhyperlink{h4}, let $x \in X$ and let $c \in \partial \blob(x,X,Y,h)$. Then there is a unique $f \in \emb(X,Y,h)$ such that $\lift_x(f)=c$. Moreover, $f$ is a slit mapping rel $x$. If $\vphi$ and $\psi$ are initial and terminal quadratic differentials for $f$ rel $x$, then $\vphi$ has a simple pole at $x$ and $\psi$ has a simple pole at $f(x)$.
\end{proposition}

\begin{proof}
Since $\blob(x,X,Y,h)$ is closed, $c$ belongs to $\blob(x,X,Y,h)$. Hence there exists some $f \in \emb(X,Y,h)$ such that $\lift_x(f)=c$. 

Let $\{c_n\}_{n=1}^\infty$ be a sequence in $\wtilde Y \setminus \blob(x,X,Y,h)$ such that $c_n \to c$ as $n \to \infty$. Then  $Y \setminus \{\pi_Y(c_n)\}$ converges to $Y \setminus \{\pi_Y(c)\}$ as $n \to \infty$, where we can take markings to be the identity outside of a small neighborhood of $\pi_Y(c)$. By the previous lemma, there exists a Teichm\"uller embedding $f_n:X \to Y$ rel $x$ homotopic to $h$ such that $\lift_x(f_n)=c_n$. By Lemma \ref{lem:limteich}, we can extract a subsequence of the restrictions $f_n^\star: X \setminus \{x\} \to Y \setminus \{\pi_Y(c_n)\}$ converging to a Teichm\"uller embedding $g: X\setminus \{x\} \to Y \setminus \{\pi_Y(c)\}$ which is homotopic to the restriction $f^\star$ of $f$. As $f^\star$ is conformal, it is a slit mapping by Theorem \ref{thm:IoffeUniqueness}, so that $f$ is a slit mapping rel $x$. For any $g \in \emb(X,Y,h)$ realizing $c$, we have that $g$ is homotopic to $f$ rel $x$. Hence $g=f$ by Theorem \ref{thm:IoffeUniqueness} and Remark \ref{rem:nopunctures2}.

Suppose that $\psi$ does not have a pole at $y$. Then $f$ is an honest slit mapping from $X$ to $Y$. But we assumed that $\emb(X,Y,h)$ does not contain any slit mapping. Therefore $\psi$ has a simple pole at $f(x)$ and $\vphi=f^*\psi$ has a simple pole at $x$.
 \end{proof}

Assume hypothesis \myhyperlink{h4} and let $x \in X$. Suppose that $c\in \partial\blob(x,X,Y,h)$ and  that $f\in \emb(X,Y,h)$ realizes $c$. We say that any quadratic differential $\psi \in \qd^+(Y \setminus \{f(x)\})$ with respect to which $f$ is a slit mapping rel $x$ \emph{realizes} $c$ as well. Even though $f$ is unique by Proposition \ref{prop:bdryblob}, $\psi$ need not be unique up to scale (see Remark \ref{rem:termindiff}). Nevertheless, the set of quadratic differentials realizing a given point $c\in \partial \blob(x,X,Y,h)$ is convex.

\begin{lemma} \label{lem:qdconvex}
Assume hypothesis \myhyperlink{h4}. Suppose that $\psi_0$ and $\psi_1$ in $\qd^+(Y\setminus\{y\})$ realize $c\in \partial \blob(x,X,Y,h)$ where $y = \pi_Y(c)$. Then for every $\alpha,\beta>0$ the quadratic differential $\alpha \psi_0 + \beta \psi_1$ belongs to $\qd^+(Y\setminus\{y\})$, realizes $c$, and has a simple pole at $y$. 
\end{lemma}
\begin{proof}
Let $f$ be the slit mapping realizing $c$. Then $\psi :=\alpha \psi_0 + \beta \psi_1$ is $\geq 0$ along any arc in $\what Y \setminus f(X)$ so that $f$ is a slit mapping rel $x$ with respect to $\psi$. If $\psi$ does not have a pole at $y$, then $f$ is a slit mapping from $X$ to $Y$, contradicting hypothesis \myhyperlink{h4}.
 \end{proof}

We will see that any $\psi\in \qd^+(Y\setminus\{\pi_Y(c)\})$ which realizes $c$ tells us something about the shape of $\blob(x,X,Y,h)$ near $c$. More precisely, $\psi$ can be used to find a region $U\subset \wtilde Y$ of forbidden values for the maps in $\emb(X,Y,h)$, that is, disjoint from $\blob(x,X,Y,h)$. The idea is that if we push $\pi_Y(c)$ in some directions, then a certain quantity goes up, while it cannot increase under conformal embedding. This quantity is extremal length.

\section{Extremal length of partial measured foliations} \label{sec:EL}

There are several equivalent ways to define the extremal length of a measured foliation. We follow the approach developed in \cite{GardinerLakic2} and \cite{GardinerLakic}. 

A map will be called \emph{almost-smooth} if it is continuous and continuously differentiable except perhaps at finitely many points.

\begin{definition}
A \emph{partial measured foliation} $F=\{(U_j,v_j)\}_{j \in J}$ on a Riemann surface $X$ is a collection of open sets $U_j \subset X$ together with almost-smooth functions $v_j : U_j \to \RR$ satisfying
$$
v_j = \pm v_k + c_{jk}
$$
on $U_j \cap U_k$, where $c_{jk}$ is locally constant. The set $U = \bigcup_j U_j$ is called the \emph{support} of $F$. Since $|\mathrm{d}v_j|=|\mathrm{d}v_k|$ on $U_j\cap U_k$, these patch up to a well-defined object $|\mathrm{d}v|$ on $U$ called the \emph{transverse measure} of $F$. The \emph{leaves} of $F$ are the maximal connected subsets of $U$ on which each $v_j$ is locally constant. We will write $F$ or $|\mathrm{d}v|$ interchangeably.
\end{definition}

\begin{remark}
This is much weaker that usual notion of measured foliation \cite{FLP}. For one thing, the support does not have to be the complement of a finite set. Moreover, the leaves are not required to be $1$-dimensional submanifolds; they can be thick. In practice, we will only deal with partial measured foliations which are measured foliations on a subsurface.
\end{remark}

\begin{remark}
We could relax the regularity condition on the functions $v_j$ and only assume that they belong to the Sobolev space $W^{1,2}(U_j)$. This would be more natural from the point of view of quasiconformal maps. For the sake of simplicity we will stick to the almost-smooth condition.  
\end{remark}

For example, if $\varphi$ is a quadratic differential and $U_j$ is a simply connected domain on which $\varphi$ does not have any singularities, then the function
$$
v_j(z) = \im  \int_{z_0}^z \sqrt{\varphi}
$$
is well-defined on $U_j$ up to an additive constant and sign. The resulting partial measured foliation $\mathcal{F}^h\varphi=|\mathrm{d}v|=|\im \sqrt{\varphi}|$ is the \emph{horizontal foliation} of $\vphi$. Its leaves are the regular horizontal trajectories of $\varphi$. 
 
The \emph{Dirichlet energy} of a partial measured foliation $F=|\mathrm{d}v|$ with support $U$ is
$$
\dir(F):= \int_U |\nabla v|^2 \mathrm{d}A 
$$
where the gradient $\nabla v$ (only defined up to sign) is computed with respect to any conformal Riemannian metric on $X$ with corresponding area form $\mathrm{d}A$. Alternatively, we can write
$$\dir(F) = \int_U  \left(\frac{\partial v}{\partial x}\right)^2 +\left(\frac{\partial v}{\partial y}\right)^2  \mathrm{d}x\mathrm{d}y$$
where $z=x+iy$ is any conformal coordinate.

For example, if $\varphi$ is a holomorphic quadratic differential on $X$ and $\mathcal{F}^h\varphi$ is its horizontal foliation, then 
$$
\dir(\mathcal{F}^h\varphi) = \int_X |\varphi| = \|\varphi\|
$$
as can be seen by computing the Dirichlet energy in natural coordinates where $\varphi=\mathrm{d}z^2$ and $v(x+iy)=\pm y + c$.

A \emph{cross-cut} on a Riemann surface $X$ is a proper arc $\alpha : (0,1) \to X$. Two cross-cuts are \emph{homotopic} if there is a homotopy through cross-cuts between them. The \emph{height} of a homotopy class $c$ of closed curves or cross-cuts in $X$ with respect to a partial measured foliation $F=|\mathrm{d}v|$ is the quantity 
$$\hgt(c,F):=\inf_{\alpha \in c} \int_\alpha |\mathrm{d}v|,$$
where the infimum is taken over all piecewise smooth curves $\alpha \in c$, and where $|\mathrm{d}v|$ is extended to be zero outside its support.

Given partial measured foliations $F$ and $G$ on a Riemann surface $X$, we say that $G$ \emph{dominates} $F$ if
$$
\hgt(c,G) \geq \hgt(c,F)
$$
for every homotopy class $c$ of closed curves or cross-cuts in $X$. Two partial measured foliations are \emph{measure equivalent} if they dominate each other, i.e., if they have the same heights.

\begin{definition}
The \emph{extremal length} of a partial measured foliation $F$ on a Riemann surface $X$ is defined as $$\el(F , X):= \inf \left\{\, \dir(G) : G \text{ dominates } F \,\right\}.$$
\end{definition}

This is not the standard definition of extremal length. However, if $F$ is a measured foliation (not partial) on a finite Riemann surface $X$, then $\el(F , X)$ is the norm of the unique holomorphic quadratic differential $\Phi(F,X)$ on $X$ with the same heights as $F$ (this is often taken as a definition). The existence and uniqueness of $\Phi(F,X)$ is due to Hubbard and Masur \cite{HubbardMasur}, but what interests us here is its minimizing property.

\begin{theorem}[The minimal norm property] \label{thm:QDsrealizeEL}
Let $X$ be a finite Riemann surface, let $\vphi \in \qd(X) \setminus \{0\}$, and let $\mathcal{F}^h\vphi$ be the horizontal foliation of $\vphi$. Then $$\el(\mathcal{F}^h\vphi , X)=\dir(\mathcal{F}^h\vphi)=\|\vphi\|.$$ That is,
$$
\dir(\mathcal{F}^h\vphi) \leq \dir(G)
$$
 for every partial measured foliation $G$ on $X$ which dominates $\mathcal{F}^h\vphi$.
\end{theorem}

This is proved in \cite[Theorem 3.2]{MardenStrebel} (see also \cite[Chapter VII]{Strebel}). There the result is stated for partial measured foliations $G$ whose support has a discrete complement, but this hypothesis is not used anywhere in the proof. The idea of the proof is to look at the vertical foliation $\mathcal{F}^v\vphi$ which splits up into rectangles, cylinders and minimal components. If $\lambda$ is a vertical leaf which is a cross-cut or a closed curve, then $\int_\lambda \mathrm{d}\mathcal{F}^h\vphi = \hgt(\lambda,\mathcal{F}^h\vphi) \leq \int_\lambda \mathrm{d} G .$ There is a similar inequality for minimal components, but one needs to pick a horizontal transversal and do a surgery on vertical leafs to obtain closed curves. The result is obtained by integrating these $1$-dimensional inequalities against the leaf space and applying the Cauchy-Schwarz inequality. This is a sophisticated version of the so-called ``length-area argument''.

\subsection{Pushing forward}

If $f : X \to Y$ is an almost-smooth embedding between Riemann surfaces and $F=\{(U_j,v_j)\}_{j \in J}$ is a partial measured foliation on $X$, then the \emph{push-forward} $f_*F = \{(f(U_j),v_j \circ f^{-1})\}_{j \in J}$ is a partial measured foliation on $Y$. This justifies our preference for partial measured foliations over measured foliations: we will push them forward by conformal embeddings.
 
\begin{remark} 
The notion of heights induces a topology on the set of (measure equivalence classes of) partial measured foliations on a finite Riemann surface. The reader should be warned that the push-forward operator is not continuous with res\-pect to that topology \cite[Example 4.4]{DylanEtAl}.
\end{remark}

A standard calculation shows that Dirichlet energy increases by a factor at most $K$ under almost-smooth $K$-quasiconformal embeddings \cite[p.14]{AhlforsLectures}. As a consequence, extremal length increases by a factor at most $K$ under such maps. In particular, extremal length does not increase under conformal embeddings.

\begin{lemma} \label{lem:eldecreases}
Let $f : X \to Y$ be an almost-smooth $K$-quasi\-con\-for\-mal embedding between finite Riemann surfaces and let $\vphi \in \qd(X) \setminus \{0\}$. Then $$\el(f_*\mathcal{F}^h\vphi, Y) \leq K \el(\mathcal{F}^h\vphi , X).$$
\end{lemma}

\begin{proof}
Let $\zeta = \sigma+ i \tau$ be a conformal coordinate on $f(X)$ and let $z=x+iy$ be a conformal coordinate at $f^{-1}(\zeta)$. Since $f$ is $K$-quasiconformal, we have $$\mathrm{d}x\mathrm{d}y = \det(\mathrm{d}f^{-1}) \mathrm{d}\sigma\mathrm{d}\tau \geq K^{-1}\|\mathrm{d} f^{-1} \|^2\mathrm{d}\sigma\mathrm{d}\tau.$$
Let $v(z)= \im \int_{z_0}^z \sqrt \vphi$ so that $|\mathrm{d} v| = |\im \sqrt \vphi| = \mathcal{F}^h\vphi$. We compute
\begin{align*}
\dir(f_*|\mathrm{d}v|) & = \int_{f(X)}  |\nabla(v\circ f^{-1})(\zeta)|^2 \,\mathrm{d}\sigma\mathrm{d}\tau \\ & =  \int_{f(X)}|(\mathrm{d}f^{-1})(\nabla v)(f^{-1}(\zeta))|^2 \,\mathrm{d}\sigma\mathrm{d}\tau \\
& \leq  \int_{f(X)}\|\mathrm{d}f^{-1}\|^2 |(\nabla v)(f^{-1}(\zeta))|^2 \,\mathrm{d}\sigma\mathrm{d}\tau \\ 
& \leq K \int_{f(X)}\det(\mathrm{d}f^{-1}) |(\nabla v)(f^{-1}(\zeta))|^2 \,\mathrm{d}\sigma\mathrm{d}\tau \\
& = K \int_{X}\ |(\nabla v)(z)|^2 \,\mathrm{d}x\mathrm{d}y = K \dir(|\mathrm{d}v|).
\end{align*}
It follows that
$$
\el(f_*|\mathrm{d}v|, Y) \leq \dir(f_*|\mathrm{d}v|) \leq K \dir(|\mathrm{d}v|) = K \el(|\mathrm{d}v|,X),
$$
where the last equality holds by Theorem \ref{thm:QDsrealizeEL}.
 \end{proof}

\begin{remark}
It would be desirable to have this for any partial measured foliation $F$ on $X$, but one runs into a difficulty: it is not clear that $f_* G$ dominates $f_* F$ if $G$ dominates $F$. See Lemma \ref{lem:pushforward} for a similar statement which sounds obvious but is not entirely straightforward.
\end{remark}

The previous inequality is sharp as the case of Teichm\"uller embeddings illustrates.

\begin{lemma} \label{lem:elslit}
Let $f:X\to Y$ be a Teichm\"uller embedding of dilatation $K$ with initial and terminal quadratic differentials $\vphi$ and $\psi$. Then $f_*\mathcal{F}^h\vphi$ is measure equivalent to $\mathcal{F}^h\psi$ on $Y$ and we have $$\el(\mathcal{F}^h\psi,Y)=\el(f_*\mathcal{F}^h\vphi,Y)=K \el(\mathcal{F}^h\vphi,X).$$
\end{lemma}
\begin{proof}
We have $\mathcal{F}^h\psi=f_*\mathcal{F}^h\vphi$ on $f(X)$ since $f(x+iy)=Kx+iy$ in natural coordinates. Moreover, since $Y \setminus f(X)$ is a finite union of horizontal arcs and points, the integral of $\mathcal{F}^h\psi$ is zero along any piece of curve contained in this set. For any piecewise smooth curve $\alpha$ in $Y$, we thus have
\begin{align*}
\int_\alpha \mathrm{d}\mathcal{F}^h\psi & = \int_{\alpha} \chi_{f(X)} \mathrm{d}\mathcal{F}^h\psi + \int_{\alpha} \chi_{Y \setminus f(X)} \mathrm{d}\mathcal{F}^h\psi \\
& = \int_{\alpha} \mathrm{d}f_*\mathcal{F}^h\vphi + \int_{\alpha} \chi_{Y \setminus f(X)} \mathrm{d}\mathcal{F}^h\psi \\
& = \int_{\alpha} \mathrm{d}f_*\mathcal{F}^h\vphi
\end{align*}
where $\chi_A$ is the characteristic function of the set $A$. Therefore $\mathcal{F}^h\psi$ and $f_*\mathcal{F}^h\vphi$ are measure equivalent, which implies that they have the same extremal length. By Lemma \ref{lem:eldecreases}, the inequality
$$
\el(f_*\mathcal{F}^h\vphi,Y)\leq K \el(\mathcal{F}^h\vphi,X)
$$
holds and by Theorem \ref{thm:QDsrealizeEL} we have $\el(\mathcal{F}^h\psi,Y)=\|\psi\|$ and $\el(\mathcal{F}^h\vphi,X)=\|\vphi\|$.
We also have $\|\psi\|=K\|\vphi\|$ since $f$ stretches horizontally by a factor $K$ and $f(X)$ has full measure in $Y$. Putting everything together, we get
$$
\|\psi\|=\el(\mathcal{F}^h\psi,Y) = \el(f_*\mathcal{F}^h\vphi,Y) \leq K \el(\mathcal{F}^h\vphi,X) = K \|\vphi\| = \|\psi \|
$$
and hence equality holds.
 \end{proof}

We will need a sufficient condition for when the push-forwards of a measured foliation by two homotopic embeddings are measure equivalent. We say that an embedding $f: X \to Y$ between finite Riemann surfaces is \emph{tame} if it is almost-smooth and extends to a continuous map $\what f : \what X \to \what Y$  which is piecewise smooth along $\partial X$. For example, every Teichm\"uller embedding is tame.

\begin{lemma} \label{lem:pushforward}
Let $f,g : X \to Y$ be homotopic tame embeddings between finite Riemann surfaces and let $\vphi \in \qd^+(X)$. Then $f_*\mathcal{F}^h\vphi$ and $g_*\mathcal{F}^h\vphi$ are measure equivalent on $Y$.
\end{lemma}

It is important that we take $\vphi$ in $\qd^+(X)$ and not just in $\qd(X)\setminus \{0\}$. We want the transverse measure of each boundary component to be zero.

\begin{proof}
Let $F = \mathcal{F}^h\vphi$. By symmetry, it suffices to show that $f_*F$ dominates $g_*F$. Let $\alpha$ be a piecewise smooth closed curve or cross-cut in $Y$. Given $\eps > 0$, we have to find a closed curve or cross-cut $\beta$ homotopic to $\alpha$ such that
\begin{equation} \label{eq:measure}
\int_\beta g_*F \leq \int_\alpha f_* F + \eps.
\end{equation}
If $\alpha$ is contained in $f(X)$, then we can just take $\beta = g \circ f^{-1}(\alpha)$. The difficulty is when $\alpha$ intersects $Y \setminus f(X)$. 

If $\int_\alpha f_* F=\infty$ there is nothing to show, so assume the integral is finite.  Then there is only a finite number $n$ of components of $f^{-1}(\alpha)$ which are essential in $X$ in the sense that they cannot be homotoped into an arbitrarily small neighborhhood of $\partial X \cup \dot{X}$. Let $\delta = \eps/(8n)$ and let $X^\delta$ be $\what X$ minus a $\delta$-neighborhood of $\partial X \cup \dot{X}$ in the metric $|\vphi|$. The map $g\circ f^{-1}: f(X^\delta) \to g(X^\delta)$ can be extended to an almost-smooth homeomorphism $H : Y \to Y$ homotopic to the identity (see \cite[Lemma 2]{Masur}).

We will take $\beta$ to be a modified version of $H(\alpha)$. First, each inessential component of $H(\alpha) \cap g(X)$ can be homotoped within $g(X)$ to an arc contained in $g(X \setminus X^\delta)$ having arbitrarily small transverse measure with respect to $g_* F$. Even if there are infinitely many such inessential components, we can make sure that their total transverse measure is at most $\eps /2$. Let $\gamma$ be the curve obtained after having done this. 

Next, we homotope each essential component $c$ of $\gamma\cap g(X \setminus X^\delta)$, within $g(X \setminus X^\delta)$, to an arc having transverse measure at most $2 \delta$, while keeping endpoints fixed. This is possible since the height of each component of $X \setminus X^\delta$ is $\delta$. Moreover, no matter how many times $g^{-1}(c)$ winds around the annular component of $X \setminus X^\delta$ in which it is contained, we can push this winding part toward $\partial X \cup \dot{X}$. In doing so, the transverse measure of the winding part tends to zero and is thus eventually less than $\delta$, for a total of at most $2 \delta$. There are $2n$ such components---two for each essential component of $\gamma \cap g(X)$---for a total transverse measure of $(2n)(2\delta) = \eps / 2$.

Let $\beta$ be the curve obtained after having done these two modifications. Then inequality \eqref{eq:measure} holds since $\beta \setminus g(X)$ has measure zero, $\beta \cap g(X \setminus X^\delta) $ contributes at most $\eps/2+\eps/2 = \eps$, and $\beta \cap g(\overline{X^\delta})$ is the image by $g \circ f^{-1}$ of a subset of $\alpha$ so its transverse measure with respect to $g_*F$ is at most the transverse measure of $\alpha$ with respect to $f_* F$.
 \end{proof}

\section{The blob is semi-smooth} \label{sec:semismooth}

In this section, we show that each quadratic differential realizing a point on the boundary of $\blob(x,X,Y,h)$ determines a vector normal to the blob and vice versa. We use this to prove that the blob satisfies a certain regularity condition near its boundary which we call semi-smoothness. We first need the following formula for the derivative of extremal length \cite[Theorem 8]{Gardiner}.

\begin{theorem}[Gardiner's formula] \label{thm:formula}
Let $S$ be a finite Riemann surface and let $F$ be a measured foliation on $S$. Then the function\footnote{This is a slight abuse of notation. Recall that a point in $\teich(S)$ is an equivalence class of marking $f:S \to Z$, which we may assume is almost-smooth. By the expression $\el(F,Z)$ we really mean $\el(f_*F,Z)$.\label{ftnote:abuse}} $Z\mapsto \el(F,Z)$ is diffe\-ren\-tiable on $\teich(S)$. Its derivative at $Z \in \teich(S)$ in the direction $\mu \in T_Z \teich(S)$ is 
\begin{equation} \label{eq:GardFormula}
2 \langle \mu , \vphi \rangle
\end{equation}
where $\vphi \in \qd(Z)$ is the unique quadratic differential such that $\mathcal F^h \vphi$ is measure equivalent to $F$.
\end{theorem}

We apply this formula along a curve $Z_t= Y \setminus \{\gamma(t)\}$ for some analytic path $\gamma: I \to Y$, where $Y$ is a finite hyperbolic surface. It is implicit here that the change of marking $Z_0 \to Z_t$ is given by point-pushing along $\gamma$. In this case, the pairing $\langle \mu , \vphi \rangle$ is proportional to the real part of the residue of $\vphi$ in the direction of $\gamma'(0)$. See \cite{McMullenTwist} for a similar but more general calculation. 

\begin{lemma} \label{lem:green}
Let $Y$ be a finite hyperbolic surface, let $\gamma: (-1,1) \to Y$ be an analytic arc and let $Z_t=Y \setminus \{\gamma(t)\}$. Then the derivative $\mu$ of $Z_t$ at $t=0$ satisfies 
\begin{equation} \label{eq:pairing}
\langle \mu , q \rangle = -\pi \re[ \res_{\gamma(0)}(q \cdot \gamma'(0))]
\end{equation}
for every $q \in \qd(Z_0)$.
\end{lemma}
\begin{proof}
We may assume that $\gamma$ is embedded by restricting to a subinterval. Since $\gamma$ is analytic, the vector field $\gamma_*(\partial / \partial t)$ can be extended to a holomorphic vector field $\vec{v}$ on an embedded round disk $D \subset Y$ centered at $y=\gamma(0)$. Let $E$ be another round disk of smaller radius centered at $y$ and let $\phi$ be a smooth bump function which is equal to $1$ on $E$ and $0$ outside $D$. Then $\phi \mathbf{v}$ is a smooth vector field defined on all of $Y$. Let $\Phi_t$ be the time-$t$ flow for $\phi \mathbf{v}$ and let $\nu(t) = \overline{\partial} \Phi_t /\partial \Phi_t $. Then for small enough $t$, we have $\Phi_t(\gamma(0))=\gamma(t)$ by construction so that $\mu=\nu'(0) = \overline{\partial}(\phi \mathbf{v})$. We compute
\begin{align*}
 \int_{Z_0} \mu q &= \int_D q \overline{\partial}(\phi \mathbf{v}) = \int_{D\setminus E} q \overline{\partial}(\phi \mathbf{v}) = \int_{D\setminus E} \overline{\partial}(\phi q \mathbf{v})  \\ & = -\frac{i}{2} \int_{\partial(D\setminus E)} \phi q \mathbf{v} = \frac{i}{2} \int_{\partial E} q\mathbf{v} = -\pi \res_{y}(q\mathbf{v}).
\end{align*}
The equality from the first line to the second is by Green's theorem and the change of sign in the next equality comes from reversing orientation on $\partial E$. To conclude the proof, recall that $\langle \mu , q \rangle = \re \int_{Z_0} \mu q$ by definition.
 \end{proof}

Let $Y$ be a finite hyperbolic surface and let $y \in Y$. If $\psi \in \qd(Y\setminus \{y\})$ has a simple pole at $y$, then there exists a tangent vector $\mathbf{v} \in T_{y} Y$ (unique up to rescaling) such that $\res_y(\psi \mathbf{v}) < 0$. We say that $\mathbf{v}$ is \emph{vertical} for $\psi$. For example, if $\psi = \frac{1}{z}\mathrm{d} z^2$ then $\mathbf{v}= - \frac{\partial}{\partial z}$ is vertical at $0$ since 
$$\res_0(\psi \mathbf{v})= \frac{1}{2\pi i }\oint \frac{-1}{z}\,\mathrm{d}z = -1.$$ 
A vector $\vec{w}\in T_{y} Y$ is \emph{mostly vertical} for $\psi$ if $\re[\res_y(\psi \vec{w})] < 0$. If $\mathbf{v}$ is vertical and $\lambda \in \CC$ then $\lambda\mathbf{v}$ is mostly vertical if and only if $\re \lambda > 0$. Thus the mostly vertical vectors are those that make an angle less than $\pi/2$ with the vertical direction (here angles are measured with respect to the Riemann surface structure on $Y$, not the cone metric coming from $|\psi|$). 

\begin{figure}[htp] 
\centering
\includegraphics[scale=.7]{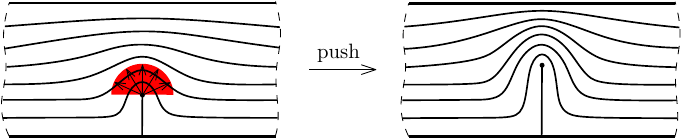}
\caption{Pushing a simple pole of a quadratic differential in a mostly vertical direction increases the extremal length of its horizontal foliation.}
\label{fig:pushup}
\end{figure}

If we push $y$ in a mostly vertical direction, then the extremal length of $\mathcal{F}^h \psi$ will increase according to equations \eqref{eq:GardFormula} and \eqref{eq:pairing}, at least for small time $|t|<\delta$. The intuition for this is that if we push in mostly vertical directions, we stretch the leaves of $\mathcal{F}^h\psi$ and compress them together, thereby increasing Dirichlet energy (see Figure \ref{fig:pushup}).

 We will need a more uniform statement where $\delta$ can be chosen independently of the direction. We can do this as long as we restrict to a compact set of angles. Given a point $y \in Y$ and a tangent vector $\vec{v} \in T_y Y$, let $\gamma_{\vec{v}} : \RR \to Y$ be the hyperbolic geodesic such that $\gamma_{\vec{v}}'(0)=\vec{v}$. We will push $y$ along these specific paths. 

\begin{lemma} \label{lem:ELgoesup}
Let $Y$ be a finite hyperbolic surface and let $y \in Y$. Suppose that $\psi \in \qd(Y\setminus \{y\})$ has a simple pole at $y$ and let $\vec{v} \in T_y Y$ be vertical for $\psi$. Then for every $\theta_0 \in (0, \frac{\pi}{2})$ there exists a number $\delta>0$ such that 
$$
 \el(\mathcal{F}^h \psi, Y \setminus \{ \gamma_{e^{i\theta}\mathbf{v}}(t)\}) > \el(\mathcal{F}^h \psi, Y \setminus \{ y \})
$$ 
for every $\theta \in [-\theta_0 , \theta_0]$ and every $t \in (0,\delta)$.
\end{lemma}
\begin{proof}
Suppose not and let $F =\mathcal{F}^h \psi$. Then there exist convergent sequences $t_n \searrow 0$ and $\theta_n \to \theta_\infty \in [-\theta_0 , \theta_0]$ as $n \to \infty$ such that 
$$\el(F, Y \setminus \{ \gamma_{e^{i\theta_n}\mathbf{v}}(t_n)\}) \leq \el(F, Y \setminus \{ y \})$$
for every $n \in \NN$.

For $\vec{w} \in T_y Y$ let $\mu_{\vec{w}}$ be the derivative of the path $t \mapsto Y \setminus \{ \gamma_{\vec{w}}(t)\}$ at $t=0$. Then $ \mu_{e^{i\theta_n}\mathbf{v}} \to \mu_{e^{i\theta_\infty}\mathbf{v}}$ as $n \to \infty$ by equation \eqref{eq:pairing}. Since $Z \mapsto \el(F, Z)$ is differentiable at $Z=Y \setminus \{ y \}$ we have that
$$
\frac{\el(F, Y \setminus \{ \gamma_{e^{i\theta_n}\mathbf{v}}(t_n)\}) -  \el(F, Y \setminus \{ y \})}{t_n } \to  2\langle \mu_{e^{i\theta_\infty}\mathbf{v}}, \psi \rangle 
$$
as $n \to \infty$. But the left-hand side is non-positive for each $n$ while
$$
2\langle \mu_{e^{i\theta_\infty}\mathbf{v}}, \psi \rangle  = -2\pi \re[ \res_y (\psi e^{i\theta_\infty}\mathbf{v}) ] = -2\pi \cos(\theta_\infty) \res_y (\psi \mathbf{v}) > 0.
$$
This is a contradiction.
\end{proof}

Let us introduce some more terminology. Given $z \in \CC$, $\mathbf{v}\in T_z \CC \setminus \{0\}$, $\theta \in (0,\pi)$, and $\delta>0$, we denote by $\sect(\mathbf{v},\theta,\delta)$ the open angular sector based at $z$ with radius $\delta$ and angle $\theta$ on either side of $\mathbf{v}$. In symbols,
$$
\sect(\mathbf{v},\theta,\delta) = \left\{\, z+te^{i\phi}\frac{\mathbf{v}}{|\mathbf{v}|} : \phi \in(-\theta,\theta) \text{ and } t \in (0, \delta ) \,\right\} .
$$

\begin{definition}
Let $B\subset \CC$ be closed. A vector $\mathbf{v} \in T_z \CC$ with $z\in B$ is \emph{normal to $B$} if $\mathbf{v}\neq 0$ and if there are angular sectors arbitrarily close to half-disks pointing in the direction of $\mathbf{v}$ which are disjoint from $B$. More precisely, $\mathbf{v}$ is normal to $B$ if $\mathbf{v}\neq 0$ and if for every $\theta \in (0,\frac{\pi}{2})$, there exists a $\delta >0$ such that $\sect(\mathbf{v},\theta,\delta)\cap B = \varnothing$. 
\end{definition}

\begin{figure}[htp]  \centering
\includegraphics[scale=1.25]{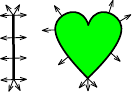}
\caption{A closed subset of the plane and some of its normal vectors. There is no normal vector at the inward corner and several normal vectors at outward corners.}
\end{figure}

Note that conformal homeomorphisms preserve normal vectors. Thus the definition extends to any Riemann surface. Furthermore, instead of Euclidean rays we can use geodesic rays for any conformal Riemannian metric (for example the hyperbolic metric if the surface is hyperbolic).

We are ready to prove that any quadratic differential realizing a point on the boundary of the blob determines a vector normal to the blob.

\begin{proposition}[Vertical vectors are normal] \label{prop:vertnormal}
Assume hypothesis \myhyperlink{h4}. Let $x \in X$, let $c \in \partial \blob(x,X,Y,h)$ and suppose that $\psi \in \qd^+(Y \setminus \{\pi_Y(c)\})$ realizes $c$. If $\mathbf{v}\in T_{c} \wtilde Y$ is vertical for $\pi_Y^* \psi$, then $\mathbf{v}$ is normal to $\blob(x,X,Y,h)$.
\end{proposition}
\begin{proof}
Let $f\in \emb(X,Y,h)$ be the slit mapping rel $x$ realizing $c$ provided by Proposition \ref{prop:bdryblob}. Let $y = \pi_Y(c) = f(x)$, let $g: X \setminus \{x\} \to Y \setminus \{y\}$ be the restriction of $f$, let $G = \mathcal{F}^h \psi$ on $Y \setminus \{y\}$, let $\vphi = g^* \psi$ and let $F = \mathcal{F}^h \vphi$ on $X \setminus \{x\}$. Also let $\vec{w} = \mathrm{d}\pi_Y(\vec{v})$ so that $\vec{w}$ is vertical for $\psi$ at $y$.

Fix $\theta_0 \in (0,\frac{\pi}{2})$. By Lemma \ref{lem:ELgoesup}, there exists a $\delta >0$ such that for every $\theta \in [-\theta_0,\theta_0]$ and every $t \in (0,\delta)$ we have
$$
\el(G, Y \setminus \{ \gamma_{e^{i\theta}\mathbf{w}}(t) \} ) > \el(G, Y \setminus \{y\} ).
$$
Since $g$ is a slit mapping, $G$ is measure equivalent to $g_* F$ on $Y \setminus \{y\}$ by Lemma \ref{lem:elslit}. We also have $\el(G, Y \setminus\{y\} )= \|\psi\|=\|\vphi\|=\el(F, X \setminus \{x\} )$ by Theorem \ref{thm:QDsrealizeEL}, so the above inequality translates to
\begin{equation} \label{eq:abuse}
\el(g_* F, Y \setminus \{ \gamma_{e^{i\theta}\mathbf{w}}(t) \} ) > \el(F, X \setminus \{x\}).
\end{equation}

Fix $\theta \in [-\theta_0,\theta_0]$ and $t \in (0,\delta)$. Let $\push$ be a quasiconformal point-pushing diffeomorphism along $\gamma_{e^{i\theta}\mathbf{w}}([0,t])$. Recall that the left-hand side of \eqref{eq:abuse} really stands for the extremal length of $(\push \circ g)_* F$ since $\push : Y \setminus \{y\} \to Y \setminus \{ \gamma_{e^{i\theta}\mathbf{w}}(t) \}$ is our change of marking by convention (see footnote \ref{ftnote:abuse}). 

Suppose that there is a conformal embedding $M: X \setminus \{x\} \to Y \setminus \{ \gamma_{e^{i\theta}\mathbf{w}}(t) \}$ homotopic to $\push \circ g$. Then by Corollary \ref{cor:tameapprox}, either $M$ is a slit mapping or it is homotopic to a conformal embedding which extends analytically to the ideal boundary. In either case, $M$ is homotopic to a tame conformal embedding $N$. By Lemma \ref{lem:pushforward}, $N_*F$ and $(\push \circ g)_*F$ are measure equivalent on $Y \setminus \{ \gamma_{e^{i\theta}\mathbf{w}}(t) \}$ which means that
$$
\el(N_*F, Y \setminus \{ \gamma_{e^{i\theta}\mathbf{w}}(t) \}) = \el((\push \circ g)_*F, Y \setminus \{ \gamma_{e^{i\theta}\mathbf{w}}(t)\}) > \el(F, X \setminus \{x\}).
$$
On the other hand, $\el(N_*F, Y \setminus \{ \gamma_{e^{i\theta}\mathbf{w}}(t) \}) \leq \el(F, X \setminus \{x\})$ by Lemma \ref{lem:eldecreases}. We conclude that no such $M$ exists. 

Equivalently, the concatenation of $c$ with $\gamma_{e^{i\theta}\mathbf{w}}([0,t])$ does not belong to $\blob(x,X,Y,h)$. Since this holds for every $\theta \in [-\theta_0,\theta_0]$ and every $t \in (0,\delta)$, the open hyperbolic sector of radius $\delta\|\vec{v}\|_{\wtilde Y}$ (the norm is with respect to the hyperbolic metric) and angle $\theta_0$ on either side $\vec{v}$ in $\wtilde Y$ is disjoint from $\blob(x,X,Y,h)$. Thus $\vec{v}$ is normal to $\blob(x,X,Y,h)$, as $\theta_0 \in (0, \frac{\pi}{2})$ was arbitrary.
 \end{proof}

As a corollary, we obtain the converse of Proposition \ref{prop:bdryblob}.

\begin{corollary} \label{cor:blobcomplement}
Assume hypothesis \myhyperlink{h4} and let $x \in X$. If $f:X \to Y$ is a slit mapping rel $x$ homotopic to $h$, then $\lift_x(f) \in \partial \blob(x,X,Y,h)$.
\end{corollary}
\begin{proof}
Let $c=\lift_x(f)$. Since $f$ is conformal, $c$ belongs to $\blob(x,X,Y,h)$. Let $\psi$ be a terminal quadratic differential for $f$ rel $x$. If $\psi$ does not have a simple pole at $y=f(x)$ then $f$ is a slit mapping, contradicting the assumption that $\emb(X,Y,h)$ does not contain any. Thus $\psi$ has a simple pole at $y$ so that there exists a vector $\mathbf{v} \in T_{c} \wtilde Y$ which is vertical for $\pi_Y^*\psi$. We can apply the same reasoning as in the proof of Proposition \ref{prop:vertnormal} to conclude that $\vec{v}$ is normal to $\blob(x,X,Y,h)$. In particular, $c \in \partial \blob(x,X,Y,h)$.
 \end{proof}

We will also prove the converse of Proposition \ref{prop:vertnormal}, namely that every vector normal to $\blob(x,X,Y,h)$ is vertical with respect to some realizing quadratic differential. The first step is to strengthen Lemma \ref{lem:ELgoesup} by allowing the point $y \in Y$ and the quadratic differential $\psi$ to vary inside a compact set. We have to be a little more careful about how to compare two surfaces $Y \setminus \{ y_1\}$ and $Y \setminus \{ y_2\}$ for $y_1,y_2 \in Y$ though.

\begin{remark} \label{rem:theta}
Let $S$ be a finite hyperbolic surface and let $s \in S$. Given $Y \in \teich(S)$ and $y \in Y$, there is no canonical way to define a marking $S \setminus \{s\} \to Y \setminus \{y\}$. However, if $b \in Y$ is some basepoint, $\wtilde Y$ is the universal cover of $Y$ based at $b$, and $f: S \setminus \{s\} \to Y \setminus \{b\}$ is a fixed marking, then point-pushing provides a continuous map $\Theta : \wtilde Y \to \teich(S \setminus \{s\})$ as follows. We can represent any $c \in \wtilde Y$ by a smooth immersed curve $\gamma : [0,1] \to Y$ with $\gamma(0)=b$. Let $\push_\gamma$ be a quasiconformal point-pushing diffeomorphism along $\gamma$, then define $\Theta(c)$ as the marked surface $\push_\gamma \circ f : S \setminus \{s\} \to Y \setminus \{\pi_Y(c)\}$ in $\teich(S \setminus \{s\})$. We use these markings implicitly below.
 \end{remark}

\begin{lemma} \label{lem:deltaupper}
Let $Y$ be a finite hyperbolic surface. Let $\theta_0 \in (0, \frac{\pi}{2})$ and let $K$ be a compact set of pairs $(c, \psi)$ where $c\in \wtilde Y$ and $\psi \in \qd(Y \setminus \{\pi_Y(c)\})$ has a simple pole at $\pi_Y(c)$. There exists a number $\delta>0$ such that $\forall (c,\psi) \in K$, $\forall \theta \in [-\theta_0,\theta_0]$ and $\forall t \in (0,\delta)$ we have
$$
\el(\mathcal{F}^h \psi, Y \setminus \{\gamma_{e^{i\theta} \vec{v}}(t)\}) > \el(\mathcal{F}^h \psi, Y \setminus \{\pi_Y(c)\}) 
$$ 
where $\vec{v}$ is the vertical vector for $\psi$ at $\pi_Y(c)$ rescaled to have norm $1$ with respect to the hyperbolic metric.
\end{lemma}
\begin{proof}
The idea of the proof is the same as for Lemma \ref{lem:ELgoesup}. If the statement fails, then there exist convergent sequences $t_n \searrow 0$, $(c_n, \psi_n) \to (c,\psi) \in K$, $\theta_n \to \theta \in [-\theta_0,\theta_0]$ such that 
$$
\el(\mathcal{F}^h \psi_n, Y \setminus \{\gamma_{e^{i\theta_n} \vec{v}_n}(t_n)\}) \leq \el(\mathcal{F}^h \psi_n, Y \setminus \{\pi_Y(c_n)\})
$$ 
for every $n \in \NN$, where $\vec{v_n}$ is the unit vertical vector for $\psi_n$ at $\pi_Y(c_n)$. Since $\psi_n \to \psi$ as $n \to \infty$ we have that $\vec{v_n} \to \vec{v}$ where $\vec{v}$ is the unit vertical vector for $\psi$ at $\pi_Y(c)$. If $\mu_n$ denotes the derivative of the path $t \mapsto Y \setminus \{\gamma_{e^{i\theta_n} \vec{v}_n}(t)\}$ at $t=0$ and $\mu$ the derivative of $t \mapsto Y \setminus \{\gamma_{e^{i\theta} \vec{v}}(t)\}$ at $t=0$, then $\mu_n \to \mu$ as $n \to \infty$ since $\theta_n \to \theta$ and $\vec{v}_n \to \vec{v}$. In the limit, we have
$$
2\langle \mu , \psi \rangle = -2\pi \re[\res_{\pi_Y(c)}(\psi e^{i\theta}\vec{v})] = -2\pi \cos(\theta)\re[\res_{\pi_Y(c)}(\psi\vec{v})] > 0.
$$
The difference with Lemma \ref{lem:ELgoesup} is that the measured foliations $\mathcal{F}^h \psi_n$ are not constant, so we are looking at a sequence of functions rather than a single one.

Let $b\in Y$ be our basepoint for the universal cover $\pi_Y : \wtilde Y \to Y$. Recall that any $c\in \wtilde Y$ determines a marking $f_c : Y \setminus \{b\} \to Y \setminus \{\pi_Y(c)\}$. Thus the measured foliations $F_n = f_{c_n} ^* (\mathcal{F}^h \psi_n)$ and $F=f_c^*(\mathcal{F}^h \psi)$ all live on the same surface $S = Y \setminus \{b\}$. Moreover $F_n \to F$ in measure since $(c_n,\psi_n) \to (c,\psi)$ as $n \to \infty$. Given a measured foliation $G$ on $S$ and $Z \in \teich(S)$, let $\Phi(G,Z)$ be the quadratic differential on $Z$ whose horizontal foliation is measure equivalent to $G$. By Hubbard--Masur \cite{HubbardMasur}, the map $\mathcal{MF}(S) \times \teich(S) \to T^* \teich(S)$ sending $(G,Z)$ to $\Phi(G,Z)$ is continuous. Thus the function $$Z \mapsto \Lambda_G(Z) := \el(G, Z) = \|\Phi(G,Z)\|$$ depends continuously on $G$. Furthermore, its derivative $$(Z,\nu) \mapsto \mathrm{d}_Z\Lambda_G (\nu) = 2 \langle \nu , \Phi(G,Z) \rangle$$ is continuous on $T \teich(S)$ and depends continuously on $G$.

Let $\{x_n\}_{n=1}^\infty$ and $\{y_n\}_{n=1}^\infty$ be sequences in $\RR^k$ such that $x_n \to 0$, $y_n \to 0$ and $\frac{y_n - x_n}{|y_n - x_n|} \to v$ as $n \to \infty$. Suppose that $g_n,g : \RR^k \to \RR$ are $C^1$ functions such that $g_n \to g$ and $\mathrm{d} g_n \to \mathrm{d}g$ uniformly on compact sets as $n \to \infty$. Then
$$
\frac{g_n(y_n)-g_n(x_n)}{|y_n - x_n|} \to \mathrm{d}g(v)
$$
as $n \to \infty$. One can show this using the fundamental theorem of calculus along the line segment between $x_n$ and $y_n$ for instance. By the previous paragraph, this criterion applies to the extremal length functions $\Lambda_{F_n}$ and $\Lambda_{F}$ on $\teich(S)$ (which is diffeomorphic to some $\RR^k$) to conclude that
$$
\frac{\el(\mathcal{F}^h \psi_n, Y \setminus \{\gamma_{e^{i\theta_n} \vec{v}_n}(t_n)\}) - \el(\mathcal{F}^h \psi_n, Y \setminus \{\pi_Y(c_n)\})}{t_n} \to 2\langle \mu , \psi \rangle 
$$
as $n \to \infty$. This is a contradiction, since left-hand side is non-positive for every $n \in \NN$ while the right-hand side is positive.
 \end{proof}

The next observation that we need is that that any non-zero limit of a sequence of terminal quadratic differentials is a terminal quadratic differential for the limiting Teichm\"uller embedding.

\begin{lemma} \label{lem:limqd}
Let $h : X \to Y$ be a generic embedding between finite Riemann surfaces and let $x \in X$. Let $f_n: X \to Y$ be a Teichm\"uller embedding rel $x$ homotopic to $h$ with unit norm terminal quadratic differential $\psi_n$. Suppose that $\lift_x(f_n)$ converges to some $c \in \wtilde Y$. Then there is a Teichm\"uller embedding $f$ rel $x$ homotopic to $h$ with unit norm terminal quadratic differential $\psi$ such that $f_n \to f$ and $\psi_n \to \psi$ as $n \to \infty$ after passing to a subsequence.
\end{lemma}
\begin{proof}
Let $K_n$ be the dilatation of $f_n$. By hypothesis, the marked surface $\Theta(\lift_x(f_n))$ converges to $\Theta(c)$ in the Teichm\"uller space $\teich$ of $Y$ minus a point, where $\Theta : \wtilde Y \to \teich$ is the map from Remark \ref{rem:theta}. In particular, the sequence $K_n$ is bounded so we may assume it converges to some $K \geq 1$. By Lemma \ref{lem:limteich}, there is a Teichm\"uller embedding $f$ rel $x$ homotopic to $h$ such that $f_n \to f$ after passing to a subsequence.

Let $c_n = \lift_x(f_n)$. The set $\{c\} \cup\{c_n\}_{n=1}^\infty$ is compact in $\wtilde Y$ and so is its image by $\Theta$ in $\teich$. Since the set of unit cotangent vectors over a compact set in $\teich$ is compact, there is a quadratic differential $\psi \in \qd(Y \setminus \pi_Y(c))$ such that $\psi_n \to \psi$ after passing to a subsequence. We have to show that $\psi$ is a terminal quadratic differential for $f$.

Let $\vphi_n$ be the initial quadratic differential of $f_n$ corresponding to $\psi_n$. We may assume that $\vphi_n$ converges to some $\vphi \in \qd^+(X \setminus \{x\})$ since its norm $1/K_n$ is bounded above and below. Suppose that $z_0 \in X \setminus\{ x \}$ is not a zero of $\vphi$. Then there is a compact simply connected neighborhood $U$ of $z_0$ on which $\vphi$ does not vanish. If $n$ is large enough, then $\vphi_n$ does not have any zeros in $U$ either. If $V = f(U)$, then $\psi$ and $\psi_n$ do not have zeros in $V$ when $n$ is large enough. We can choose square roots consistently so that $\sqrt{\vphi_n} \to \sqrt{\vphi}$ uniformly on $U$ and $\sqrt{\psi_n} \to \sqrt{\psi}$ uniformly on $V$. Then for every $z \in U$ and every $n$ we have
$$
\int_{f_n(z_0)}^{f_n(z)} \sqrt{\psi_n} = K_n \re\left( \int_{z_0}^z \sqrt{\vphi_n} \right) + i \im\left( \int_{z_0}^z \sqrt{\vphi_n} \right)
$$
since $f_n$ restricted to $X \setminus \{x\}$ is a Teichm\"uller embedding of dilatation $K_n$ with respect to $\vphi_n$ and $\psi_n$. Taking the limit as $n \to \infty$ we get
$$
\int_{f(z_0)}^{f(z)} \sqrt{\psi} = K \re\left( \int_{z_0}^z \sqrt{\vphi} \right) + i \im\left( \int_{z_0}^z \sqrt{\vphi} \right)
$$
which means that $f$ is locally of Teichm\"uller form with respect to $\vphi$ and $\psi$.
 \end{proof}

Lastly, we will need the fact that the set of vectors which are vertical for some realizing quadratic differential at a given point is convex.

\begin{lemma} \label{lem:vertconvex}
Assume hypothesis \myhyperlink{h4}. Let $x \in X$, let $c \in \partial \blob(x,X,Y,h)$ and suppose that $\psi_0,\psi_1 \in \qd^+(Y\setminus \{\pi_Y(c)\})$ realize $c$. If $\mathbf{v}_0$ and $\mathbf{v}_1$ are vertical for $\psi_0$ and $\psi_1$ respectively at $y=\pi_Y(c)$, then there exist $\alpha,\beta > 0$ such that $\mathbf{v}_0+\mathbf{v}_1$ is vertical for $\alpha \psi_0 + \beta \psi_1$.
\end{lemma}
\begin{proof}
Take $\alpha = - \frac{|\mathbf{v}_0|}{|\mathbf{v}_1|} \res_y(\psi_1 \mathbf{v}_1)$ and $\beta = - \frac{|\mathbf{v}_1|}{|\mathbf{v}_0|} \res_y(\psi_0 \mathbf{v}_0)$. A calculation shows that $\res_y((\alpha \psi_0 + \beta \psi_1)(\mathbf{v}_0+\mathbf{v}_1))\leq 0$. By Lemma \ref{lem:qdconvex}, the quadratic differential $\alpha \psi_0 + \beta \psi_1$ has a simple pole at $y$. This implies that $\mathbf{v}_0+\mathbf{v}_1\neq 0$ and hence that $\res_y((\alpha \psi_0 + \beta \psi_1)(\mathbf{v}_0+\mathbf{v}_1))<0$.
 \end{proof}

We are now able to show that normal vectors are vertical.

\begin{proposition}[Normal vectors are vertical]\label{prop:normaldirections}
Assume \myhyperlink{h4} and let $x \in X$. Suppose that $\mathbf{v}$ is normal to $\blob(x,X,Y,h)$ at some point $c$. Then there exists a quadratic differential $\psi \in \qd^+(Y \setminus \pi_Y(c))$ realizing $c$ such that $\mathbf v$ is vertical for $\pi_Y^*\psi$.
\end{proposition}

\begin{proof}
Let $V_c \subset T_c \wtilde Y$ denote the set of vectors which are vertical for the pull-back of some quadratic differential realizing $c$. By Lemma \ref{lem:qdconvex} and Lemma \ref{lem:vertconvex}, $V_c$ is convex. Moreover $V_c \cup \{0_c\}$ is closed by Lemma \ref{lem:limqd}. Suppose that $\mathbf v$ is not in $V_c$. Then there is an open half-plane $H$ through the origin in $T_c \wtilde Y$ containing $V_c$ such that $\mathbf v$ is not in the closure $\overline H$. Let $c_n$ be a sequence converging to $c$ along the geodesic ray $r$ which is normal to $\overline H$ at $c$. Since $r$ makes an angle strictly less than $\pi /2$ with $\mathbf v$ and since $\mathbf v$ is normal to $\blob(x,X,Y,h)$, we may assume that $c_n$ is not in $\blob(x,X,Y,h)$ for any $n$. Let $f_n$ be the Teichm\"uller embedding rel $x$ realizing $c_n$ provided by Lemma \ref{lem:bdryblobteich} and let $\psi_n$ be its terminal quadratic differential, normalized to have norm $1$.

We may assume that $\psi_n$ converges to some $\psi \in \qd^+(Y \setminus \pi_Y(c))$ and that $f_n$ converges to a slit mapping $f$ rel $x$ with terminal quadratic differential $\psi$ by Lemma \ref{lem:limqd}. Indeed, the limit $f$ is conformal since $c$ belongs to $\blob(x,X,Y,h)$. If $\psi$ is holomorphic at $\pi_Y(c)$, then $f$ is a slit mapping on $X$, contrary to the assumption that $\emb(X,Y,h)$ does not contain any.  Therefore $\psi$ has a simple pole at $\pi_Y(c)$, which implies that $\psi_n$ has a simple pole at $\pi_Y(c_n)$ for all but finitely many indices. Let $\mathbf{w}_n$ be the unit vertical vector for $\pi_Y^* \psi_n$ at $c_n$. Then $\mathbf{w}_n$ converges to the unit vertical vector $\mathbf{w}$ for $\pi_Y^*\psi$ at $c$. We have $\mathbf{w} \in V_c \subset H$. We will see that this yields a contradiction.

Let $\phi\in(0,\frac{\pi}{2}]$ be the angle between $\mathbf{w}$ and the line $\partial H$, and let $\theta = \frac{\pi}{2}-\frac{\phi}{2}$. By Lemma \ref{lem:deltaupper}, there exists a number $\delta>0$ such that $\forall n\in \NN$, $\forall \alpha \in [-\theta, \theta]$ and $\forall t \in (0, \delta)$ we have
$$
\el(\mathcal{F}^h \psi_n, Y \setminus \{\pi_Y(\gamma_{e^{i\alpha}\vec{w}_n}(t))\}) > \el(\mathcal{F}^h \psi_n, Y \setminus \{\pi_Y(c_n)\}).
$$
On the other hand, since $\mathbf{w}_n \to \mathbf{w}$, the angle between $\mathbf{w}_n$ and the geodesic between $c_n$ and $c$ converges to the angle between $\vec{w}$ and the ray opposite to $r$, that is, to $\frac{\pi}{2} - \phi < \theta$. Thus if $n$ is large enough then $c$ belongs to the open hyperbolic sector of radius $\delta$ and angle $\theta$ on either side of $\vec{w}_n$ so that
$$
\el(\mathcal{F}^h \psi_n, Y \setminus \{\pi_Y(c)\}) > \el(\mathcal{F}^h \psi_n, Y \setminus \{\pi_Y(c_n)\}).
$$

\begin{figure}[htp]  \centering
\includegraphics[scale=.95]{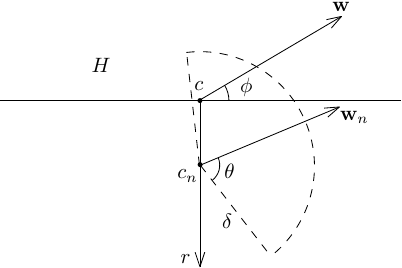}
\caption{The point $c$ is eventually contained in the sector of angle $\theta$ and radius $\delta$ about the vector $\mathbf{w}_n$.}
\label{fig:sector}
\end{figure}

Let $\vphi_n$ be the initial quadratic differential of $f_n$ corresponding to $\psi_n$ and let $K_n>1$ be the dilatation of $f_n$. By Lemma \ref{lem:elslit}, $\mathcal{F}^h\psi_n$ is measure equivalent to $(f_n)_* \mathcal{F}^h \vphi_n$ on $Y \setminus \pi_Y(c_n)$. If we let $\push_n : Y\setminus \{ \pi_Y(c_n)\} \to Y\setminus \{ \pi_Y(c)\}$ be the change of marking given by point-pushing, we get
\begin{align*}
\el((\push_n \circ f_n)_* \mathcal{F}^h \vphi_n, Y \setminus \{\pi_Y(c)\}) &= \el((\push_n)_* \mathcal{F}^h \psi_n, Y \setminus \{\pi_Y(c)\})\\
& > \el(\mathcal{F}^h\psi_n, Y \setminus \{\pi_Y(c_n)\}) \\
& = K_n \el(\mathcal{F}^h \vphi_n, X \setminus \{x \})\\
& > \el(\mathcal{F}^h \vphi_n ,X \setminus \{x\})
\end{align*}
when $n$ is large enough, by the above inequality and Lemma \ref{lem:elslit}. However, since $f$ and $\push_n \circ f_n$ are homotopic tame embeddings from $X \setminus \{x\}$ to $Y \setminus \{\pi_Y(c)\}$ and $f$ is conformal, we have
\begin{align*}
\el((\push_n \circ f_n)_* \mathcal{F}^h \vphi_n, Y \setminus \{\pi_Y(c)\}) &= \el(f_* \mathcal{F}^h \vphi_n, Y \setminus \{\pi_Y(c)\}) \\ & \leq \el(\mathcal{F}^h \vphi_n ,X \setminus \{x\})
\end{align*}
by Lemma \ref{lem:pushforward} and Lemma \ref{lem:eldecreases}. This is a contradiction, from which we conclude that $\mathbf v$ belongs to $V_c$.
 \end{proof}

We finally come to the main result of this section, which is that the blob is semi-smooth.

\begin{definition}
A closed subset $B\subset \CC$ is \emph{semi-smooth} if 
\begin{itemize}
\item for every  $z \in \partial B$, the set of normal vectors to $B$ at $z$ is non-empty and convex;
\item any non-zero limit of vectors normal to $B$ (not necessarily at the same point) is normal to $B$.
\end{itemize}
\end{definition}

For example, a $2$-dimensional manifold $M\subset \CC$ with $C^1$ boundary is semi-smooth, as is any convex set with non-empty interior. A line segment is not semi-smooth because at interior points the set of normal vectors is not convex. 

\begin{theorem}[The blob is semi-smooth] \label{thm:blobsemismooth}
Assume hypothesis \myhyperlink{h4}. Then $\blob(x,X,Y,h)$ is semi-smooth for any $x \in X$.
\end{theorem}
\begin{proof}
For every $c \in \partial \blob(x,X,Y,h)$ the set of vectors which are normal to $\blob(x,X,Y,h)$ at $c$ coincides with the set $V_c$ of vectors which are vertical for the pull-back of some quadratic differential realizing $c$, according to Propositions \ref{prop:vertnormal} and \ref{prop:normaldirections}. The set $V_c$ is non-empty by Proposition \ref{prop:bdryblob} and convex by Lemma \ref{lem:qdconvex} and Lemma \ref{lem:vertconvex}. Suppose that $\{c_n\}_{n=1}^\infty \subset \partial \blob(x,X,Y,h)$ is such that $c_n \to c$, that  $\mathbf{v}_n$ is vertical for $\pi_Y^* \psi_n$ where $\psi_n$ realizes $c_n$, and that $\mathbf{v}_n \to \mathbf{v} \neq 0$ as $n \to \infty$. By Lemma \ref{lem:limteich} we can rescale $\psi_n$ and pass to a subsequence such that it converges to some $\psi$ realizing $c$. We have $\res_c((\pi_Y^*\psi)\mathbf{v}) \leq 0$ since $\res_{c_n}((\pi_Y^*\psi_n)\mathbf{v}_n) < 0$ for every $n$. Moreover, $\psi$ must have a simple pole at $\pi_Y(c)$ for otherwise $\emb(X,Y,h)$ would contain a slit mapping. This means that  $\res_c((\pi_Y^*\psi)\mathbf{v}) \neq 0$. Therefore $\mathbf{v}$ is vertical for $\pi_Y^*\psi$, hence normal to $\blob(x,X,Y,h)$ by Proposition \ref{prop:vertnormal}.
 \end{proof}

\section{The blob is a disk}  \label{sec:noholes}

In this section, we complete the proof of Theorem \ref{thm:blobisdisk} which says that the blob is homeomorphic to a closed disk under hypothesis \myhyperlink{h4}. We first prove that every semi-smooth set is a manifold.

\begin{theorem}[The aquatic theorem] \label{thm:aquatic}
Every closed semi-smooth subset of $\CC$ is a $2$-dimensional manifold with boundary.
\end{theorem}

\begin{proof}
Let $B$ be a closed semi-smooth set. Every interior point of $B$ has a neighborhood homeorphic to an open subset of $\CC$, namely the interior of $B$. Thus we only have to show that every boundary point $z\in \partial B$ has a neighborhood homeomorphic to a half-disk. By applying an isometry of the plane, we may assume that $z=0$ and that the vector $i$ bissects the cone $N_0$ of vectors normal to $B$ at $0$. Let $\phi$ be half the angle of $N_0$, let $\alpha = \phi + \frac{\pi}{2}$ and let $\beta = \pi - \alpha$. 

\begin{figure}[htp]  \centering
\includegraphics[scale=.9]{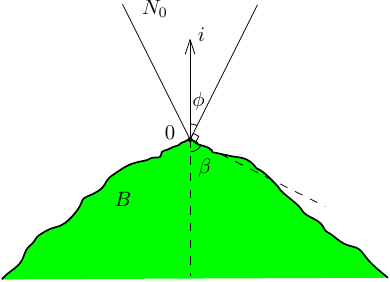}
\caption{The cone of normal vectors $N_0$ and the angles $\alpha$ and $\beta$.}
\label{fig:normalvectors}
\end{figure}

Since $B$ is semi-smooth, we have $\phi<\frac{\pi}{2}$ and hence $\beta > 0$. Moreover, for every $\theta \in (0,\alpha)$ there exists a $\delta>0$ such that the open angular sector $\sect(i,\theta,\delta)$ is disjoint from $B$ (notation is as in the previous section). We now show the existence of closed sectors pointing downwards contained in $B$.

\begin{claim} \hypertarget{clm:bubbles}{}
For every $\theta \in (0,\beta)$, there exists a $\delta > 0$ such that the closed angular sector $\overline{\sect(-i,\theta,\delta)}$ is contained in $B$.
\end{claim}
\begin{proof}[Proof of Claim]
Suppose not. Then there exists a $\theta \in (0,\beta)$ and a sequence $\delta_n \searrow 0$ for which the corresponding angular sector $S_n=\overline{\sect(-i,\theta,\delta_n)}$ intersects the complement of $B$ for every $n$. Let $D_n$ be a closed disk in $S_n$ disjoint from $B$. Slide the center of $D_n$ in a straight line towards $0$ until the disk first hits $B$, and let $D_n^*$ be the resulting disk. The intersection points of $D_n^*$ with $B$ all lie on the half of $\partial D_n^*$ which is closest to $0$. Let $z_n$ be any point in this intersection. Then $z_n$ is on the boundary of $B$ and the unit vector $\mathbf{v}_n$ pointing from $z_n$ to the center of $D_n^*$ is normal to $B$. Since $S_n$ shrinks to $0$, we have $z_n \to z$. Each vector $\mathbf{v}_n$ makes an angle at most $\theta +\frac{\pi}{2}$ with the downward direction. Therefore, the vectors $\mathbf{v}_n$ can only accumulate onto vectors forming an angle at least $\beta - \theta$ with the cone $N_0$. This contradicts the hypothesis that every limit of normal vectors is normal. 
 \end{proof}

\begin{figure}[htp]  \centering
\includegraphics[scale=.9]{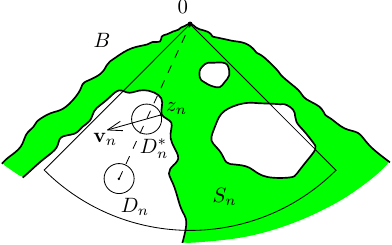}
\caption{Bubbles floating to the surface of $B$.}
\label{fig:bubbles}
\end{figure}

Let $\theta^+ \in (0, \alpha)$, let $S^+=\sect(i,\theta^+,\delta^+)$ be disjoint from $B$, let $\theta^- \in(0,\beta )$, and let $S^-=\overline{\sect(-i,\theta^-,\delta^-)}$ be contained in $B$. Let $I \subset S^+ $ be a compact horizontal segment symmetric about the vertical line through $0$ and lying entirely above $S^-$. We define a map $p : I \to \partial B$ as follows. For $z \in I$, let $z$ fall straight down until it first hits $B$, and let $p(z)$ be this first hitting point. Note that $p(x+iy)=x+iq(x,y)$ for some function $q$ so that $p$ is injective.

\begin{claim}
The map $p$ is continuous on some subinterval $J \subset I$ centered at the midpoint of $I$.
\end{claim}
\begin{proof}[Proof of Claim]
It is easy to see that $p$ is continuous at the midpoint $p^{-1}(0)$. This is because $p$ keeps the $x$-coordinate unchanged and moreover, $p(z)$ is below $S^+$ and above $S^-$. Thus the $y$-coordinate of $p(z)$ converges to $0$ as $z \to p^{-1}(0)$.
 
\begin{figure}[htp]  \centering
\includegraphics[scale=.9]{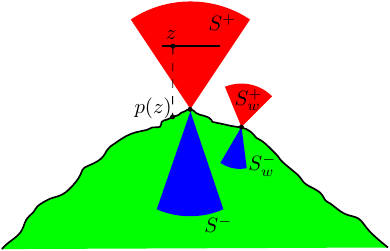}
\caption{The vertical projection $p$ is akin to rain falling on the surface of $B$. The projection is continuous (in fact Lipschitz) by the squeeze theorem.}
\label{fig:rain}
\end{figure} 
 
Let $0<\eps < \beta /2$. By semi-smoothness, there exists a $\delta > 0$ such that if $w\in \partial B$ is within distance $\delta$ of $0$, then every vector in $N_w$ is within angle $\phi + \eps$ of the upward direction. For every $w \in \partial B$ with $|w|<\delta$, let $\mathbf{v}_w$ be the bisector of $N_w$ and let $\phi_w$ be half the angle of $N_w$. For every $\theta_w^+ \in (0,\phi_w+\frac{\pi}{2})$ there is an open sector $S_w^+=\sect(\mathbf{v}_w,\theta_w^+, \delta^+_w)$ disjoint from $B$ by definition of $N_w$. Since
$$
\phi+\eps = \alpha - \frac{\pi}{2} + \eps < \alpha + \frac{\beta - \pi}{ 2}  = \frac{\alpha}{2} < \frac{\pi}{2} \leq \phi_w+\frac{\pi}{2},
$$
we may choose $\theta_w^+$ so that $S_w^+$ contains the vertical direction in its span. By the previous \hyperlink{clm:bubbles}{Claim}, there is also a closed sector $S_w^-=\overline{\sect(-\mathbf{v}_w,\theta_w^-,\delta_w^-)}$ contained in $B$ for every $\theta_w^- \in (0, \frac{\pi}{2} - \phi_w)$. The angle that $-\mathbf{v}_w$ makes with the downward direction is equal to the angle that $\mathbf{v}_w$ makes with the vertical direction, which is at most $\phi+\eps-\phi_w$ hence strictly less than $\frac{\pi}{2} - \phi_w$. Thus we may choose $\theta_w^-$ so that $S_w^-$ contains the downward direction in its interior.

By continuity of $p$ at $p^{-1}(0)$, there is a closed interval $J \subset I$ centered at $p^{-1}(0)$ such that $p(J)$ is contained in the ball of radius $\delta$ about $0$. Let $z \in J$, let $w=p(z)$, and let $S_w^+$ and $S_w^-$ be angular sectors as described in the previous paragraph. Also let $K \subset S_w^+$ be a compact horizontal segment crossing the vertical line through $w$ and lying entirely above $S_w^-$. By construction, the vertical segment from $z$ to $w$ intersects $B$ only at $w$. Since $B$ is closed, we may assume that the rectangle with bottom edge $K$ and upper edge $L\subset I$ is disjoint from $B$, by making $K$ shorter if necessary. For every $\zeta \in L$, the image $p(\zeta)$ is thus squeezed between $S_w^+$ and $S_w^-$, so that $p$ is continuous at $z$.
 \end{proof}

Thus $p(J)$ is the graph of a continuous function. Let $0<\delta<|J|/2$. For every $z \in J$ with $|x| < \delta$, draw the open vertical segment of length $2 \delta$ centered at $p(z)$, and let $U_\delta$ be the union of those segments. The continuity of $p$ implies that $U_\delta$ is open.

\begin{claim}
If $\delta$ is small enough, then the component of $U_\delta \setminus p(J)$ above $p(J)$ is disjoint from $B$ and the component below $p(J)$ is contained in $B$. 
\end{claim}
\begin{proof}[Proof of Claim]
If $\delta$ is small enough, then the component $C^+$ of $U_\delta \setminus p(J)$ above $p(J)$ lies below $J$ itself. By definition of $p$, for every $z \in J$ the open vertical segment between $z$ and $p(z)$ is disjoint from $B$, so that $C^+$ is disjoint from $B$.

For the component lying below $p(J)$, we use the same idea as in the proof of the first \hyperlink{clm:bubbles}{Claim}. Suppose that the result does not hold. Then there is a sequence $\delta_n \searrow 0$ such that for every $n$, there is a closed disk $D_n$ contained in the component of $U_{\delta_n} \setminus p(J)$ below $p(J)$. Slide the center of $D_n$ upwards until the disk first hits $B$, and let $D_n^*$ be this hitting disk. Every intersection point of $D_n^*$ with $B$ is on the upper half of $\partial D_n^*$. Let $z_n$ be any point in that intersection. Then $z_n$ is on the boundary of $B$ and the unit vector $\mathbf{v}_n$ pointing from $z_n$ towards the center of $D_n^*$ is normal to $B$. As $n \to \infty$, we have $z_n \to 0$. Moreover, the vectors $\mathbf{v}_n$ only accumulate onto vectors forming an angle at least $\frac{\pi}{2}$ with the upwards direction at $0$, hence outside $N_0$. This contradicts the semi-smoothness of $B$ at $0$. 
 \end{proof}

By this last claim, if $\delta$ is small enough then $U_\delta \cap B$ is equal to the union of $U_\delta \cap p(J)$ with the component of $U_\delta \setminus p(J)$ below $p(J)$. This neighborhood of $z$ in $B$ is clearly homeomorphic to the rectangle $(-\delta,\delta)\times (0,\delta]$, which in turn is homeomorphic to a half-disk. Thus $B$ is a $2$-manifold with boundary. 
\end{proof}

Since $\blob(x,X,Y,h)$ is semi-smooth under appropriate hypotheses, it is a $2$-manifold with boundary.

\begin{corollary} \label{cor:blob2manifold}
Assume hypothesis \myhyperlink{h4}. Then $\blob(x,X,Y,h)$ is a compact, connected, $2$-manifold with boundary for any $x \in X$.
\end{corollary}
\begin{proof}
$\blob(x,X,Y,h)$ is compact and connected because $\emb(X,Y,h)$ is (Lemma \ref{lem:compactness} and Theorem \ref{thm:connected}) and $\lift_x$ is continuous (Lemma \ref{lem:cont}). It is semi-smooth by Theorem \ref{thm:blobsemismooth}, hence a $2$-manifold with boundary by Theorem \ref{thm:aquatic}.
\end{proof}

In particular, the blob is the closure of its interior. We use this to show that the blob depends continuously on parameters. We first need to define a topology on closed subsets of a space. 

Let $S$ be a topological space and let  $\fell(S)$ be the hyperspace of closed subsets of $S$. The \emph{Fell topology} on $\fell(S)$ is the topology generated by neighborhoods of the form $N(K, \mathcal U)$, where $K\subset S$ is compact, $\mathcal U$ is a finite collection of open subsets of $S$, and $N(K, \mathcal U)$ is the set of $A \in \fell(S)$ such that $A \cap K = \varnothing$ and $A \cap U \neq \varnothing$ for every $U \in \mathcal U$. 

\begin{theorem}[Fell \cite{Fell}]
For any topological space $S$, the hyperspace $\fell(S)$ is compact. If $S$ is locally compact, then $\fell(S)$ is Hausdorff.
\end{theorem}

If $S$ is first-countable and Hausdorff, then a sequence $\{A_n\} \subset \fell(S)$ converges to $A \in \fell(S)$ if and only if
\begin{itemize}
\item for every $a \in A$, there exist $a_n \in A_n $ such that $a_n \to a$;
\item for every sequence $\{a_n\}$ with $a_n \in A_n$, if $\{a_n\}$ accumulates onto $a \in S$, then $a \in A$.
\end{itemize}

We use the Fell topology on closed subsets of $\wtilde Y \cong \DD$ where $Y$ is a finite hyperbolic surface. To prove convergence, we mostly rely on Fell's compactness theorem and the above criterion for sequences. 

It is fairly clear that the blob depends upper semi-continuously on para\-meters. The same holds for its boundary.

\begin{lemma} \label{lem:upperblob}
Suppose that $h:X \to Y$ is a generic embedding between finite Riemann surfaces. Then $\blob(z,Z,Y,h)$ and $\partial \blob(z,Z,Y,h)$  depend upper semi-continuously on $Z\setminus \{z\} \in \teich(X\setminus \{x\})$. More precisely, suppose that $Z_n \setminus \{z_n \} \to Z \setminus \{z\}$ in $\teich(X\setminus \{x\})$, that $\blob(z_n,Z_n,Y,h) \to A$  in $\fell(\wtilde Y)$, and that $\partial \blob(z_n,Z_n,Y,h) \to B$  in $\fell(\wtilde Y)$ as $n \to \infty$. Then $A \subset \blob(z,Z,Y,h)$ and $B \subset \partial \blob(z,Z,Y,h)$.
\end{lemma}
\begin{proof}
Let $c \in A$. By hypothesis there exist $c_n \in \blob(z_n,Z_n,Y,h)$ such that $c_n \to c$ as $n \to \infty$. Let $f_n \in \emb(Z_n,Y,h)$ be such that $\lift_{z_n}(f_n)=c_n$. Let $\sigma_n^\star : Z\setminus \{z\} \to Z_n \setminus \{z_n\}$ be the quasiconformal homeomorphism homotopic to the change of marking with minimal dilatation $K_n$ and let $\sigma_n: Z \to Z_n$ be its extension. By hypothesis, $K_n\to 1$ as $n\to \infty$. By Lemma \ref{lem:compactness}, we may pass to a subsequence such that $f_n \circ \sigma_n$ converges to a conformal embedding $f : Z \to Y$ homotopic to $h$. Then
$$
\lift_z(f) = \lim_{n \to \infty} \lift_z(f_n\circ \sigma_n) = \lim_{n\to \infty} \lift_{z_n}(f_n) = \lim_{n\to \infty} c_n = c,  
$$  
so that $c \in \blob(z,Z,Y,h)$.

Now let $c \in B$ and let $c_n \in \partial \blob(z_n,Z_n,Y,h)$ be such that $c_n \to c$. By Proposition \ref{prop:bdryblob}, there exists a slit mapping $f_n$ rel $z_n$ from $Z_n$ to $Y$ homotopic to $h$ such that $\lift_{z_n}(f_n)=c_n$. By Lemma \ref{lem:limqd}, we can pass to a subsequence such that $f_n$ converges to some slit mapping $f$ rel $z$ from $Z$ to $Y$. Then $\lift_z(f) = c$ so that $c \in \blob(z,Z,Y,h)$. Moreover, $y\in \partial\blob(z,Z,Y,h)$ by Corollary \ref{cor:blobcomplement}. 
\end{proof}

We do not know if the blob moves continuously in general, but it does when there are no slit mappings at the limiting parameters.

\begin{lemma} \label{lem:convifnoslit}
Assume hypothesis \myhyperlink{h4} and let $x \in X$. If $X_n \setminus \{x_n\} \to X\setminus \{x\}$ in $\teich(X\setminus \{x\})$ as $n \to \infty$, then $\blob(x_n,X_n,Y,h) \to \blob(x,X,Y,h)$ and $\partial \blob(x_n,X_n,Y,h) \to \partial \blob(x,X,Y,h)$ in $\fell(\wtilde Y)$ as $n \to \infty$.
\end{lemma}
\begin{proof}
 By compactness of $\fell(\wtilde Y)$, it suffices to prove that if $\blob(x_n,X_n,Y,h)$ converges to some closed set $A$ and $\partial \blob(x_n,X_n,Y,h)$ converges to some closed set $B$ as $n \to \infty$, then $A=\blob(x,X,Y,h)$ and $B=\partial \blob(x,X,Y,h)$. 
 
We prove convergence of the blobs first. By Lemma \ref{lem:upperblob}, the inclusion $A \subset \blob(x,X,Y,h)$ holds. We claim that the interior of $\blob(x,X,Y,h)$ is contained in $A$. Let $c$ be in the interior of $\blob(x,X,Y,h)$ and suppose that there is an infinite set $J \subset \NN$ such that $c$ is not contained in $\blob(x_n,X_n,Y,h)$ for every $n\in J$. Then for every $n \in J$, there exists a Teichm\"uller embedding $f_n$ rel $x_n$ with $\lift_{x_n}(f_n)=c$. After passing to a subsequence in $J$, we get that $f_n \to f$ for some Teichm\"uller embedding $f$ rel $x$ by Lemma \ref{lem:limqd}. We have $\lift_x(f)=c$ by continuity of $\lift_x$. By Corollary \ref{cor:blobcomplement}, $c$ is in the complement of the interior of $\blob(x,X,Y,h)$. This is a contradiction, which means that $c$ is contained in $\blob(x_n,X_n,Y,h)$ for all but finitely many indices, and hence $c \in A$. Since $A$ is closed and $\blob(x,X,Y,h)$ is the closure of its interior, we have $\blob(x,X,Y,h) \subset A$ and hence $A=\blob(x,X,Y,h)$.

By Lemma \ref{lem:upperblob}, we have $B \subset \partial \blob(x,X,Y,h)$. Let $c \in \partial \blob(x,X,Y,h)$. Let $U$ be any connected neighborhood of $c$ in $\wtilde Y$. We claim that if $n$ is large enough, then $U$ intersects both the complement of $\blob(x_n,X_n,Y,h)$ and the interior of $\blob(x_n,X_n,Y,h)$.  Suppose on the contrary that $U$ is contained in $\blob(x_n,X_n,Y,h)$ for every $n$ in an infinite set $J \subset \NN$. Then $U \subset A = \blob(x,X,Y,h)$, which is nonsense since $c$ is on the boundary of $\blob(x,X,Y,h)$. Similarly, suppose that $U$ is contained in the complement of $\blob(x_n,X_n,Y,h)$ for every $n$ in an infinite set $J \subset \NN$. Then for every $z\in U$ and every $n \in J$ there is a Teichm\"uller embedding $f_n : X_n \to Y$ rel $x_n$ homotopic to $h$ such that $\lift_{x_n}(f)=z$. By Lemma \ref{lem:limteich}, $f_n$ converges to a Teichm\"uller embedding $f$ rel $x$ after passing to a subsequence. Then $\lift_x(f)=z$ so that $z \in \partial \blob(x,X,Y,h)$ by Corollary \ref{cor:blobcomplement}. This is a contradiction, which proves the claim. Let $n$ be large enough so that $U$ intersects both the interior and the complement of $\blob(x_n,X_n,Y,h)$. Since $U$ is connected, it also intersects $\partial \blob(x_n,X_n,Y,h)$. Since $U$ can be chosen arbitrarily small, this shows that $c \in B$.
 \end{proof}

Similarly, nested families of blobs move continuously. In what follows, $\{X_r\}_{r\in [0,\infty]}$ is a $1$-parameter family of enlargements of $X$ as in Section \ref{sec:localmax}, and $\m$ is the associated modulus of extension. 

\begin{lemma}
Let $h:X \to Y$ be a generic embedding between finite Riemann surfaces, let $R$ be the maximum of $\m$ on $\emb(X,Y,h)$, and let $x \in X$. Then the maps $r \mapsto \blob(x,X_r,Y,h)$ and $r \mapsto \partial\blob(x,X_r,Y,h)$ are continuous on the interval $[0,R]$.
\end{lemma}
\begin{proof}
We may assume that $R>0$ since otherwise there is nothing to show. If $r \in [0,R)$ and $\rho \to r$, then $\blob(x,X_\rho,Y,h) \to \blob(x,X_r,Y,h)$ and $\partial\blob(x,X_\rho,Y,h) \to \partial \blob(x,X_r,Y,h)$ by Lemma \ref{lem:convifnoslit}, since $\emb(X_r,Y,h)$ does not contain any slit mapping. It remains to prove continuity at $r = R$. By compactness of the hyperspace $\fell(\wtilde Y)$ and Lemma \ref{lem:upperblob}, it suffices to show that if $r_n \nearrow R $, if $\blob(x,X_{r_n},Y,h) \to A$, and if $\partial \blob(x,X_{r_n},Y,h) \to B$ as $n \to \infty$, then $A \supset \blob(x,X_R,Y,h)$ and $B \supset \partial \blob(x,X_R,Y,h)$. By Proposition \ref{prop:globalmax}, every element of $\emb(X_R,Y,h)$ is a slit mapping\footnote{This is unless $Y$ is the triply-punctured sphere, in which case $R=\infty$ so that $\emb(X_R,Y,h)$ contains only one element anyway by Lemma \ref{lem:emptyidealbdry}.} so that $\blob(x,X_R,Y,h)$ is homeomorphic to a point or an interval (Proposition \ref{prop:blobslit}). Let $c \in \blob(x,X_R,Y,h)= \partial \blob(x,X_R,Y,h)$. Then $c\in \blob(x,X_{r_n},Y,h)$ for every $n$. Indeed, since $r_n \leq R$ there is a canonical inclusion $X_{r_n} \subset X_R$ which means that $\emb(X_R,Y,h) \subset \emb(X_{r_n},Y,h)$ and similarly for the blobs. It follows that $c \in A$. Let $U$ be a connected neighborhood of $c$ in $\wtilde Y$. Then $U$ intersects $\blob(x,X_{r_n},Y,h)$ since $c \in \blob(x,X_{r_n},Y,h)$. Thus $U$ intersects the interior of $\blob(x,X_{r_n},Y,h)$ because $\blob(x,X_{r_n},Y,h)$ is the closure of its interior. Suppose that $U$ is contained in $\blob(x,X_{r_n},Y,h)$ for every $n$ in an infinite set $J \subset \NN$. Then $U$ is contained in $A$ and hence in $\blob(x,X_R,Y,h)$. This is absurd since $\blob(x,X_R,Y,h)$ has empty interior. Thus $U$ intersects $\partial \blob(x,X_{r_n},Y,h)$ for all large enough $n$ and hence $c \in B$. 
\end{proof}

We use continuity to show that the blob has no holes and is thus homeomorphic to a closed disk, under hypothesis \myhyperlink{h4}.

\begin{proof}[Proof of Theorem \ref{thm:blobisdisk}]
We will show that the complement of $\blob(x,X,Y,h)$ is connected, which is sufficient in view of Corollary \ref{cor:blob2manifold}. Let $z_1$ and $z_2$ be any two points in  $\wtilde Y \setminus \blob(x,X,Y,h)$. Let $R$ be the maximum value of the modulus of extension $\m$.  Note that $z_1$ and $z_2$ are contained in $\wtilde Y \setminus \blob(x,X_\rho,Y,h)$ for every $\rho \in [0,R]$ as the blobs are nested. Let $r$ be the infimum of the set of $\rho \in [0,R]$ such that $z_1$ and $z_2$ are in the same component of $\wtilde Y \setminus \blob(x,X_\rho,Y,h)$. The set of such $\rho$ is non-empty since $\blob(x,X_R,Y,h)$ is a point or a compact interval, and hence has connected complement.

Suppose that $z_1$ and $z_2$ are in different components of $\wtilde Y \setminus \blob(x,X_r,Y,h)$. Then $r < R$. In particular, $\blob(x,X_r,Y,h)$ is a compact $2$-manifold so that each boundary component of $\blob(x,X_r,Y,h)$ is a simple closed curve. Let $C_1$ be the component of $\partial \blob(x,X_r,Y,h)$ surrounding $z_1$, let $C_2$ be the one surrounding $z_2$, and let $\gamma$ be a simple closed curve in the interior of $\blob(x,X_r,Y,h)$ separating $C_1$ from $C_2$. For all $\rho$ close enough to $r$ we have that $\partial \blob(x,X_\rho,Y,h)$ is disjoint from $\gamma$. On the other hand, there is a sequence $\rho_n \searrow r$ such that $z_1$ and $z_2$ belong to the same component of $\wtilde Y \setminus \blob(x,X_{\rho_n},Y,h)$. Let $\alpha_n$ be a path in $\wtilde Y \setminus \blob(x,X_{\rho_n},Y,h)$ connec\-ting $z_1$ and $z_2$. For every $n$, $\alpha_n$ intersects $\gamma$, say at $w_n$. Since $\gamma$ is compact, we may pass so a subsequence so that $w_n \to w$ for some $w \in \gamma$. Now $w$ is in the interior of $\blob(x,X_r,Y,h)$. Let $U$ be an open disk centered at $w$ whose closure is contained in the interior of $\blob(x,X_r,Y,h)$. Since $\blob(x,X_{\rho_n},Y,h) \to \blob(x,X_r,Y,h)$ as $n \to \infty$, the open set $U$ must intersect $\blob(x,X_{\rho_n},Y,h)$  for all large enough $n$. Since $w_n \in \gamma \setminus \blob(x,X_{\rho_n},Y,h) $ and since $\gamma \cup U$ is connected, the intersection of $\gamma \cup U$ with $\partial \blob(x,X_{\rho_n},Y,h)$ is non-empty. Let $\zeta_n$ be in the intersection. After passing to a subsequence, $\zeta_n$ converges to some point $\zeta$ in $\gamma \cup \overline{U}$. This is a contradiction since $\partial \blob(x,X_{\rho_n},Y,h) \to \partial\blob(x,X_r,Y,h)$ as $n \to \infty$ but $\gamma \cup \overline{U}$ is disjoint from $\partial\blob(x,X_r,Y,h)$. Therefore $z_1$ and $z_2$ belong to the same component of $\wtilde Y \setminus \blob(x,X_r,Y,h)$. 

Suppose that $r > 0$. Let $\gamma$ be a path joining $z_1$ to $z_2$ in $\wtilde Y \setminus \blob(x,X_r,Y,h)$. Since $\gamma$ is compact and $\blob(x,X_\rho,Y,h)$ depends continuously on $\rho$, the two are disjoint for all $\rho$ sufficiently close to $r$. Then $z_1$ and $z_2$ belong to the same component of $\wtilde Y \setminus \blob(x,X_\rho,Y,h)$ for all $\rho < r$ sufficiently close to $r$, which contradicts the minimality of $r$. We conclude that $r=0$ and that $z_1$ and $z_2$ belong to the same component of $\wtilde Y \setminus \blob(x,X_0,Y,h)$. Since $z_1$ and $z_2$ were arbitrary, the complement of $\blob(x,X,Y,h)=\blob(x,X_0,Y,h)$ is connected. Thus $\blob(x,X,Y,h)$ is homeomorphic to a closed disk. 
\end{proof}

\section{The deformation retraction} \label{sec:movingpoints}

The goal of this section is to prove that $\emb(X,Y,h)$ is contractible under \myhyperlink{h4}, which is the main case of Theorem \ref{thm:defret}. Recall that \myhyperlink{h4} stands for the hypothesis that $h:X \to Y$ is a generic embedding between finite Riemann surfaces (hence $X$ and $Y$ are hyperbolic), $X$ has non-empty ideal boundary, and $\emb(X,Y,h)$ is non-empty and does not contain any slit mapping.

Fix once and for all a countable dense set $\{x_1,x_2,x_3,...\}$ of distinct points in $X$ and a universal covering map $\pi_Y : \DD \to Y$. For each $n \in \NN$, let $b_n \in \DD$ be such that $\pi_Y(b_n) = h(x_n)$ and define a lift $L_n= \lift_{x_n}: \maps(X,Y,h) \to \DD$ as in Section \ref{sec:blob}. For each $f\in \maps(X,Y,h)$, $L_n(f)$ is the endpoint of the lift of $t \mapsto H(x_n,t)$ based at $b_n$, where $H$ is any homotopy from $h$ to $f$. Let $\blob(x_j) = L_j(\emb(X,Y,h))$.

Let us also fix some $F \in \emb(X,Y,h)$ which maximizes the modulus of extension $\m$ from Section \ref{sec:localmax}. Note that $\m(F) > 0$ by the hypothesis that $\emb(X,Y,h)$ does not contain any slit mapping. It follows that $L_1(F)$ does not lie on the boundary of $\blob(x_1)$ by Proposition \ref{prop:bdryblob}.  We will construct a (strong) deformation retraction of $\emb(X,Y,h)$ into $\{F\}$.

Given any $f \in \emb(X,Y,h)$, we define a sequence of paths $\gamma_n : [0,1] \to \DD$ inductively as follows. Let $ G[1] : \overline \DD \to \blob(x_1)$ be the Riemann map with $G[1](0)=L_1(F)$ and $ G[1]'(0)>0$, and let
$$
\gamma_1(t) = \begin{cases}  L_1(f) & \text{if } t \in [0,1/2) \\ G[1]((2-2t){G}[1]^{-1}(L_1(f))) & \text{if } t \in [1/2,1]. \end{cases}
$$
In words, $\gamma_1$ stays at $L_1(f)$ for half the time and then moves at constant speed along the conformal ray towards the ``center'' $L_1(F)$ of $\blob(x_1)$. In particular, $\gamma_1(t)$ belongs to $\blob(x_1)$ for every $t \in [0,1]$ so that there exists some $g \in \emb(X,Y,h)$ such that $L_1(g)=\gamma_1(t)$.

Let $n \geq 2$. Suppose that paths $\gamma_1,\ldots \gamma_{n-1}$ have been defined in such a way that
\begin{itemize}
\item the points $\pi_Y(\gamma_1(t)),\ldots,\pi_Y(\gamma_{n-1}(t))$ are distinct for every $t \in [0,1]$;
\item $\gamma_j$ is constant on the interval $[0,2^{-j}]$ for every $j \in \{1,\ldots, n-1\}$;
\item $\gamma_j(0)=L_j(f)$ and $\gamma_j(1)=L_j(F)$ for every $j \in \{1,\ldots, n-1\}$.
\end{itemize} 
Then let 
$$
X[n]=X \setminus \{x_1,\ldots, x_{n-1}\}, \quad  Y[n,t]=Y \setminus \{\pi_Y(\gamma_1(t)),\ldots,\pi_Y(\gamma_{n-1}(t))\},
$$
and let $h[n,t]= \push[n,t] \circ h$ where $\push[n,t] : Y \to Y$ is a multi-point-pushing diffeomorphism chosen so that $L_j(h[n,t])=\gamma_j(t)$ for every $j\in \{1,\ldots,n-1\}$ and every $t \in [0,1]$. The embedding $h[n,t]$ is generic for every $t$. Define
$$
E[n,t] = \emb(X[n],Y[n,t],h[n,t]) \quad \text{and} \quad \blob[n,t]= L_n(E[n,t]).
$$
We assume that $E[n,t]$ is non-empty for every $t \in [0,1]$ as part of the induction hypothesis. Note that $\blob[n,t]$ is either a closed disk or a point. Indeed, if $E[n,t]$ contains a slit mapping then it is homeomorphic to a point or an interval by Theorem \ref{thm:AlmostRigidity}. But since $h[n,t]$ has an essential puncture, there is at most one slit mapping in $E[n,t]$ by Remark \ref{rem:nopunctures2}. If $E[n,t]$ does not contain any slit mapping, then $\blob[n,t]$ is a closed disk by Theorem \ref{thm:blobisdisk}. Also, since we chose the paths $\gamma_1, \ldots, \gamma_{n-1}$ to be constant on $[0,2^{1-n}]$, the set $\blob[n,t]$ does not change for $t$ in that interval. The next step is to choose a conformal center for $\blob[n,t]$. 

\begin{lemma} \label{lem:center}
For every $t \in [0,1]$ there is is a unique map $g[n,t]$ maximizing $\m$ within $E[n,t]$. The map $t \mapsto g[n,t]$ is continuous, constant on $[0,2^{1-n}]$, and satisfies $g[n,1]=F$.
\end{lemma}
\begin{proof}
The map $\m$ is upper semi-continuous on the compact space $E[n,t]$. It thus attains its maximum at some $g[n,t]$ say with value $R$. By Proposition \ref{prop:globalmax}, the maximal extension of $g[n,t]$ is a slit mapping from $X_R\setminus \{x_1,\ldots, x_{n-1}\}$ to $Y[n,t]$. By Remark \ref{rem:nopunctures2}, the map $g[n,t]$ is unique since it sends a puncture to a puncture. Any limit of the maximal extension of $g[n,t]$ as $t \to s$ is a slit mapping from $X_R\setminus \{x_1,\ldots, x_{n-1}\}$ to $Y[n,s]$ by Lemma \ref{lem:limteich} and thus its restriction to $X[n]$ maximizes $\m$ in $E[n,s]$ by Lemma \ref{lem:localisglobal}. Thus $g[n,t] \to g[n,s]$ as $t \to s$. The paths $\gamma_j$ for $j \in \{1,\ldots , n-1\}$ are all constant on $[0,2^{1-n}]$ so $g[n,t]$ does not change on that interval. Finally, since $F$ maximizes $\m$ on $\emb(X,Y,h)$, it maximizes $\m$ on the subset $E[n,1]$ as well, and we have $g[n,1]=F$.   
\end{proof}

Let $ G[n] : \overline \DD \to \blob[n,0]$ be the Riemann map normalized in such a way that $ G[n](0)=L_n(g[n,0])$ and $ G[n]'(0)>0$. Then let
$$
\gamma_n(t) = \begin{cases} L_n(f) & \text{if } t \in [0,2^{-n}) \\
  G[n]((2-2^nt){G}[n]^{-1}(L_n(f))) & \text{if } t \in [2^{-n},2^{1-n}) \\
  L_n(g[n,t]) & \text{if } t \in [2^{1-n},1]. \end{cases}
$$
This means that $\gamma_n$ stays put at $L_n(f)$ for some time, then travels along the conformal ray towards the center of $\blob[n,0]=\blob[n,2^{1-n}]$, and then follows the center for the rest of the time. It is possible that $\blob[n,0]$ is a point if $L_{n-1}(f)$ is in the boundary of $\blob[n-1,0]$. In that case we let $ G[n]: \overline \DD \to \blob[n,0]$ be the constant map, i.e., we keep $\gamma_n$ constant on $[0,2^{1-n}]$. By construction we have $\pi_Y(\gamma_n(t))\in Y[n,t]$ which means that the points $\pi_Y(\gamma_1(t)),\ldots,\pi_Y(\gamma_{n-1}(t)),\pi_Y(\gamma_{n}(t))$ are distinct for every $t\in [0,1]$. Moreover the path $\gamma_n$ is constant on the interval $[0,2^{-n}]$. Finally, $E[n+1,t]$ is non-empty since $\gamma_n(t) \in \blob[n,t]$ for every $t \in [0,1]$. This finishes the induction scheme.

We now show that the paths $\{\gamma_n\}$ automatically define a path from $f$ to $F$ inside the space $\emb(X,Y,h)$.

\begin{lemma} \label{lem:uniquepath}
For every $t \in [0,1]$, there exists a unique $f_t\in\emb(X,Y,h)$ such that $L_n(f_t)=\gamma_n(t)$ for every $n \in \NN$. The map $t \mapsto f_t$ is continuous and satisfies $f_0 = f$ and $f_1=F$.
\end{lemma}
\begin{proof}
Observe that $E[n,t]$ is a non-empty closed subset of $\emb(X,Y,h)$ and is thus compact. Therefore, for each $t \in [0,1]$, the nested intersection $\bigcap_{n=1}^\infty E[n,t]$ is non-empty. Any two functions in the intersection agree on the dense set $\{x_1,x_2,...\}$ and hence on all of $X$. Therefore, there is a unique function $f_t$ in the intersection. Moreover, $f_t$ varies continuously with $t$. Indeed, if $g$ is any limit of any subsequence of $f_t$ as $t \to s$, then for every $n \in \NN$ we have
$$
g(x_n)= \lim_{t \to s} f_t(x_n) = \lim_{t \to s} \pi_Y( \gamma_n(t)) = \pi_Y (\gamma_n(s)) = f_s(x_n)
$$ 
so that $g=f_s$. It follows that $f_t \to f_s$ as $t\to s$. By construction we have $L_n(f)=\gamma_n(0)$ and $L_n(F)=\gamma_n(1)$ for every $n \in \NN$ so that $f_0=f$ and $f_1=F$. 
\end{proof}

We thus have a map $H:\emb(X,Y,h) \times [0,1] \to \emb(X,Y,h)$ defined by $H(f,t)=f_t$. This map is such that 
\begin{itemize}
\item $t \mapsto H(f,t)$ is continuous for every $f \in \emb(X,Y,h)$;
\item $H(f,0)=f$ and $H(f,1)=F$ for every $f \in \emb(X,Y,h)$;
\item $H(F,t)=F$ for every $t \in [0,1]$.
\end{itemize}
The last point holds because if $f=F$, then every path $\gamma_n$ is constant and hence $f_t=F$ for every $t$. It remains to prove that $H$ is continuous in both variables.

\begin{lemma} \label{lem:Hcontinuous}
Suppose that for every $n \in \NN$, the path $\gamma_n \in \maps([0,1],\DD)$ depends continuously on $f \in\emb(X,Y,h)$, where each space is equipped with the compact-open topology. Then the map $H$ defined above is continuous.
\end{lemma}
\begin{proof}
If for every $n \in \NN$ the map $(f,t) \mapsto H(f,t)(x_n)$ is continuous, then $H$ is continuous. This is because of the compactness of $\emb(X,Y,h)$ and the fact that $\{x_1,x_2,\ldots\}$ is dense in $X$ (the proof is the same as in Lemma \ref{lem:uniquepath}). Since $H(f,t)(x_n)=f_t(x_n)=\pi_Y(\gamma_n(t))$, it thus suffices that $(f,t)\mapsto \gamma_n(t)$ be continuous. This condition is equivalent to the requirement that $f \mapsto \gamma_n$ is continuous since $[0,1]$ is locally compact Hausdorff (see \cite[p.287]{Munkres}). 
\end{proof}

Since the map $f \mapsto L_1(f)$ is continuous and the Riemann map ${G}[1]$ is continuous, it is easy to see that $f \mapsto \gamma_1$ is continuous. We proceed by induction for the rest. Let $n \geq 2$ and suppose that the maps $f \mapsto \gamma_{j}$ are all continuous for $j = 1,...,n-1$. We will prove that the map $f \mapsto {G}[n]$ is continuous, which obviously implies that $f \mapsto \gamma_n$ is continuous. We use the following theorem of Rad\'o, a proof of which is given in \cite[p.26]{Pommerenke}.

\begin{theorem}[Rad\'o] \label{thm:rado}
Let $(D_k,w_k)$ and $(D,w)$ be topological closed disks in $\CC$, each with a marked point in the interior. Suppose that $w_k \to w$ and that $D_k \to D$ in the Fell topology. Suppose also that there are parametrizations  $c_k : S^1 \to \partial D_k$ and $ c : S^1 \to \partial D$ such that $c_k \to c$ uniformly. Then the normalized Riemann map $(\overline \DD,0) \to (D_k,w_k)$ converges uniformly on $\overline \DD$ to the normalized Riemann map $(\overline \DD,0) \to (D,w)$.  
\end{theorem}

By a slight generalization\footnote{The codomain $Y[n,t]$ is not fixed but depends continuously on $(f,t)$. The results of Section \ref{sec:noholes} generalize easily to this situation.} of Lemma \ref{lem:upperblob}, the maps $(f,t) \mapsto \blob[n,t]$ and $(f,t) \mapsto \partial \blob[n,t]$ are upper semi-continuous in $t$. Moreover, they are continuous at every $(f,t)$ for which $E[n,t]$ does not contain any slit mapping by Lemma \ref{lem:convifnoslit}. But if $E[n,t]$ contains a slit mapping, then $\blob[n,t]=\partial \blob[n,t]$ is a single point and thus upper semi-continuity at $(f,t)$ implies continuity. By Lemma \ref{lem:center}, the conformal center $L_n(g[n,t])$ of $\blob[n,t]$ also depends continuously on $(f,t)$. The only thing that remains to be checked is that the boundary $\partial \blob[n,t]$ can be parametrized as to converge uniformly.

\begin{definition}
Let $\{c_k\}_{k=1}^\infty$ be a sequence of simple closed curves in $\CHAT$. We say that $\{c_k\}_{k=1}^\infty$ has a \emph{collapsing finger} if after passing to a subsequence, there exist $x_k,y_k,z_k, w_k \in S^1$ in cyclic order and $x,y \in \CHAT$ with $x\neq y$  such that $c_k(x_k)\to x$, $c_k(y_k) \to y$, $c_k(z_k) \to x$, and $c_k(w_k)\to y$.   
\end{definition}

\begin{figure}[htp]  \centering
\includegraphics[scale=.95]{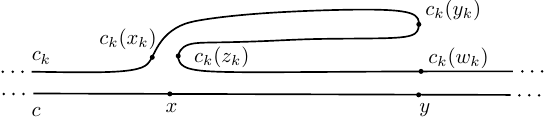}
\caption{A piece of a curve with a finger about to collapse.}
\label{fig:finger}
\end{figure}

We now show that collapsing fingers are the only obstructions to uniform convergence of simple closed curves.

\begin{lemma} \label{lem:nofingersuniform}
Let $c_k$ and $c$ be simple closed curves in $\CHAT$ such that $c_k(S^1)$ converges to  $c(S^1)$ in the Fell topology. If $\{c_k\}_{k=1}^\infty$ does not have any collapsing finger, then we can reparametrize $c_k$ such that $c_k \to c$ uniformly.
\end{lemma}

\begin{proof}
By the Jordan--Schoenflies Theorem, $c$ can be extended to a homeo\-morphism $\hat{c} : \CHAT \to \CHAT$. Then $\hat{c}^{-1}\circ c_k(S^1) \to S^1$ in the Fell topology and the sequence $\{\hat{c}^{-1}\circ  c_k\}_{k=1}^\infty$ does not have any collapsing finger. Moreover, if $\sigma_k:S^1 \to S^1$ is a homeomorphism such that $\hat{c}^{-1}\circ  c_k \circ \sigma_k$ converges uniformly to the inclusion map $S^1 \hookrightarrow \CC$, then $c_k\circ \sigma_k$ converges uniformly to $c$. We may thus assume that $c$ is the inclusion map $S^1 \hookrightarrow \CC$. 

If $k$ is large enough, then $c_k(S^1)$ is disjoint from $0$ and $\infty$. We claim that if $k$ is large enough, then the winding number of $c_k$ around the origin is $\pm 1$. Since $c_k$ is simple, its winding number is either $-1$, $0$, or $1$. Suppose the claim is false. Then after passing to a subsequence, the winding number of $c_k$ is $0$ for every $k$. Let $\arg(c_k)=c_k / |c_k|$ and let $\Arg(c_k) : S^1 \to \RR$ be a lift of $\arg(c_k)$ under the universal covering map $\RR \to S^1$. This lift exists because the winding number is zero. Let $[a_k,b_k]$ be the image of $\Arg(c_k)$. Since $c_k(S^1)$ converges to $S^1$, it follows that the image $\arg(c_k)(S^1)$ converges to $S^1$ as well, and hence $\liminf_{n\to \infty} b_k - a_k \geq 2\pi$. Thus if $k$ is large enough, then $b_k-a_k> \pi$ . Let $x_k$ and $z_k$ in $S^1$ be such that $\Arg(c_k)(x_k)=a_k$ and $\Arg(c_k)(z_k)=\min(b_k,a_k+2 \pi)$. Also let $y_k\in\overline{x_k z_k}$ and $w_k \in \overline{z_k x_k}$ be such that $$\Arg(c_k)(y_k)=\Arg(c_k)(w_k)=a_k+\pi.$$ Since $c_k(S^1) \to S^1$, we may pass to a subsequence so that $c_k(x_k)$, $c_k(z_k)$, $c_k(y_k)$, and $c_k(w_k)$ converge to some $x$, $y$, $z$, and $w$ in $S^1$. Then $x=z$, $y=w$, and $x \neq y$, i.e. $\{c_k\}_{n=1}^\infty$ has a collapsing finger. This is a contradiction, which proves the claim.

If the winding number of $c_k$ around the origin is $-1$, then we reverse the parametrization so that it becomes $+1$. Let $\zeta_1^k < \zeta_2^k < \ldots < \zeta_k^k $ be a partition of $S^1$ into $k$ congruent arcs. Since $c_k$ has winding number $1$, we can find $\xi_1^k < \xi_2^k < \ldots < \xi_k^k $ in $S^1$ such that $\arg(c_k)(\xi_j^k)=\zeta_j^k$ for every $j \in \{1,\ldots,k\}$. Let $\sigma_k: S^1 \to S^1$ be any homeomorphism such that $\sigma_k(\zeta_j^k)=\xi_j^k$ for every $j \in \{1,\ldots,k\}$. We claim that $c_k\circ \sigma_k$ converges uniformly to the inclusion map $c: S^1 \hookrightarrow \CC$.

To simplify notation, we assume that $c_k$ was parametrized correctly from the start, i.e. we assume that for every $k\gg 0$ and every $j \in \{1,\ldots,k\}$, we have $\arg(c_k)(\zeta_j^k)=\zeta_j^k$.
If $c_k$ does not converge uniformly to $c$, then there exists an $\eps >0$ and an infinite set $J \subset \NN$ such that for every $k \in J$, there exists a $y_k \in S^1$ such that $|c_k(y_k)-y_k| \geq \eps$. Since $S^1$ is compact and $c_k(S^1) \to S^1$, we can pass to a subsequence such that $y_k \to x$ and $c_k(y_k) \to y$ for some $x$ and $y$ in $S^1$. Note that $|y - x| \geq \eps$ and in particular $y \neq x$.  Let $j \in \{1,\ldots, k\}$ be such that $\zeta_j^k \leq y_k < \zeta_{j+1}^k$, where we define $\zeta_{k+1}^k=\zeta_1^k$. Then let $x_k=\zeta_j^k$ and $z_k=\zeta_{j+1}^k$. Also let $w_k \in \{\zeta_1^k,\ldots,\zeta_k^k\}$ be the closest point to $y$ which comes after $z_k$ but before $x_k$ in the cyclic order on $S^1$. We have $c_k(x_k)=x_k \to x$, $c_k(y_k)\to y$, $c_k(z_k)=z_k\to x$, and $c_k(w_k)=w_k \to y$. In other words, the sequence $\{c_k\}_{k=1}^\infty$ has a collapsing finger, which is a contradiction. Therefore, $c_k$ converges uniformly to $c$.  
\end{proof}

To conclude the proof that $H$ is continuous, we show that $\partial \blob[n,t]$ does not have any collapsing fingers. The reason for this is that the blobs $\blob[n,t]$ are uniformly semi-smooth, meaning that any non-zero limit of a sequence of vectors normal to some blob is normal to the limiting blob\footnote{The proof of this is a straightforward generalization of Theorem \ref{thm:blobsemismooth}}. Now if there was a collapsing finger somewhere, then we would see two normal vectors pointing opposite to each other in the limit, which is forbidden by the definition of semi-smoothness.

\begin{lemma}
Suppose that $(f_k , t_k) \to (f,t)$ in $\emb(X,Y,h)\times [0,1]$. Then $\partial \blob[n,t_k]$ converges to $\partial \blob[n,t]$  without collapsing fingers.  
\end{lemma}
\begin{proof}
Let $B_k= \blob[n,t_k]$, $B= \blob[n,t]$, $c_k=\partial B_k$, and $c=\partial B$. Suppose that after passing to a subsequence we can find $x_k,y_k,z_k,w_k \in c_k$ in cyclic order and $x,y \in c$ with $x\neq y$ such that $x_k,z_k \to x$ and $y_k,w_k\to y$. Rotate and translate the picture in such a way that $x=0$ and that the upward direction $i$ bisects the cone $N_0$ of vectors normal to $B$ at $0$. 

By the proof of Theorem \ref{thm:aquatic}, there exists a rectangle $Q$ centered at $0$ with sides parallel to the coordinate axes such that $Q \cap c$ is the graph of a continuous function.  Since $y \neq 0$, we can  shrink $Q$ so that it does not contain $y$. Let $\delta>0$ be such that the vertical $\delta$ neighborhood $U_\delta$ of $Q \cap c$ is contained in $Q$. Then $Q\setminus U_\delta$ is compact and disjoint from $c$. Let $k$ be large enough so that $x_k$ and $z_k$ are in $Q$,  $y_k$ and $w_k$ are not in $Q$, and $c_k$ is disjoint from $Q \setminus U_\delta$. Then of the three subarcs $\overline{x_k y_k}$, $\overline{y_k z_k}$, and $\overline{z_k w_k}$ of $c_k$, at least two must cross the same vertical side $S$ of $Q$. This implies that $S \setminus B_k$ is disconnected. Hence there is an open subinterval $I$ of $S \setminus B_k$ whose highest point is contained in $B_k$. Let $D$ be a closed round disk centered on $I$ and contained in $\CC \setminus B_k$. Move the center of $D$ upwards until the boundary of the translated disk $D^*$ first hits $B_k$. Any intersection point $p_k$ of $D^*$ with $B_k$ is on the top half of $\partial D^*$. Moreover, the unit vector $\mathbf{v}_k$ based at $p_k$ and pointing in the direction of the center of $D^*$ is normal to $B_k$. 

\begin{figure}[htp]  \centering
\includegraphics[scale=1]{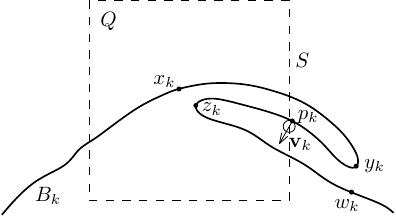}
\caption{If the sequence of blobs has a collapsing finger, then we can find a sequence of normal vectors which accumulate to a vector which is not normal to the limiting blob.}
\label{fig:collapse}
\end{figure}

Since we can choose $Q$ to be arbitrarily small, we can arrange so that $p_k \to 0$. Then the normal vectors $\mathbf{v}_k$ accumulate onto vectors pointing towards the lower half-plane at $0$. This is a contradiction since the cone of normal vectors $N_0$ is contained in the upper half-plane. 
\end{proof}

Thus by Lemma \ref{lem:nofingersuniform}, the boundary of $\blob[n,t]$ can be parametrized in a way that depends uniformly continuously on $(f,t)$. By Theorem \ref{thm:rado}, this implies that the Riemann map $G[n] : \overline \DD \to \blob[n,t]$ depends uniformly continuously on $(f,t)$. Therefore the path $\gamma_n$ depends uniformly continuously on $f$, and hence $H$ is continuous by Lemma \ref{lem:Hcontinuous}. This shows that $\emb(X,Y,h)$ is contractible under hypothesis \myhyperlink{h4}. 

More generally, $\emb(X,Y,h)$ is contractible whenever $h$ is generic and $\emb(X,Y,h)$ is non-empty. The proof for $\partial X = \varnothing$ is in Subsection \ref{subsec:aut} and if $\emb(X,Y,h)$ contains a slit mapping then this follows from Theorem \ref{thm:AlmostRigidity}.

\section{The remaining cases} \label{sec:trivialcases}

In this section, we complete the proof of Theorem \ref{thm:defret}, that is, we describe the homotopy type of $\emb(X,Y,h)$ when $h$ is not generic, always assuming that $\emb(X,Y,h)$ is non-empty. As before, we let $\{x_1,x_2,x_3,\ldots\}$ be a dense set of distinct points in $X$. 

\subsection{$h$ is cyclic but not parabolic}

Suppose that $h: X \to Y$ is cyclic but not parabolic and that $Y$ is hyperbolic. This case is analogous to the case where $h$ is generic. We only provide an outline of the proof. 

If $\emb(X,Y,h)$ contains a slit mapping, then it is homeomorphic to a circle by Theorem \ref{thm:AlmostRigidity}. So we may assume that $\emb(X,Y,h)$ does not contain any slit mapping. Form the annulus cover $\pi_A: A \to Y$ corresponding to the image of $\pi_1(h)$. We can define a continuous map $\lift_{x_1}: \maps(X,Y,h) \to A$ in a similar fashion as in Section \ref{sec:blob}. Let $\blob(x_1,X,Y,h)$  be the image of $\emb(X,Y,h)$ by $\lift_{x_1}$. This $\blob(x_1,X,Y,h)$ is compact and connected since $\emb(X,Y,h)$ is, according to Lemma \ref{lem:compactness} and Theorem \ref{thm:connected}. All of the results from Sections \ref{sec:blob}, \ref{sec:semismooth} and \ref{sec:noholes} pertai\-ning to the local geometry of $\blob(x_1,X,Y,h)$ extend to the current setting. The conclusion is that $\blob(x_1,X,Y,h)$ is a $2$-dimensional manifold with boundary.

 Let $\{X_r\}_{r \in [0,\infty]}$ be a $1$-parameter family of enlargements of $X$ as in Section \ref{sec:localmax}, let $\m$ be the associated modulus of extension, and let $R$ be the ma\-xi\-mum value of $\m$. By Proposition \ref{prop:globalmax} and Theorem \ref{thm:AlmostRigidity}, $\emb(X_R,Y,h)$ is homeomorphic to a circle via the evaluation map. It follows that $\blob(x_1,X_R,Y,h)$ is homeo\-morphic to a circle so that its complement in $A$ has two connected components. By a similar argument as in the proof of Theorem \ref{thm:blobisdisk}, the complement of $\blob(x_1,X,Y,h)$ has two connected components as well. Since $\blob(x_1,X,Y,h)$ is a planar $2$-dimensional manifold with two boundary components, it is homeomorphic to a closed annulus.

Let $D$ be any deformation retraction of $\blob(x_1,X,Y,h)$ into the circle $\blob(x_1,X_R,Y,h).$ The strategy after this step is the same as in Section \ref{sec:movingpoints}. Given $f \in \emb(X,Y,h)$, we let $\gamma_1(t)= \lift_{x_1}(f)$ for $t \in [0,1/2)$ and $\gamma_1(t) = D(\lift_{x_1}(f),2t-1)$ for $t \in [1/2,1]$. The map $h[1,t]: X \setminus \{x_1\} \to Y \setminus \pi_A(\gamma_1(t))$ obtained by composing $h$ with point-pushing along $\pi_A \circ \gamma_1$ is now generic so that we can construct the next paths $\gamma_2, \gamma_3, \ldots$ in the same way as in the previous section. The end result is a deformation retraction of $\emb(X,Y,h)$ into the circle $\emb(X_R,Y,h)$.

\subsection{$Y$ is the punctured disk} \label{subsec:puncdisk}

Let $h: X \to \DD \setminus \{0\}$ be a non-trivial (hence parabolic) embedding. By Montel's theorem, any sequence in $\emb(X, \DD \setminus \{0\},h)$ has a subsequence converging locally uniformly to either an element in $\emb(X, \DD \setminus \{0\},h)$ or to a constant map into $\overline \DD$. As $h$ is non-trivial, the only possible constant limit is $0$. Thus the set of $f \in \emb(X, \DD \setminus \{0\},h)$ which map $x_1$ outside of a fixed open neighborhood of $0$ is compact. Since $\overline \DD \setminus \{0\}$ acts by multiplication on $\emb(X, \DD \setminus \{0\},h)$, the image $V(x_1)$ of the evaluation map at $x_1$ is equal to a punctured disk $r \overline \DD \setminus \{0\}$ for some $r\in(0,1)$. By a slight modification of Proposition \ref{prop:bdryblob}, for every $y_1 \in \partial V(x_1)$ there is a unique $f \in \emb(X, \DD \setminus \{0\}, h)$ such that $f(x_1)=y_1$. We do not need to lift here: if two maps $f,g \in \maps(X,\DD \setminus \{0\},h)$ agree at $x_1$ they are homotopic rel $x_1$ because the pure mapping class group of the twice-punctured disk is trivial. Thus the inverse image of $\partial V(x_1)= \{ z \in \DD : |z|=r \}$ in $\emb(X, \DD \setminus \{0\},h)$ is homeomorphic to a circle.

Given any $f \in \emb(X, \DD \setminus \{0\}, h)$, let $\gamma_1:[0,1] \to \DD \setminus \{0\}$ be cons\-tant equal to $f(x_1)$ on $[0,1/2)$ followed by the radial ray from $f(x_1)$ to $rf(x_1)/|f(x_1)|$ on $[1/2,1]$. Note that the map $h[1,t]: X \setminus \{x_1\} \to \DD \setminus \{0,\gamma_1(t)\}$ (notation as in Section \ref{sec:movingpoints}) is generic for every $t \in [0,1]$, so we can construct the next paths $\gamma_2, \gamma_3, \ldots$ as before. This gives a deformation retraction of $\emb(X , \DD \setminus \{0\},h)$ into the circle of slit mappings rel $x_1$.

\subsection{$Y$ is the disk}

Suppose that $Y$ is the unit disk $\DD$. Consider the map $\DD \to \aut_0(\DD)$ which sends $a\in \DD$ to the automorphism $M_a(z)=\frac{z-a}{1-\bar{a}z}$. We have a homeomorphism
$$
\emb(X,\DD,h)\to \DD \times \emb(X \setminus \{ x_1 \} , \DD \setminus  \{0\},M_{h(x_1)} \circ h^\star)
$$
given by $f \mapsto (f(x_1),M_{f(x_1)} \circ f)$, where $h^\star:X \setminus \{ x_1 \}\to \DD \setminus \{h(x_1)\}$ is the restriction of $h$. Since $\DD$ is contractible, 
$\emb(X,\DD,h)$ is homotopy equivalent to the second factor. Note that $h$ is trivial and $M_{h(x_1)} \circ h^\star$ is parabolic, with codomain the once-punctured disk. By the previous case, $\emb(X \setminus \{ x_1 \} , \DD \setminus  \{0\},M_{h(x_1)} \circ h^\star)$ is homotopy equivalent to $S^1$. Thus $\emb(X,\DD,h)$ is homotopy equivalent to a circle, which in turn is homotopy equivalent to the unit tangent bundle of $\DD$.

\subsection{$X$ is the disk} \label{subsec:disk}

Suppose that $X=\DD$ and that $Y\neq \DD$ is hyperbolic. Here $h$ is trivial so we may drop it from the notation. We first define a map from the unit tangent bundle $T^1 Y$ to $\emb(\DD, Y)$ as follows. Given $\mathbf{v} \in T_y ^1Y$, let $D_\mathbf{v} \subset Y$ be the largest embedded ball centered at $y$ in the hyperbolic metric, and let $F_\mathbf{v} : \DD \to D_\mathbf{v}$ be the Riemann map with $F_\mathbf{v}(0)=y$ and $F_\mathbf{v}'(0)=\lambda \mathbf{v}$ for some $\lambda > 0$. The map $\mathbf{v} \mapsto F_\mathbf{v}$ is an embedding from $T^1 Y$ to $\emb(\DD, Y)$. We will construct a deformation retraction of $\emb(\DD, Y)$ into the image of that map.

Given $f \in \emb(\DD, Y)$, let $\mathbf{v} \in T^1 Y$ be the unique vector such that $f'(0)= \lambda \mathbf{v}$ for some $\lambda > 0$. Then let $r \in (0,1]$ be the largest number such that $f(r\DD) \subset D_\mathbf{v}$ and let $f^\dagger(z)=f(rz)$. Then $F_\mathbf{v}^{-1} \circ f^\dagger : \DD \to \DD$ is a conformal embedding which fixes the origin and has positive derivative there.

Let $g:\DD \to \DD$ be a conformal embedding with $g(0)=0$ and $g'(0)>0$. For every $t \in (0,1]$, define $\rho_t = \inf \{ \rho>0 : g(t\DD) \subset \rho\DD \}$ and $g_t(z)= g(tz) / \rho_t$. By Koebe's distortion theorem \cite[p.33]{Duren} we have 
$$
\frac{t}{(1+t)^2} \leq \frac{\rho_t}{|g'(0)|} \leq \frac{t}{(1-t)^2}
$$
and it follows that $g_t \to \id$ as $t \to 0$.

We define a deformation retraction of $\emb(\DD, Y)$ into $\{F_{\vec{v}} : \vec{v} \in T^1 Y \}$ by the formula
$$
H(f,t) = \begin{cases}  z \mapsto f((1-(1-r)2t) z) & \text{if } t \in [0,1/2) \\ F_\mathbf{v} \circ g_{(2-2t)} & \text{if } t \in [1/2,1] \end{cases}
$$
where $r$, $f^\dagger$ and $\vec{v}$ are defined in terms of $f$ as above and $g=F_\mathbf{v}^{-1} \circ f^\dagger$.

\subsection{$h$ is trivial}

Suppose that $X \neq \DD$, that $Y\neq \DD$ is hyperbolic, and that $h$ is trivial. Given $f \in \emb(X,Y,h)$, let $D_f \subset Y$ be the smallest topological disk containing the image of $f$. We can define $D_f$ by filling the holes of $f(X)$. Then let $F : \DD \to D_f$ be the Riemann map with $F(0)=f(x_1)$ and $F^{-1}(f(x_2))>0$. We thus obtain an embedding 
$$
 \emb(X,Y,h) \to \emb(\DD , Y) \times W
$$
defined by $f\mapsto(F, F^{-1} \circ f)$, where $W$ is the set of all conformal embeddings from $X$ to $\DD$ sending $x_1$ to $0$ and $x_2$ to a positive real number. There is an obvious left inverse 
$$ \emb(\DD , Y) \times W \to \emb(X,Y,h)$$
given by $(G,g)\mapsto G \circ g$. 

By the previous subsection, there is a deformation retraction $H_1$ from $\emb(\DD , Y)$ into a subset homeomorphic to $T_1 Y$. As for $W$, it is homeomorphic to the quotient $\emb(X\setminus\{x_1\},\DD\setminus\{0\},g^\star)/S^1$ where $g:X \to \DD$ is any embedding with $g(x_1)=0$, $g^\star$ is the restriction of $g$, and $S^1$ acts by multiplication. Subsection \ref{subsec:puncdisk} provides a deformation retraction $H_2$ of $W$ into a singleton. Therefore $\emb(X,Y,h)$ deformation retracts into a subset homeomorphic to $T^1 Y$ via the formula $(f,t)\mapsto H_1(F,t)\circ H_2(F^{-1}\circ f,t)$.

\subsection{$h$ is parabolic} \label{subsec:para}

Suppose that $h:X \to Y$ is parabolic, where $Y$ is hyperbolic and not the once-punctured disk. Let $p$ be the puncture around which $h$ wraps non-trivially. Given $f \in \emb(X,Y,h)$, we define a disk $D_f \subset Y \cup \{p\}$ by filling the holes of $f(X)$ in $Y \cup \{p\}$. Then we define $F : \DD \to D_f$ to be the Riemann map with $F(0)=p$ and $F^{-1}(f(x_1))>0$. This yields an embedding 
$$
\emb(X,Y,h) \to \emb(\DD\setminus \{0\}, Y, G^\star) \times W
$$
defined by $f \mapsto (F , F^{-1} \circ f)$, where $G^\star$ is the restriction of some embedding $G : \DD \to Y \cup \{p\}$ satisfying $G(0)=p$ and $W$ is the set of all conformal embeddings $g:X \to \DD \setminus \{0\}$ such that $g(x_1)>0$ which are homotopic to $F_0^{-1}\circ f_0$ for any $f_0\in \emb(X,Y,h)$ . The first factor deformation retracts into a circle by Subsection \ref{subsec:disk} whereas the second factor deformation retracts into a point by Subsection \ref{subsec:puncdisk}. By applying the left inverse of the above embedding (the composition map), we get a deformation retraction of $\emb(X,Y,h)$ into a circle.\medskip

It remains to treat the cases where $Y$ is not hyperbolic. In those cases, we can quotient $\emb(X,Y,h)$ by the action of $\aut_0(Y)$ to reduce to the hyperbolic case.

\subsection{$Y$ is a torus}

Suppose that $Y$ is a torus. Then $\emb(X,Y,h)$ is homeomorphic to $$\aut_0(Y) \times \emb(X\setminus \{x_1\},Y\setminus \{h(x_1)\},h^\star)$$ where $h^\star:X\setminus \{x_1\} \to Y\setminus \{h(x_1)\}$ is the restriction of $h$. Recall that $\aut_0(Y)$ is homeomorphic to $Y$ itself.

If $h$ is trivial, then $h^\star$ is parabolic and its codomain is a hyperbolic surface. Subsection \ref{subsec:para} shows that $\emb(X\setminus \{x_1\},Y\setminus \{h(x_1)\},h^\star)$ is then homotopy equivalent to $S^1$. This means that $\emb(X,Y,h)$ is homotopy equivalent to a $3$-dimensional torus, or the unit tangent bundle of $Y$.

If $h$ is non-trivial, then $h^\star$ is generic so that $\emb(X\setminus \{x_1\},Y\setminus \{h(x_1)\},h^\star)$ is contractible. Thus $\emb(X,Y,h)$ is homotopy equivalent to $\aut_0(Y)\approx Y$.

\subsection{$Y$ is the sphere with at most $2$ punctures}

Suppose that $Y$ is the Riemann sphere $\CHAT$. Let $P=\{x_1,x_2,x_3\}$. Then $$\emb(X,Y,h) \approx \aut_0(Y) \times \emb(X\setminus P,Y\setminus h(P),h^\star)$$
where $h^\star:X \setminus P \to Y \setminus h(P)$ is the restriction of $h$. Observe that $h$ is automa\-tically trivial and that $h^\star$ is generic. Thus $\emb(X,Y,h)$ is homotopy equivalent to $\aut_0(Y)$. We stated in Subsection \ref{subsec:aut} that $\aut_0(\CHAT)$ is homotopy equivalent to the unit tangent bundle of $\CHAT$. Here is a proof. First, $\aut_0(\CHAT)$ is homeomorphic to the set of triples $(a,\mathbf{v},b)$ where $a,b \in \CHAT$ are distinct and $\mathbf{v} \in T_a \CHAT$ is non-zero. The homeomorphism is given by $f \mapsto (f(0),f'(0),f(\infty))$. This set of triples deformation retracts into $T\CHAT \setminus \vec{0}$ (the complement of the zero section in the tangent bundle of $\CHAT$) by moving the point $b$ along the spherical geodesic to the antipode of $a$. Now $T\CHAT \setminus \vec{0}$ clearly deformation retracts into the unit tangent bundle $T^1\CHAT$.

Suppose that $Y$ is the Riemann sphere minus a point. Then we can repeat the same trick with $P=\{x_1,x_2 \}$ instead. Again, $h$ is trivial and its restriction $h^\star:X \setminus P \to Y \setminus h(P)$ is generic. Thus $\emb(X,Y,h)$ is homotopy equivalent to $\aut_0(\CC)$. The latter is homeomorphic to the complement of the zero section in $T\CC$. This deformation retracts onto the unit tangent bundle of $\CC$ (which deformation retracts further into a circle).

Suppose that $Y$ is the Riemann sphere minus two points. Then we can puncture at one point to factor out the action of $\aut_0(Y)$. That group is homeo\-morphic to $S^1 \times \RR$, hence homotopy equivalent to $S^1$. If $h$ is trivial, then its restriction $h^\star : X \setminus \{ x_1 \} \to Y \setminus \{h(x_1)\}$ is a parabolic embedding into a hyperbolic surface and $\emb(X\setminus \{x_1\},Y\setminus \{h(x_1)\},h^\star)$ is homotopy equivalent to $S^1$ by Subsection \ref{subsec:para}. Thus $\emb(X,Y,h)$ is homotopy equi\-va\-lent to $S^1 \times S^1$, which is in turn homotopy equivalent to the unit tangent bundle of $Y$. If $h$ is non-trivial, then it is cyclic. In this case, $h^\star$ is generic so that $\emb(X\setminus \{x_1\},Y\setminus \{h(x_1)\},h^\star)$ is contractible and $\emb(X,Y,h)$ is homotopy equivalent to $S^1$.\medskip

The reader can check that we have exhausted all possibilities for the embedding $h:X\to Y$, which concludes the proof of Theorem \ref{thm:defret}. The latter obviously implies Theorem \ref{thm:CouchThm}.

\providecommand{\bysame}{\leavevmode\hbox to3em{\hrulefill}\thinspace}
\providecommand{\MR}{\relax\ifhmode\unskip\space\fi MR }
\providecommand{\MRhref}[2]{%
  \href{http://www.ams.org/mathscinet-getitem?mr=#1}{#2}
}
\providecommand{\href}[2]{#2}

\end{document}